# The Tits alternative for $\mathrm{Out}(F_n)$ I: Dynamics of exponentially-growing automorphisms

By Mladen Bestvina, Mark Feighn, and Michael Handel

## Contents

## 1. Introduction

A group satisfies the *Tits alternative* if each of its subgroups either contains a free group of rank two or is virtually solvable. The Tits alternative derives its name from the result of J. Tits [Tit] that finitely generated linear groups satisfy this alternative. N. Ivanov [Iva] and J. McCarthy [McC] have shown that mapping class groups of compact surfaces also satisfy this alternative. J. Birman, A. Lubotzky, and J. McCarthy [BLM] and N. Ivanov [Iva] complement the Tits alternative for surface mapping class groups by showing that solvable subgroups of such are virtually finitely generated free abelian of bounded index. The analog fails for linear groups since, for example, $GL(3; \mathbb{Z})$ contains the Heisenberg group.

The outer automorphism group $Out(F_n)$ of a free group $F_n$ of finite rank $n$ reflects the nature of both linear and mapping class groups. Indeed, it maps onto $GL(n; \mathbb{Z})$ and contains the mapping class group $MCG(S)$ of a compact surface $S$ with fundamental group $F_n$. E. Formanek and C. Procesi [FP] have shown that $Out(F_n)$ is not linear if $n > 3$. It is unknown if mapping class groups of compact surfaces are all linear. In a series of two papers we prove:

THEOREM 1.0.1. *The group* $Out(F_n)$ *satisfies the Tits alternative.*

In a third paper [BFH2], we prove the following complementary result.

THEOREM 1.0.2. *A solvable subgroup of* $Out(F_n)$ *has a finitely generated free abelian subgroup of index at most* $3^{5n^2}$.

The rank of an abelian subgroup of $Out(F_n)$ is bounded by $vcd(Out(F_n)) = 2n - 3$ for $n > 1$ [CV]. With regard to the relationship between solvable and abelian subgroups, $Out(F_n)$ behaves like $MCG(S)$. H. Bass and A. Lubotzky [BL] showed that solvable subgroups of $Out(F_n)$ are virtually polycyclic.

Theorem 1.0.1 is divided into two parts according to the growth rate of the automorphisms being considered. An element $\mathcal{O}$ of $Out(F_n)$ has *polynomial growth* if for each conjugacy class $[[\gamma]]$ of an element in $F_n$ the word length of $\mathcal{O}^i([[\gamma]])$ with respect to some fixed finite generating set for $F_n$ grows at most polynomially in $i$. An element $\mathcal{O}$ of $Out(F_n)$ has *exponential growth* if for some conjugacy class this sequence grows at least exponentially in $i$. An element of $Out(F_n)$ has either polynomial or exponential growth (see for example [BH1]). The set of outer automorphisms that has polynomial growth is denoted $PG(F_n)$; the set that have polynomial growth and unipotent image in $GL(n; \mathbb{Z})$ is denoted $UPG(F_n)$. A subgroup of $Out(F_n)$ is said to be PG (respectively UPG) if all of its elements are contained in $PG(F_n)$ (respectively UPG$(F_n)$). Every PG subgroup contains a finite index UPG subgroup (Corollary 5.7.6).

Every UPG subgroup of $\mathrm{MCG}(S)$ is abelian. In fact, each UPG subgroup of $\mathrm{MCG}(S)$ is contained in a group generated by Dehn twists in a set of pairwise disjoint simple closed curves [Iva], [BLM]. The structure of UPG subgroups of $\mathrm{Out}(F_n)$ is richer. In particular, they may contain free subgroups of rank 2; see Remark 1.2 of [BFH1]. The second paper in this series [BFH1] is a study of UPG subgroups of $\mathrm{Out}(F_n)$. It contains a proof of the following theorem.

THEOREM 1.0.3. *A UPG subgroup of $\mathrm{Out}(F_n)$ that does not contain a free subgroup of rank 2 is solvable.*

This, the first paper in this series, culminates in the following theorem. Theorem 1.0.1 is an immediate consequence of it and Theorem 1.0.3.

THEOREM 7.0.1. *Suppose that $\mathcal{H}$ is a subgroup of $\mathrm{Out}(F_n)$ that does not contain a free subgroup of rank 2. Then there is a finite index subgroup $\mathcal{H}_0$ of $\mathcal{H}$, a finitely generated free abelian group $A$, and a map*

$$\Phi : \mathcal{H}_0 \to A$$

*such that $\mathrm{Ker}(\Phi)$ is UPG.*

In [BFH4] (see also [BFH3]) which is independent of the current series, there is a proof of the Tits alternative for a special class of subgroups of $\mathrm{Out}(F_n)$.

Although our work focuses on the Tits alternative, our approach has always been toward developing a general understanding of subgroups of $\mathrm{Out}(F_n)$ and their dynamics on certain spaces of trees and bi-infinite paths. In the remainder of this section and in the introduction to [BFH1], we take up this general viewpoint.

We establish our dynamical point of view by recalling an experiment described by Thurston. Suppose that $S$ is a compact surface equipped with a complete hyperbolic metric and that $\phi$ is an element of the mapping class group $\mathrm{MCG}(S)$. Each free homotopy class of closed curves in $S$ is represented by a unique closed geodesic. This determines a natural action of $\phi$ on the set of closed geodesics in $S$ and we denote the image of the geodesic $\sigma$ under this action by $\phi_\#(\sigma)$.

Choose a closed geodesic $\sigma$ and positive integer $k$. Using a fine point, draw $\phi_\#^k(\sigma)$ on $S$ and step back so that you can no longer see individual drawn lines but only the places where lines accumulate. If $\sigma$ is periodic under the action of $\phi$, then you will not see anything. In all other cases, as $k$ increases the image will stabilize and you will see a nonempty closed set $V(\sigma)$ of disjoint simple geodesics. Most $\sigma$ produce the same stabilized image and we denote this by $V(\phi)$. The exceptional cases produce $V(\sigma)$ that are subsets of $V(\phi)$.

This experiment neatly captures the essential features of Thurston's normal form for elements of MCG($S$) ([Thu]; see also [FLP]). For each $\phi \in$ MCG($S$), there is a canonical decomposition of $S$ along a (possibly empty) set of disjoint annuli $A_j$ into subsurfaces $S_i$ of negative Euler characteristic. The mapping class $\phi$ restricts to a mapping class on each $S_i$ that either has finite order or is pseudo-Anosov. On each $A_j$, $\phi$ restricts to a (possibly trivial) Dehn twist. If $\phi|S_i$ is pseudo-Anosov, denote the associated attracting geodesic lamination by $\Gamma_i^+$; if $\phi|A_j$ is a nontrivial Dehn twist, denote the core geodesic of $A_j$ by $\alpha_j$. Then $V(\phi)$ is the union of the $\Gamma_i^+$'s and $\alpha_j$'s. Each $V(\sigma)$ is a union of $\Gamma_i^+$'s and $\alpha_j$'s; more precisely, $\Gamma_i^+$ (respectively $\alpha_j$) is contained in $V(\sigma)$ if and only if $\Gamma_i^+$ (respectively $\alpha_j$) has nonempty transverse intersection with $\sigma$.

Relative train track maps $f : G \to G$ were introduced in [BH1] as the Out($F_n$) analog of the Thurston normal form. An outer automorphism $\mathcal{O}$ is represented by a homotopy equivalence $f : G \to G$ of a marked graph and a filtration $\emptyset = G_0 \subset G_1 \subset \cdots \subset G_K = G$ by $f$-invariant subgraphs. Thus we view $\mathcal{O}$ as being built up in stages. The marked graph $G$ is broken up into strata $H_i$ (the difference between $G_i$ and $G_{i-1}$) that are, in some ways, analogous to the $S_i$'s and $A_j$'s that are part of the Thurston normal form for $\phi \in$ MCG($S$).

Associated to each stratum $H_r$ is an integral transition matrix which records the multiplicity of the $i^{\text{th}}$ edge of $H_r$ in the image under $f$ of the $j^{\text{th}}$ edge of $H_r$. If the filtration is sufficiently refined, then these matrices are either identically zero or irreducible (see subsection 2.5). We call strata of the former type zero strata and strata of the latter type irreducible. To each irreducible stratum we associate the Perron-Frobenius eigenvalue of its transition matrix.

Irreducible strata are said to be non-exponentially-growing or exponentially-growing according to whether their associated Perron-Frobenius eigenvalues are, respectively, equal to one or greater than one. Exponentially-growing strata correspond to pseudo-Anosov components. There are three types of non-exponentially-growing strata. If $f$ acts periodically on the edges of $H_i$, then $H_i$ is analogous to a subsurface $S_i$ on which $\phi$ acts periodically. If the lengths of the edges of $H_i$ grow linearly under iteration by $f$, then $H_i$ is analogous to an annulus with nontrivial Dehn twisting. If the lengths of the edges of $H_i$ have a faster than linear growth rate, then $H_i$ has no surface counterpart. Zero strata play a lesser role in the theory.

The space $\mathcal{B}(G)$ (see §2) of bi-infinite unoriented paths (hereafter referred to as lines) in a marked graph $G$ is the $F_n$ analog of the space $\mathcal{G}(S)$ of complete geodesics in $S$. Periodic lines are called circuits and correspond to closed geodesics. There is a natural action of $\mathcal{O}$ on $\mathcal{B}(G)$. Since one cannot directly 'see' lines in $G$, we pose the analogy for the experiment as follows. Given a circuit $\gamma$, what are the accumulation points in $\mathcal{B}(G)$ of the forward $\mathcal{O}$-orbit of $\gamma$? This is not a completely faithful translation. Geodesics that are contained in $S_i \setminus \Gamma_i^+$ occur as accumulation points for the forward $\phi_\#$-orbit of certain $\sigma$ but are not contained in $V(\phi)$.

An exhaustive study of the action of $\mathcal{O}$ on $\mathcal{B}(G)$ is beyond the scope of any single paper. Our goal is to build a general framework for the subject with sufficient detail to prove Theorem 7.0.1. In some cases we develop an idea beyond what is required for the Tits alternative and in some cases we do not. Our decisions are based not only on the relative importance of the idea but also on the number of pages required to do the extra work.

The key dynamical invariant introduced in this paper is the *attracting lamination* associated to an exponentially growing stratum of a relative train track map $f : G \to G$. It is the analog of the unstable measured geodesic lamination $\Gamma_i^+$ associated to a pseudo-Anosov component of a mapping class element. We take a purely topological point of view and define these laminations to be closed sets in $\mathcal{B}(G)$; measures are not considered in this paper. To remind the reader that we are not working in a more structured space (and in fact are working in a non-Hausdorff space), we use the term *weak attraction* when describing limits in $\mathcal{B}(G)$. Thus a line $L_1$ is weakly attracted to a line $L_2$ under the action of $\mathcal{O}$ if for every neighborhood $U$ of $L_2$ in $\mathcal{B}(G)$, there is a positive integer $K$ so that $\mathcal{O}_\#^k(L_1) \in U$ for all $k > K$.

The set $\mathcal{L}(\mathcal{O})$ of attracting laminations associated to the exponentially growing strata of a relative train track map $f : G \to G$ representing $\mathcal{O}$ is finite (Lemma 3.1.13) and is independent of the choice of $f : G \to G$. After passing to an iterate if necessary, we may assume that each element of $\mathcal{L}(\mathcal{O})$ is $\mathcal{O}$-invariant.

An attracting lamination $\Lambda^+$ has preferred lines, called *generic lines*, that are dense in $\Lambda^+$ (Lemma 3.1.15). All generic lines have the same neighborhoods in $\mathcal{B}(G)$ (Corollary 3.1.11) and so weakly attract the same lines. We refer to this common set of weakly attracted lines as the basin of weak attraction for $\Lambda^+$. An element of $\mathcal{L}(\mathcal{O})$ is *topmost* if it is not contained in any other element of $\mathcal{L}(\mathcal{O})$.

Of central importance to our study is the following question. Which circuits (and more generally which bierecurrent lines (Definition 3.1.3)) are contained in the basin of weak attraction for a topmost $\Lambda^+$?

A first guess might be that a circuit $\gamma$ is weakly attracted to $\Lambda^+$ if and only if it intersects the stratum $H_r$ that determines $\Lambda^+$. This fails in two ways. First, strata are not invariant; the $f$-image of an edge in $H_i$, $i > r$, may contain edges in $H_r$. Thus $\Lambda^+$ may attract circuits that do not intersect $H_r$. For the second, suppose that $\phi : S \to S$ is a pseudo-Anosov homeomorphism of a compact surface with one boundary component. If $\pi_1(S)$ is identified with $F_n$, then the outer automorphism determined by $\phi$ is represented by a relative train track map with a single stratum. This stratum is exponentially growing and so determines an attracting lamination $\Gamma^+$. The only circuit not attracted to $\Gamma^+$ is the one, say $\rho$, determined by $\partial S$. Since $\rho$ crosses every edge in $G$ twice, one cannot expect to characterize completely the basin of weak attraction in terms of a subgraph of $G$.

The key ingredient in analyzing the basin of attraction for a pseudo-Anosov lamination $\Gamma^+$ is intersection of geodesics: a circuit is attracted to $\Gamma^+$ if and only if it intersects the dual lamination $\Gamma^-$. Unfortunately, intersection of geodesics has no analog in $\mathrm{Out}(F_n)$. Indeed, this is a frequently encountered stumbling block in generalizing from $\mathrm{MCG}(S)$ to $\mathrm{Out}(F_n)$. We overcome this by modifying and improving the relative train track methods of [BH1] and by a detailed analysis of the action of $f$ on paths in $G$. Most of this analysis is contained in Section 5. A very detailed statement of our improved relative train track maps is given in Theorem 5.1.5 and we refer the reader to the introduction of Section 5 for an overview of its contents. We believe that improved relative train tracks are important in their own right and will be useful in solving other problems (see, for example, [Mac1], [Mac2] and [Bri]).

It is shown in subsection 3.2 that there is a pairing between elements of $\mathcal{L}(\mathcal{O})$ and elements of $\mathcal{L}(\mathcal{O}^{-1})$ that is analogous to the pairing between stable and unstable pseudo-Anosov laminations. The paired laminations are denoted $\Lambda^+$ and $\Lambda^-$ with $\Lambda^+ \in \mathcal{L}(\mathcal{O})$ and $\Lambda^- \in \mathcal{L}(\mathcal{O}^{-1})$.

Theorem 5.1.5 shows that any outer automorphism $\mathcal{O}$ can be realized by an 'improved relative train track map' and an associated filtration. Each exponentially-growing stratum $H_r$ of an improved relative train track map has a canonically associated (see Lemma 4.2.5) finite set $P_r$ of paths in $G_r$. If $P_r$ is nonempty, then it contains a preferred element $\rho_r$. We assign a path $\hat{\rho}_r$ to each exponentially-growing stratum $H_r$ as follows. If $P_r \neq \emptyset$, then $\hat{\rho}_r = \rho_r$. If $P_r = \emptyset$, then choose a vertex in $H_r$ and define $\hat{\rho}_r$ to be the trivial path at that vertex.

For any subgraph $X$ of $G$ and finite path $\rho \subset G$, define $\langle X, \rho \rangle$ to be the groupoid of paths that can be decomposed into a concatenation of subpaths that are either entirely contained in $X$ or are equal to $\rho$ or $\bar{\rho}$. If $\rho$ is a trivial path, then a nontrivial path in $G$ is contained in $\langle X, \rho \rangle$ if and only if it is contained in $X$. The following theorem is one of the two main results in this paper. Although the statement is completely analogous to a well-known result about mapping classes, the proof is entirely different.

THEOREM 6.0.1 (Weak Attraction Theorem). *Suppose that $\Lambda^+$ is a topmost element of $\mathcal{L}(\mathcal{O})$, that $f : G \to G$ is an improved relative train track map representing $\mathcal{O}$ and that $H_r$ is the exponentially-growing stratum that determines $\Lambda^+$. Then there exists a subgraph $Z$ such that $Z \cap G_r = G_{r-1}$ and such that every bicurrent path $\gamma \subset G$ satisfies exactly one of the following.*

1. $\gamma$ *is a generic line for* $\Lambda^-$.

2. $\gamma \in \langle Z, \hat{\rho}_r \rangle$.

3. $\gamma$ *is weakly attracted to* $\Lambda^+$.

To study the iterated images of a bi-infinite path $\gamma$, we subdivide it into 'non-interacting' subpaths whose behavior under iteration is largely determined by a single stratum. This *splitting* is the subject of subsection 4.1 and parts of Section 5. Roughly speaking, one can view this as the analog of subdividing a geodesic in $S$ according to its intersections with the $S_i$'s and $A_j$'s that are part of Thurston's normal form.

There are three parts to our proof of Theorem 7.0.1. First, we use the well known 'ping-pong' method of Tits (Proposition 1.1 of [Tit]) to establish a criterion for a subgroup $\mathcal{H}$ of $\mathrm{Out}(F_n)$ to contain a free subgroup of rank two.

COROLLARY 3.4.3.   *Suppose that* $\Lambda^+ \in \mathcal{L}(\mathcal{O})$ *and* $\Lambda^- \in \mathcal{L}(\mathcal{O}^{-1})$ *are paired and* $\mathcal{O}$-*invariant, that* $\mathcal{H}$ *is a subgroup of* $\mathrm{Out}(F_n)$ *containing* $\mathcal{O}$ *and that there is an element* $\psi \in \mathcal{H}$ *such that generic lines of the four laminations* $\psi^{\pm 1}(\Lambda^\pm)$ *are weakly attracted to* $\Lambda^+$ *under the action of* $\mathcal{O}$ *and are weakly attracted to* $\Lambda^-$ *under the action of* $\mathcal{O}^{-1}$. *Then* $\mathcal{H}$ *contains a free subgroup of rank two.*

In Section 7, we combine this criterion with the weak attraction theorem and a homology argument to prove the following.

LEMMA 7.0.10.   *If* $\mathcal{H} \subset \mathrm{Out}(F_n)$ *does not contain a free subgroup of rank two, then there is a finite collection* $\mathcal{L}$ *of attracting laminations for elements of* $\mathcal{H}$ *and a finite index subgroup* $\mathcal{H}_0$ *of* $\mathcal{H}$ *that stabilizes each element of* $\mathcal{L}$ *and that satisfies the following property. If* $\psi \in \mathcal{H}_0$ *and if* $\Lambda^+ \in \mathcal{L}(\psi)$ *and* $\Lambda^- \in \mathcal{L}(\psi^{-1})$ *are paired topmost laminations, then at least one of* $\Lambda^+$ *and* $\Lambda^-$ *is in* $\mathcal{L}$.

The last ingredient of the proof of Theorem 7.0.1 is contained in subsection 3.3. Denote the stabilizer in $\mathrm{Out}(F_n)$ of an attracting lamination $\Lambda^+$ by $\mathrm{Stab}(\Lambda^+)$.

COROLLARY 3.3.1.   *There is a homomorphism* $PF_{\Lambda^+} : \mathrm{Stab}(\Lambda^+) \to \mathbb{Z}$ *such that* $\Psi \in \mathrm{Ker}(PF_{\Lambda^+})$ *if and only if* $\Lambda^+ \notin \mathcal{L}(\Psi)$ *and* $\Lambda^+ \notin \mathcal{L}(\Psi^{-1})$.

The analogous result for the mapping class group is an immediate corollary of the fact (exposé 12 of [FLP]) that the measured foliations associated to a pseudo-Anosov homeomorphism are uniquely ergodic. Any mapping class that topologically preserves the measured foliation must projectively fix its invariant transverse measure and so multiplies this transverse measure by some scalar factor. The assignment of the logarithm of this scalar factor to the mapping class defines the analogous homomorphism.

Because we are working in $\mathcal{B}$ and not a more structured space that takes measures into account, we cannot measure the attraction factor directly. Instead of an invariant measure defined on the lamination itself, we use a length function on paths in a marked graph. The length function depends on the choice of marked graph but the factor by which an element of $\mathrm{Stab}(\Lambda^+)$ expands this length does not.

The three parts to the proof of Theorem 7.0.1 are tied together at the end of Section 7.

## 2. Preliminaries

2.1. *Marked graphs and topological representatives.* A *marked graph* is a graph $G$ along with a homotopy equivalence $\tau : R_n \to G$ from the rose $R_n$ with $n$ petals and vertex $*$. We assume that $F_n$ is identified with $\pi_1(R_n, *)$ and hence also with $\pi_1(G, \tau(*))$. A homotopy equivalence $f : G \to G$ induces an outer automorphism of $\pi_1(G, \tau(*))$ and so an outer automorphism $\mathcal{O}$ of $F_n$. The set of vertices of $G$ is denoted $\mathcal{V}$. If $f(\mathcal{V}) \subset \mathcal{V}$ and if the restriction of $f$ to each edge of $G$ is an immersion, then we say that $f : G \to G$ is *a topological representative of $\mathcal{O}$.*

A *filtration* for a topological representative $f : G \to G$ is an increasing sequence of (not necessarily connected) $f$-invariant subgraphs $\emptyset = G_0 \subset G_1 \subset \cdots \subset G_K = G$. The closure $H_r$ of $(G_r \setminus G_{r-1})$ is a subcomplex called the $r^{\mathrm{th}}$ *stratum.*

*Throughout this paper, $G$ will be a marked graph, $f : G \to G$ will be a topological representative, $G_r$ will be a filtration element and $H_r$ will be a filtration stratum. The universal cover of $G$ is a tree denoted by $\Gamma$.*

2.2. *Paths, circuits and lines.* In this subsection we set notation for our treatment of 'geodesics'.

Let $\Gamma$ be the universal cover of a marked graph $G$ and let $\mathrm{pr} : \Gamma \to G$ be the covering projection. A map $\tilde{\alpha} : J \to \Gamma$ with domain a (possibly infinite) interval $J$ will be called a *path in $\Gamma$* if it is an embedding or if $J$ is finite and the image is a single point; in the latter case we say that $\tilde{\alpha}$ is *a trivial path.*

We will not distinguish between paths in $\Gamma$ that differ only by an orientation-preserving change of parametrization. Thus we are interested in the oriented image of $\tilde{\alpha}$ and not $\tilde{\alpha}$ itself. If the domain of $\tilde{\alpha}$ is finite, then the image of $\tilde{\alpha}$ has a natural decomposition as a concatenation $\tilde{E}'_1 \tilde{E}_2 \cdots \tilde{E}_{k-1} \tilde{E}'_k$ where $\tilde{E}_i$, $1 < i < k$, is a directed edge of $\Gamma$, $\tilde{E}'_1$ is the terminal segment of a directed edge $\tilde{E}_1$ and $\tilde{E}'_k$ is the initial segment of a directed edge $\tilde{E}_k$. If the endpoints of the image of $\tilde{\alpha}$ are vertices, then $\tilde{E}'_1 = \tilde{E}_1$ and $\tilde{E}'_k = \tilde{E}_k$. The sequence $\tilde{E}'_1 \tilde{E}_2 \cdots \tilde{E}'_k$ is called *the edge path associated to $\tilde{\alpha}$.* This notation extends naturally to the case that the domain is a ray or the entire line. In the former case, an edge path has the *half-infinite* form $\tilde{E}'_1\tilde{E}_2\cdots$ or $\cdots\tilde{E}_{-2}\tilde{E}'_{-1}$ and in the latter case has the *bi-infinite* form $\cdots\tilde{E}_{-1}\tilde{E}_0\tilde{E}_1\tilde{E}_2\cdots$.

If $J$ is finite, then every map $\tilde{\alpha}: J \to \Gamma$ is homotopic rel endpoints to a unique (possibly trivial) path $[\tilde{\alpha}]$; we say that $[\tilde{\alpha}]$ *is obtained from $\tilde{\alpha}$ by tightening*. If $\tilde{f}: \Gamma \to \Gamma$ is a lift of $f: G \to G$ , we denote $[\tilde{f}(\tilde{\alpha})]$ by $\tilde{f}_{\#}(\tilde{\alpha})$.

A *path in $G$* is the composition of the projection map pr with a path in $\Gamma$. Thus a map $\alpha: J \to G$ with domain a (possibly infinite) interval will be called a path if it is an immersion or if $J$ is finite and the image is a single point; paths of the latter type are said to be trivial. If $J$ is finite, then every map $\alpha: J \to G$ is homotopic rel endpoints to a unique (possibly trivial) path $[\alpha]$; we say that $[\alpha]$ *is obtained from $\alpha$ by tightening*. For any lift $\tilde{\alpha}: J \to \Gamma$ of $\alpha$, $[\alpha] = pr[\tilde{\alpha}]$. We denote $[f(\alpha)]$ by $f_{\#}(\alpha)$.

We do not distinguish between paths in $G$ that differ by an orientation-preserving change of parametrization. The *edge path associated to $\alpha$* is the projected image of the edge path associated to a lift $\tilde{\alpha}$. Thus the edge path associated to a path with finite domain has the form $E'_1E_2\cdots E_{k-1}E'_k$ where $E_i$, $1 < i < k$, is an edge of $G$, $E'_1$ is the terminal segment of an edge $E_1$ and $E'_k$ is the initial segment of an edge $E_k$.

We reserve the word *circuit* for an immersion $\alpha: S^1 \to G$. Any homotopically nontrivial map $\sigma: S^1 \to G$ is homotopic to a unique circuit $[\sigma]$. As was the case with paths, we do not distinguish between circuits that differ only by an orientation-preserving change in parametrization and we identify a circuit $\alpha$ with a *cyclically ordered edge path* $E_1E_2\ldots E_k$.

*Throughout this paper we will identify paths and circuits with their associated edge paths.*

For any path $\alpha$ in $G$ define $\bar{\alpha}$ to be $\alpha$ *with its orientation reversed*. To make this precise choose an orientation-reversing homeomorphism inv as follows. If $J$ is either finite or bi-infinite, then inv : $J \to J$; if $J = (-\infty, b]$ , then inv : $[b, \infty) \to (-\infty, b]$; if $J = [a, \infty)]$ then inv : $(-\infty, a] \to [a, \infty)$. Define $\bar{\alpha} = \alpha \circ$ inv. We sometimes refer to $\bar{\alpha}$ as the *inverse of $\alpha$*. The inverse of a path in $\Gamma$ is defined similarly.

There are times when we want to ignore a path's orientation. In these cases we will refer to $\alpha$ or $\tilde{\alpha}$ as *an unoriented path*. We reserve the word *line* for an unoriented bi-infinite path. If $\tilde{\alpha}$ contains $\tilde{\alpha}_0$ or its inverse as a subpath, then we say that $\tilde{\alpha}_0$ is an *unoriented subpath of $\tilde{\alpha}$*. If $\tilde{\alpha}$ is a line, then we sometimes simply write that $\tilde{\alpha}_0$ is a subpath of $\tilde{\alpha}$ since the lack of orientation is implicit in the fact that $\tilde{\alpha}$ is unoriented. Similar notation is used for unoriented subpaths in $G$.

*The space of lines in $\Gamma$* is denoted $\tilde{\mathcal{B}}(\Gamma)$ and is equipped with what amounts to the compact-open topology. Namely, for any finite path $\tilde{\alpha}_0 \subset \Gamma$ (with endpoints at vertices if desired), define $N(\tilde{\alpha}_0) \subset \tilde{\mathcal{B}}(\Gamma)$ to be the set of lines in $\Gamma$ that contain $\tilde{\alpha}_0$ as a subpath. The sets $N(\tilde{\alpha}_0)$ define a basis for the topology on $\mathcal{B}(\Gamma)$.

*The space of lines in $G$* is denoted $\mathcal{B}(G)$. There is a natural projection map from $\tilde{\mathcal{B}}(\Gamma)$ to $\mathcal{B}(G)$ and we equip $\mathcal{B}(G)$ with the quotient topology. A basis for the topology is constructed by considering finite paths $\alpha_0$ (with endpoints at vertices if desired) and defining $N(\alpha_0) \subset \mathcal{B}(G)$ to be the set of lines in $G$ that contain $\alpha_0$ as a subpath.

In the analogy with the mapping class group, $\mathcal{B}(G)$ corresponds to the space of complete geodesics in a closed surface $S$ equipped with a particular hyperbolic metric; $\tilde{\mathcal{B}}(\Gamma)$ corresponds to the space of complete geodesics in the universal cover $\tilde{S}$.

Nielsen's approach to the mapping class group (see [HT] for example) begins with the fact that each mapping class $\phi$ determines a homeomorphism $\phi_\#$ on the space of complete geodesics in $S$. This can be briefly described as follows. The universal cover $\tilde{S}$ is compactified by a 'circle at infinity' $S_\infty$ in such a way that complete geodesics in $\tilde{S}$ correspond to distinct pairs of points in $S_\infty$. One proves that if $\tilde{h} : \tilde{S} \to \tilde{S}$ is any lift of a homeomorphism $h : S \to S$ representing $\phi$, then $\tilde{h}$ extends to a homeomorphism of $S_\infty$. Since $\tilde{h}$ induces an equivariant homeomorphism on pairs of points in $S_\infty$, it induces an equivariant homeomorphism $\tilde{h}_\#$ on the space of geodesics in $S_\infty$ and a homeomorphism $h_\#$ on the space of complete geodesics in $S$. One then checks that $h_\#$ depends only on $\phi$ and not on the choices of $h$ and $\tilde{h}$.

There are analogous results for $\mathrm{Out}(F_n)$. The circle at infinity is replaced by the Cantor set $\partial F_n$ of ends of $F_n$. We assume from now on that the basepoint in $G$ has been lifted to a basepoint in $\Gamma$. The marking on $G$ then determines a homeomorphism between the space of ends of $\Gamma$ and $\partial F_n$ (see, for example, [Flo]. We use this identification and treat $\partial F_n$ as the space of ends of $\Gamma$.

*Definition* 2.2.1. Define $\tilde{\mathcal{B}} = (\partial F_n \times \partial F_n \setminus \Delta)/Z_2$, where $\Delta$ is the diagonal and where $Z_2$ acts on $\partial F_n \times \partial F_n$ by interchanging the factors. For any unordered pair of distinct elements $(c_1, c_2) \in \partial F_n \times \partial F_n$ and for any $\Gamma$, there is a unique line $\tilde{\sigma} \subset \Gamma$ connecting the ends $c_1$ and $c_2$. This process is reversible and defines a homeomorphism between $\tilde{\mathcal{B}}$ and $\tilde{\mathcal{B}}(\Gamma)$. We will often use this homeomorphism implicitly to identify $\tilde{\mathcal{B}}$ and $\tilde{\mathcal{B}}(\Gamma)$.

The diagonal action of $F_n$ on $\partial F_n \times \partial F_n$ defines an action of $F_n$ on $\tilde{\mathcal{B}}$. Define $\mathcal{B}$ to be the quotient space of this action. The action of $F_n$ on $\Gamma$ by covering translations defines an action of $F_n$ on $\tilde{\mathcal{B}}(\Gamma)$. The homeomorphism between $\tilde{\mathcal{B}}$ and $\tilde{\mathcal{B}}(\Gamma)$ is $F_n$-equivariant and so projects to a homeomorphism between $\mathcal{B}$ and $\mathcal{B}(G)$. We will often use this homeomorphism implicitly to identify $\mathcal{B}$ and $\mathcal{B}(G)$. If $\gamma \in \mathcal{B}(G)$ corresponds to $\beta \in \mathcal{B}$ then we say that $\gamma$ *realizes $\beta$ in $G$*. In the analogy with the mapping class group, $\mathcal{B}$ corresponds to an abstract space of complete geodesics in $S$ that is independent of the choice of hyperbolic metric.

*Definition* 2.2.2. Assume that the space of ends of $\Gamma$ and the space of ends of $\Gamma'$ have been identified with $\partial F_n$. If $\tilde{h} : \Gamma \to \Gamma'$ is a lift of a homotopy equivalence $h : G \to G'$ then (page 208 of [Flo]) $\tilde{h}$ determines a homeomorphism $\tilde{h} : \partial F_n \to \partial F_n$. There are *induced homeomorphisms* $\tilde{h}_{\#} : \tilde{\mathcal{B}}(\Gamma) \to \tilde{\mathcal{B}}(\Gamma')$ *and* $h_{\#} : \mathcal{B}(G) \to \mathcal{B}(G')$. If $\tilde{\alpha}$ is a line in $\Gamma$ with endpoints $P, Q \in \partial F_n$, then $\tilde{h}_{\#}(\tilde{\alpha})$ is the line in $\Gamma'$ with endpoints $\tilde{h}(P), \tilde{h}(Q)$.

Circuits correspond to periodic bi-infinite paths in $G$. We sometimes use this correspondence to think of the circuits as a subset of $\mathcal{B}$. Since every finite path $\alpha_0 \subset R_n$ extends to a circuit, the circuits form a dense set in $\mathcal{B}$. One may also identify the circuits with the set of conjugacy classes $[[a]]$ in $F_n$. (This is analogous to the fact that every free homotopy class of closed curves in a hyperbolic surface contains a unique geodesic.) An outer automorphism $\mathcal{O}$ determines an action $\mathcal{O}_{\#}$ on conjugacy classes in $F_n$ and hence on the set of circuits.

Our various definitions are tied together by the following lemma.

LEMMA 2.2.3. *Suppose that $h : G \to G'$ is a homotopy equivalence of marked graphs and that $\mathcal{O}$ is the outer automorphism determined by $h$. Then*

1. *The action induced by $h_{\#} : \mathcal{B}(G) \to \mathcal{B}(G')$ on circuits is given by $\alpha \mapsto [h(\alpha)]$.*

2. *The action induced by $h_{\#} : \mathcal{B} \to \mathcal{B}$ on conjugacy classes in $F_n$ is given by $[[a]] \mapsto \mathcal{O}_{\#}([[a]])$.*

3. *$h_{\#} : \mathcal{B} \to \mathcal{B}$ is determined by the action of $\mathcal{O}$ on circuits.*

*Proof.* Let $\tilde{\alpha} \subset \Gamma$ be a lift of a circuit $\alpha \subset G$ and let $\tilde{h} : \Gamma \to \Gamma'$ be a lift of $h : G \to G'$. A homotopy between $h(\alpha)$ and $\alpha' = [h(\alpha)]$ lifts to a bounded homotopy between $\tilde{h}(\tilde{\alpha})$ and a lift $\tilde{\alpha}'$ of $\alpha'$. This implies that $\tilde{h}(\tilde{\alpha})$ and $\tilde{\alpha}'$ have the same endpoints in $\partial F_n$ and hence that $\tilde{h}_{\#}(\tilde{\alpha}) = \tilde{\alpha}'$. Part 1 follows immediately.

Part 2 follows immediately from part 1 and the definitions. Part 3 follows from part 2 and the denseness of circuits in $\mathcal{B}$. $\square$

2.3. *The bounded cancellation lemma.* In this section we state the bounded cancellation lemma of [Coo] in the forms used in this paper. A generalization of the bounded cancellation lemma is given in [BFH4].

LEMMA 2.3.1. *For any homotopy equivalence $h : G \to G'$ of marked graphs there is a constant $C$ with the following properties.*

1. If $\rho = \alpha\beta$ is a path in $G$, then $h_\#(\rho)$ is obtained from $h_\#(\alpha)$ and $h_\#(\beta)$ by concatenating and by cancelling $c \leq C$ edges from the terminal end of $h_\#(\alpha)$ with $c$ edges from the initial end of $h_\#(\beta)$.

2. If $\tilde{h} : \Gamma \to \Gamma'$ is a lift to the universal covers, $\tilde{\alpha}$ is a line in $\Gamma$ and $\tilde{x} \in \tilde{\alpha}$, then $\tilde{h}(\tilde{x})$ can be connected to $\tilde{h}_\#(\tilde{\alpha})$ by a path with $c \leq C$ edges.

3. Suppose that $\tilde{h} : \Gamma \to \Gamma'$ is a lift to the universal covers and that $\tilde{\alpha} \subset \Gamma$ is a finite path. Define $\tilde{\beta} \subset \Gamma'$ by removing $C$ initial and $C$ terminal edges from $\tilde{h}_\#(\tilde{\alpha})$. Then $\tilde{h}_\#(N(\tilde{\alpha})) \subset N(\tilde{\beta})$. (In other words, if $\tilde{\gamma} \in \tilde{\mathcal{B}}(\Gamma)$ contains $\tilde{\alpha}$ as a subpath, then $h_\#(\tilde{\gamma}) \in \tilde{\mathcal{B}}(\Gamma')$ contains $\tilde{\beta}$ as a subpath.)

2.4. *Folding.* We now recall the folding construction of Stallings [Sta]. Suppose that $f : G \to G$ is a topological representative of $\mathcal{O}$. If $f$ is not an immersion, then there is a pair of distinct oriented edges $E_1$ and $E_2$ with the same initial endpoint and there are nontrivial initial segments $E_1^* \subset E_1$ and $E_2^* \subset E_2$ such that $f(E_1^*) = f(E_2^*)$ is a path with endpoints at vertices. There exists a surjective *folding* map $p : G \to G^1$ defined by identifying $E_1^*$ with $E_2^*$ in such a way that $f$ factors as $f = gp$ for some map $g : G^1 \to G$.

Since $f$ is a homotopy equivalence and $f_\#(\bar{E}_1^* E_2^*)$ is trivial, $\bar{E}_1^* E_2^*$ is not a closed path. Let $T$ be a triangle fibered by lines parallel to its base. Attach $T$ to $G$ so that the non-base sides are identified with $E_1^*$ and $E_2^*$ and so that the endpoints of each fiber are identified by $p_1$. The resulting space $X$ deformation retracts to $G$. Collapsing the fibers of $T$ to points defines a homotopy equivalence of $X$ to $G^1$. Moreover, the inclusion of $G$ into $X$ followed by the collapsing of the fibers agrees with $p$. Thus $p : G \to G^1$ is a homotopy equivalence.

We will apply this construction in two ways. In the first, we produce a new topological representative of $\mathcal{O}$ as follows. Define $f_1 : G^1 \to G^1$ by 'tightening' $pg : G^1 \to G^1$; i.e. by defining $f_1(e) = (pg)_\#(e)$ for each edge $e$ of $G_1$. If each $f_1(e)$ is nontrivial, we are done. If not, the set of edges with trivial $f_1$-image form a tree and we collapse each component of the tree to a point. After repeating this tighten and collapse procedure finitely many times, we arrive at the desired topological representation.

For the second application, the folding operation is repeated with $g : G^1 \to G$ replacing $f : G \to G$ and so on to conclude that $f = \theta p_k \ldots p_1$ where $G^0 = G$, $p_i : G^{i-1} \to G^i$ is a folding map and where $\theta : G^k \to G$ is an immersion. The immersion $\theta$ extends to a covering $\hat{\theta} : \hat{G} \to G$. Since $\theta$ is a homotopy equivalence, $\hat{\theta}$ must be degree one and $G \subset \hat{G}$ is a homotopy equivalence. In other words $\theta$ is an embedding and is a homeomorphism if $G$ has no valence-one vertices [Sta].

We also need a slight generalization of folding. Suppose that $E_2 = \mu_1\mu_2$ is a decomposition into subpaths and that $\sigma \subset G$ is a path satisfying the following properties : $\sigma$ and $E_2$ have a common initial endpoint; $\sigma$ does not intersect the interior of $E_2$; and $f(\mu_1) = f_\#(\sigma)$ is a path with endpoints at vertices. Define $G'$ to be the graph obtained from $G$ by identifying $\mu_1$ with $\sigma$ and let $p : G \to G'$ be the quotient map. We may think of $G \setminus E_2$ as a subcomplex of both $G$ and $G'$. Thus $G$ is obtained from $G \setminus E_2$ by adding $E_2$ and $G'$ is obtained from $G \setminus E_2$ by adding an edge $E_2'$ with terminal endpoint equal to the terminal endpoint of $E_2$ and initial endpoint equal to the terminal endpoint of $\sigma$. With this notation, $p|(G \setminus E_2)$ is the identity, $p(\mu_1) = \sigma$ and $p(\mu_2) = E_2'$. Define $g : G' \to G$ by $g|(G \setminus E_2) = f|(G \setminus E_2)$ and by $g(E_2') = f(\mu_2)$. Then $gp|(G \setminus E_2) = f|(G \setminus E_2)$ and $(gp)_\#(E_2) = f(E_2)$. In particular, $gp \simeq f$ rel $\mathcal{V}$ (= the vertex set of $G$). We refer to $p : G \to G'$ as a *generalized fold*.

2.5. *Relative train track maps.* We study an outer automorphism by analyzing the dynamical properties of its topological representatives. To facilitate this analysis we restrict our attention to topological representatives with special properties. In this subsection we recall some basic definitions and results from [BH1]. In Section 5 we extend these ideas to meet our current needs.

A *turn in $G$* is an unordered pair of oriented edges of $G$ originating at a common vertex. A turn is *nondegenerate* if it is defined by distinct oriented edges, and is *degenerate* otherwise. A *turn $(E_1, E_2)$ is contained in the filtration element $G_r$ (respectively the stratum $H_r$)* if both $E_1$ and $E_2$ are contained in $G_r$ (respectively $H_r$). If $E_1' E_2' \cdots E_{k-1} E_k'$ is the edge path associated to a path $\alpha$, then we say that $\alpha$ *contains the turns $(\bar{E}_i, E_{i+1})$ for $0 \le i \le k-1$*. This is consistent with our identification of a path with its associated edge path. Similarly, we say that $\alpha$ *crosses or contains each edge* that occurs in its associated edge path and we say that $\alpha$ *is contained in a subgraph $K$*, written $\alpha \subset K$, if each edge in its edge path is contained in $K$.

If $f : G \to G$ is a topological representative and $E$ is an edge of $G$, then we define $Tf(E)$ to be the first edge in (the edge path associated to) $f(E)$; for each turn $(E_i, E_j)$, define $Tf((E_i, E_j)) = (Tf(E_i), Tf(E_j))$. An important observation is that if $\alpha$ is a path and if the $Tf$-image of each turn in $\alpha$ is nondegenerate, then $f(\alpha)$ is a path.

Since $Tf$ sends edges to edges and turns to turns, it makes sense to iterate $Tf$. We say that a turn is *illegal* with respect to $f : G \to G$ if its image under some iterate of $Tf$ is degenerate; a turn is *legal* if it is not illegal. We say that a *path $\alpha \subset G$ is legal* if it contains only legal turns and that it is *$r$-legal* if it is contained in $G_r$ and all of its illegal turns are contained in $G_{r-1}$.

To each stratum $H_r$, we associate a square matrix $M_r$ called the *transition submatrix* for $H_r$; the $ij^{\mathrm{th}}$ entry of $M_r$ is the number of times that the $f$-image of the $j^{\mathrm{th}}$ edge crosses the $i^{\mathrm{th}}$ edge in either direction. A nonnegative matrix $M$ is *irreducible* if for each $i$ and $j$ there exists $n > 0$ so that the $ij^{\mathrm{th}}$ entry of $M^n$ is positive. By enlarging the filtration if necessary, we may assume that each $M_r$ is either the zero matrix or is irreducible. This gives us three kinds of strata. If $M_r$ is the zero matrix, then $H_r$ is a *zero stratum* . (These arise in the 'core subdivision' operation of [BH1].) If $M_r$ is irreducible, then it has an associated Perron-Frobenius eigenvalue $\lambda_r \geq 1$ [Sen]. If $\lambda_r > 1$, then we say that $H_r$ is an *exponentially-growing stratum*; if $\lambda_r = 1$, then we say that $H_r$ is a *non-exponentially-growing stratum*.

A topological representative $f : G \to G$ of $\mathcal{O}$ is a *relative train track map* with respect to the filtration $\phi = G_0 \subset G_1 \cdots \subset G_m = G$ if $G$ has no valence one vertices, if each nonzero $M_r$ is irreducible and if each exponentially-growing stratum satisfies the following conditions.

1. If $E$ is an edge in $H_r$, then $Tf(E)$ is an edge in $H_r$.

2. If $\beta \subset G_{r-1}$ is a nontrivial path with endpoints in $G_{r-1} \cap H_r$, then $f_\#(\beta)$ is nontrivial.

3. If $\sigma \subset H_r$ is a legal path, then $f(\sigma) \subset G_r$ is an $r$-legal path.

Complete details about relative train track maps can be found in [BH1].

The most important consequence of being a relative train track map is Lemma 5.8 of [BH1]. We repeat it here for the reader's convenience. A key point is that no cancellation of edges in $H_r$ occurs when the image $f^k(\sigma)$ of an $r$-legal path $\sigma \subset G_r$ is tightened to $f^k_\#(\sigma)$.

LEMMA 2.5.1.   *Suppose that $f : G \to G$ is a relative train track map, that $H_r$ is an exponentially-growing stratum and that $\sigma = a_1 b_1 a_2 \ldots b_l$ is a decomposition of an $r$-legal path into subpaths where each $a_i \subset H_r$ and each $b_j \subset G_{r-1}$. (Allow the possibility that $a_1$ or $b_l$ is trivial, but assume that the other subpaths are nontrivial.) Then $f_\#(\sigma) = f(a_1)f_\#(b_1)f(a_2)\ldots f_\#(b_l)$ and $f_\#(\sigma)$ is $r$-legal.*

2.6. *Free factor systems.*  Many of the arguments in this paper proceed by induction up through a filtration. In this subsection we consider filtrations from a group theoretic point of view and we show how to choose relative train track maps in which the steps between filtration elements are as small as possible.

We begin with the main geometric example.

*Example* 2.6.1. Suppose that $G$ is a marked graph and that $K$ is a sub-graph whose non-contractible components are labeled $C_1, \ldots, C_l$. Choose vertices $v_i \in C_i$ and a maximal tree $T \subset G$ such that each $T \cap C_i$ is a maximal tree in $C_i$. The tree $T$ determines inclusions $\pi_1(C_i, v_i) \to \pi_1(G, v)$. Let $F^i \subset F_n$ be the free factor of $F_n$ determined by $\pi_1(C_i, v_i)$ under the identification of $\pi_1(G, v)$ with $F_n$. Then $F^1 * F^2 * \cdots * F^l$ is a free factor of $F_n$. Without a specific choice of $T$, the $C_i$'s only determine the $F_i$'s up to conjugacy.

We reserve the notation $F^i$ for free factors of $F_n$. We use superscripts for the index so as to distinguish the index from the rank. The *conjugacy class* of $F^i$ is denoted $[[F^i]]$. If $F^1 * F^2 \cdots * F^k$ is a free factor and each $F^i$ is nontrivial (and so has positive rank), then we say that the collection $\mathcal{F} = \{[[F^1]], \ldots, [[F^k]]\}$ is a *nontrivial free factor system*. We refer to $\emptyset$ as *the trivial free factor system*.

Returning to Example 2.6.1, we write $\mathcal{F}(K)$ for the free factor system $\{[[\pi_1(C_1)]], \ldots, [[\pi_1(C_l)]]\}$ and say that $\mathcal{F}(K)$ *is realized by* $K$.

We define the *complexity of the free factor system* $\mathcal{F} = \{[[F^1]], \ldots, [[F^k]]\}$, written $\mathrm{cx}(\mathcal{F})$, to be 0 if $\mathcal{F}$ is trivial and to be the non-increasing sequence of positive integers that is obtained by rearranging the elements of $\{\mathrm{rank}(F^1), \ldots, \mathrm{rank}(F^k)\}$ if $\mathcal{F}$ is nontrivial. For any fixed $F_n$, there are only finitely many such complexities and we order them lexicographically. Thus $5, 3, 3, 1 > 4, 4, 4, 4, 4 > 4 > 0$; $\{[[F_n]]\}$ has the highest complexity and $\emptyset$ has the smallest.

The intersection of free factors is a free factor. More generally, we have the following result (Subgroup Theorem 3.14 of [SW]).

LEMMA 2.6.2. *Suppose that* $F_n = F^1 * F^2 \cdots * F^k$, *that* $H$ *is a subgroup of* $F_n$ *and that* $H(1), \ldots, H(l)$ *are the nontrivial subgroups of the form of* $H \cap (F^j)^c$ *for* $c \in F_n$. *Then* $H(1) * \cdots * H(l)$ *is a free factor of* $H$.

For any free factor systems $\mathcal{F}_1$ and $\mathcal{F}_2$, define $\mathcal{F}_1 \wedge \mathcal{F}_2$ to be the set of nontrivial elements of $\{[[F^i \cap (F^j)^c]] : [[F_i]] \in \mathcal{F}_1; [[F_j]] \in \mathcal{F}_2; c \in F_n\}$. Lemma 2.6.2 implies that $\mathcal{F}_1 \wedge \mathcal{F}_2$ is a (possibly empty) free factor system.

LEMMA 2.6.3. *If* $\mathcal{F}_1 \wedge \mathcal{F}_2 \neq \mathcal{F}_1$, *then* $\mathrm{cx}(\mathcal{F}_1 \wedge \mathcal{F}_2) < \mathrm{cx}(\mathcal{F}_1)$.

*Proof.* Each nontrivial $F^i \cap (F^j)^c$ is a free factor of $F^i$ and so either equals $F^i$ or has strictly smaller rank than $F^i$. Thus the set of ranks that occur for elements of $\mathcal{F}_1 \wedge \mathcal{F}_2$ is obtained from the set of ranks that occur for elements of $\mathcal{F}_1$ by (perhaps more than once) replacing a positive integer with a finite collection of strictly smaller integers. ∎

An outer automorphism $\mathcal{O}$ induces an action on the set of conjugacy classes of free factors. If $F^i$ is a free factor and $[[F^i]]$ is fixed by $\mathcal{O}$, then we say that $[[F^i]]$ *is* $\mathcal{O}$-*invariant*. Sometimes, we will abuse notation and say that $F^i$ is $\mathcal{O}$-invariant when we really mean that its conjugacy class is. We say that $\mathcal{F}$ *is* $\mathcal{O}$-*invariant* if each $[[F^i]] \in \mathcal{F}$ is $\mathcal{O}$-invariant. If $[[F^i]]$ is $\mathcal{O}$-invariant, then there is an automorphism $\Phi$ representing $\mathcal{O}$ such that $\Phi(F^i) = F^i$. Since $\Phi$ is well-defined up to composition with an inner automorphism determined by an element of $F^i$, $\Phi$ determines an outer automorphism of $F^i$ that we refer to as *the restriction of* $\mathcal{O}$ *to* $F^i$. Note that if $\mathcal{F}(K)$ is realized by $K$ and if $f : G \to G$ is a topological representative of $\mathcal{O}$ that setwise fixes each non-contractible component $C_i$ of $K$, then $\mathcal{F}(K)$ is $\mathcal{O}$-invariant.

We say that $\beta \in \mathcal{B}$ is *carried by* $[[F^i]]$ if it is in the closure of the circuits in $\mathcal{B}$ determined by conjugacy classes in $F_n$ of elements of $F^i$. It is an immediate consequence of the definitions that if $G$ is a marked graph and $K$ is a connected subgraph such that $[[\pi_1(K)]] = [[F^i]]$, then $\beta$ is carried by $[[F^i]]$ if and only if the realization of $\beta$ in $G$ is contained in $K$. A subset $B \subset \mathcal{B}$ is *carried* by the free factor system $\mathcal{F}$ if each element of $B$ is carried by an element of $\mathcal{F}$.

LEMMA 2.6.4.  *If $\beta \in \mathcal{B}$ is carried by both $[[F^1]]$ and $[[F^2]]$ then $\beta$ is carried by $[[F^1 \cap (F^2)^c]]$ for some $c \in F_n$.*

*Proof.* For $i = 1, 2$, choose a marked graph $G_i$ with one vertex $v_i$ and a subgraph $K_i$ so that the marking identifies $\pi_1(K_i, v_i)$ with $F^i$. Choose a homotopy equivalence $h : G_1 \to G_2$ that induces (via the markings on $G_1$ and $G_2$) the identity on $F_n$. Let $\beta_1 \subset K_1 \subset G_1$ and $\beta_2 = h_\#(\beta_1) \subset K_2 \subset G_2$ be bi-infinite paths that realize $\beta$. Part 2 of Lemma 2.3.1 implies that for each subpath $\sigma_k$ of $\beta_1$, $h_\#(\sigma_k) = c_k \tau_k d_k$ where $\tau_k \subset \beta_2 \subset K_2$ and $c_k$ and $d_k$ have uniformly bounded length. We may choose the $\sigma_k$'s to be an increasing collection whose union covers $\beta_1$ and so that $c_k = c$ and $d_k = d$ are independent of $k$. The union of the $\tau_k$'s covers $\beta_2$. Let $w_k = \sigma_k \bar{\sigma}_1$ and note that $h_\#(w_k) = [c\tau_k \bar{\tau}_1 \bar{c}]$ contains all but a uniformly bounded amount of $\tau_k$ as a subpath. The lemma now follows from the fact that the element of $F_n$ determined by both $w_k$ and $h_\#(w_k)$ is contained in $F^1 \cap (F^2)^c$. ∎

COROLLARY 2.6.5.  *For any subset $B \subset \mathcal{B}$ there is a unique free factor system $\mathcal{F}(B)$ of minimal complexity that carries every element of $B$. If $B$ has a single element, then $\mathcal{F}(B)$ has a single element.*

*Proof.* Since $[[F_n]]$ carries every element of $B$, there is at least one free factor system $\mathcal{F}_1$ of minimal complexity that carries every element of $B$. Suppose that $\mathcal{F}_2$ also carries every element of $B$ and that $\mathrm{cx}(\mathcal{F}_1) = \mathrm{cx}(\mathcal{F}_2)$. Lemma 2.6.4 implies that $\mathcal{F}_1 \wedge \mathcal{F}_2$ carries every element of $B$. Minimality and Lemma 2.6.3 therefore imply that $\mathcal{F}_1 = \mathcal{F}_2$. This proves that $\mathcal{F}(B)$ is well-defined.

Every element of $B$ is carried by some element of $\mathcal{F}(B)$. If $B$ has only one element but $\mathcal{F}(\mathcal{B})$ has more than one element, then we can reduce $\mathrm{cx}(\mathcal{F}(B))$ by reducing the number of elements in $\mathcal{F}(\mathcal{B})$. This proves the second part of the corollary. ∎

We write $[[F^1]] \sqsubset [[F^2]]$ if $F^1$ is conjugate to a free factor of $F^2$ and write $\mathcal{F}_1 \sqsubset \mathcal{F}_2$ if for each $[[F^i]] \in \mathcal{F}_1$ there exists a (necessarily unique) $[[F^j]] \in \mathcal{F}_2$ such that $[[F^i]] \sqsubset [[F^j]]$. The reader will easily check that if $K_1 \subset K_2$ are subgraphs of $G$, then $\mathcal{F}(K_1) \sqsubset \mathcal{F}(K_2)$.

In many of our induction arguments, it is important that the step between one filtration element and the next be as small as possible. This, and the fact that we sometimes replace $f : G \to G$ by an iterate, motivates the following definition and lemma.

*Definition* 2.6.6. A topological representative $f : G \to G$ and filtration $\emptyset = G_0 \subset G_1 \subset \cdots \subset G_K = G$ are *reduced* if each stratum $H_r$ has the following property : If a free factor system $\mathcal{F}'$ is invariant under the action of an iterate of $\mathcal{O}$ and satisfies $\mathcal{F}(G_{r-1}) \sqsubset \mathcal{F}' \sqsubset \mathcal{F}(G_r)$, then either $\mathcal{F}' = \mathcal{F}(G_{r-1})$ or $\mathcal{F}' = \mathcal{F}(G_r)$.

LEMMA 2.6.7. *For any $\mathcal{O}$-invariant free factor system $\mathcal{F}$, there exists a relative train track map $f : G \to G$ representing $\mathcal{O}$ and filtration $\emptyset = G_0 \subset G_1 \subset \cdots \subset G_K = G$ such that*:

- $\mathcal{F} = \mathcal{F}(G_r)$ *for some filtration element $G_r$.*

- *If $C$ is a contractible component of some $G_i$, then $f^j(C) \subset G_{i-1}$ for some $j > 0$.*

*If $\mathcal{O}$ is replaced by an iterate $\mathcal{O}^s$ then $f : G \to G$ may be chosen to be reduced.*

*Proof.* The first step in the proof is to show that for any nested sequence $\mathcal{F}_1 \sqsubset \cdots \sqsubset \mathcal{F}_l = \{[[F_n]]\}$ of $\mathcal{O}$-invariant free factor systems, there is a topological representative $f : G \to G$ of $\mathcal{O}$ and a filtration $\emptyset \subset G_1 \subset \cdots G_l = G$ so that each $\mathcal{F}_i$ is realized by $G_i$. The construction of $f : G \to G$ is very similar to the one in Lemma 1.16 of [BH1].

We argue by induction on $l$, the $l = 1$ case following from the fact that every $\mathcal{O}$ is represented by a homotopy equivalence of $R_n$. Let $\mathcal{F}_{l-1} = \{[[F^1]], \ldots, [[F^k]]\}$. Choose automorphisms $\Phi_i : F_n \to F_n$ representing $\mathcal{O}$ such that $\Phi_i(F^i) = F^i$.

For $1 \leq i \leq k$ and $1 \leq j \leq l-2$, define $\mathcal{F}_j^i = \{[[F^i]]\} \wedge \mathcal{F}_j$. Equivalently $\mathcal{F}_j^i$ consists of those elements of $\mathcal{F}_j$ that are contained (in the sense of $\sqsubset$) in $[[\mathcal{F}^i]]$. Then $\mathcal{F}_1^i \sqsubset \cdots \sqsubset \mathcal{F}_{l-2}^i \sqsubset \{[[F^i]]\}$ and $\mathcal{F}_j = \cup_{i=1}^k \mathcal{F}_j^i$. By induction on $l$, there are topological representatives $f_i : K^i \to K^i$ of the restriction of $\mathcal{O}$ to $F^i$ and there are filtrations $\emptyset = K_0^i \subset K_1^i \subset \cdots K_{l-1}^i = K^i$ so that each $\mathcal{F}_j^i$ is realized by $K_j^i$. We may assume inductively that $f_i$ fixes a vertex $v_i$ of $K^i$ and that the marking on $K^i$ identifies $F^i$ with $\pi_1(K^i, v_i)$ and identifies $\Phi_i$ with the automorphism $(f_i)_\# : \pi_1(K^i, v_i) \to \pi_1(K^i, v_i)$.

Let $F^{k+1} \cong F_{n_{k+1}}$ be a free factor such that $F^1 * \cdots * F^{k+1} \cong F_n$. Define $G$ to be the graph obtained from the disjoint union of the $K^i$'s by adding edges $E_i$, $2 \leq i \leq k$, connecting $v_1$ to $v_i$, and by adding $n_{k+1}$ loops $\{L_j\}$ based at $v_1$. Collapsing the $E_i$'s to $v_1$ gives a homotopy equivalence of $(G, v_1)$ onto a graph $(G', v')$ whose fundamental group is naturally identified with $F^1 * \cdots * F^{k+1} \cong F_n$. This provides a marking on $G$.

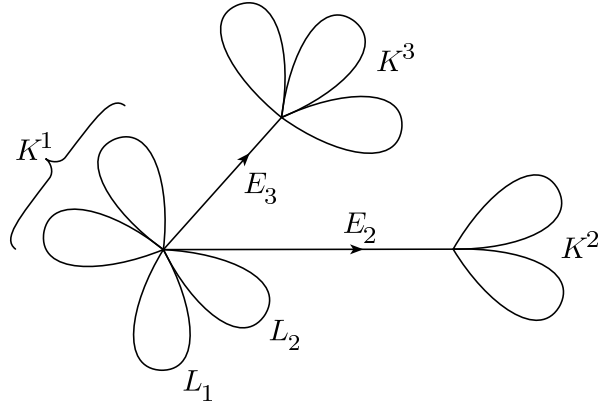

The filtration $\phi = G_0 \subset G_1 \subset \cdots \subset G_l = G$ is defined by $G_j = \cup_{i=1}^k K_j^i$. It is immediate from the definitions that $\mathcal{F}_j = \mathcal{F}(G_j)$.

There exists $c_i \in F_n$, $2 \le i \le k$, such that $\Phi_1(x) = c_i \Phi_i(x) \bar{c}_i$ for all $x \in F_n$. Let $\gamma_i \subset G$ be the loops based at $v_1$ that are identified, under the marking, with $c_i$. Extend $\cup f_i : \cup K_i \to \cup K_i$ to a topological representative $f : G \to G$ by defining $f(E_i) = \gamma_i E_i$ and by defining $f(L_j)$ according to $\Phi_1$. Then $f_\# : \pi_1(G, v_1) \to \pi_1(G, v_1)$ induces $\Phi_1 : F_n \to F_n$ and so represents $\mathcal{O}$. This completes the first part of the proof.

The second step in the proof is to promote $f : G \to G$ to a relative train track map. This may cause the filtration to be expanded but we will maintain the property that each $\mathcal{F}_j$ is realized by some filtration element. Applying this with $\mathcal{F}_1 = \mathcal{F}$ and $l = 2$ will complete the proof of the first part of the lemma.

In Section 5 of [BH1], there is an algorithm that begins with an arbitrary topological representative of $\mathcal{O}$ and filtration and produces a relative train track map and filtration that represents $\mathcal{O}$. The algorithm uses only the operations of : subdivision; folding; tightening; collapsing pre-trivial forests; valence-one homotopy; and (with some restrictions) valence-two homotopy. (See [BH1 for definitions.) It suffices to work one step at a time and show that if $\hat{f} : \hat{G} \to \hat{G}$ is obtained from $f : G \to G$ by performing one of these operations and if each $\mathcal{F}_i$ is realized by a filtration element of $G$, then each $\mathcal{F}_i$ is realized by a filtration element of $\hat{G}$.

Let $p : G \to \hat{G}$ be the natural homotopy equivalence. We will show that if $C_1$ and $C_2$ are non-contractible components of a filtration element $G_k$ in $G$, then $\hat{C}_1 = p(C_1)$ and $\hat{C}_2 = p(C_2)$ are disjoint $\hat{f}$-invariant subgraphs and $[[\pi_1(C_i)]] = [[\pi_1(\hat{C}_i)]]$. Since $p(G_k)$ is a filtration element in the induced filtration on $\hat{G}$, this will complete the second step.

We consider first the case that $\hat{G}$ is obtained from $G$ by collapsing a pre-trivial forest $X$. (A forest is pre-trivial if its image under some iterate of $f$ is a finite union of points.) Since a component of $X$ cannot intersect both $C_1$ and $C_2$, $\hat{C}_1$ and $\hat{C}_2$ are disjoint. Let $X_0$ be a component of $X$. If $C_i \cap X_0$ has more than one component then there is a circuit that is not contained in $C_i$ but whose image under some iterate of $f$ is contained in $C_i$. This contradicts the fact that $f$ and $f|C_i$ are homotopy equivalences. Thus $C_i \cap X$ is connected and $[[\pi_1(C_i)]] = [[\pi_1(\hat{C}_i)]]$.

The operations of subdivision, folding, tightening and valence-one homotopy are straightforward to check and we leave this to the reader. Suppose that $v$ is a valence-two vertex with incident edges $E_r$ and $E_s$. If $E_r$ and $E_s$ are disjoint from $C_1 \cup C_2$ then the valence-two homotopy does not affect $C_1$ or $C_2$. If $E_r$ and $E_s$ are both contained in $C_1 \cup C_2$, then they are both contained in $C_1$ or both contained in $C_2$, say $C_1$. In this case, $C_2$ is unchanged by the valence-two homotopy and it follows immediately from the definitions that $\hat{C}_1 \cap \hat{C}_2 = \emptyset$ and that $[[\pi_1(C_1)]] = [[\pi_1(\hat{C}_1)]]$. If one edge, say $E_r$, is contained in $C_1 \cup C_2$ but $E_s$ is not, then $E_s$ is not contained in $G_k$. In particular, $E_s$ belongs to a higher stratum than does $E_r$ so the valence-two homotopy is performed by sliding $v$ across $E_r$. From the point of view of $G_k$, a valence-one homotopy is being performed on one of its components. The desired properties now follow immediately from the definitions. This completes the second step.

If $C$ is a contractible component of $G_i$ and no $f^j(C)$ is contained in $G_{i-1}$, then there is a collection of components of $G_i$ that form an invariant forest. Collapse each of these components to points. After tightening the images of the remaining edges and possibly collapsing edges in other zero strata to points, there is a quotient map $q : G \to G'$, an induced topological representative $f' : G' \to G'$ and a filtration with elements of the form $q(G_j)$. By construction, each $\mathcal{F}(G_j) = \mathcal{F}(q(G_j))$ so it is still true that each $\mathcal{F}_j$ is realized by some filtration element. We will show that $f' : G' \to G'$ is a relative train track map. Thus, after repeating this operation finitely many times, we establish the first statement in the lemma.

Assume that the edges in $G'$ have the same labels as they did in $G$. The key point is that for any path $\sigma \subset G$, the edge path associated to $\sigma' = q_\#(\sigma)$ is obtained from the edge path associated to $\sigma$ by removing all occurrences of the collapsed edges. In particular, if $E$ is an edge of $G$ that does not collapse, then $f'(E)$ is obtained from $f(E)$ by removing occurrences of collapsed edges. It is straightforward to check that $q$ induces a one-to-one correspondence between the exponentially-growing (respectively non-exponentially-growing) strata of $f : G \to G$ and the exponentially-growing (respectively non-exponentially-growing) strata of $f' : G' \to G'$. If $H_s$ is exponentially growing, then conditions 1 and 3 in the definition of relative train track map for $q(H_s)$ follow immediately from conditions 1 and 3 for $H_s$. Condition 2 for $q(H_s)$ follows from condition 2 for $H_s$ and the observation that if $\sigma \subset G$ is nontrivial, then $\sigma'$ is trivial if and only if some $f^k(\sigma)$ is entirely contained in a component of $G_i$ that is collapsed by $q$. In particular if $\sigma'$ and each $f^k_\#(\sigma)$ are nontrivial, then each $(f')^k_\#(\sigma)$ is nontrivial. The remaining details are left to the reader.

For the last statement of the lemma, extend $\mathcal{F}$ to a maximal (with respect to $\sqsubset$) nested sequence $\mathcal{C}$ of distinct free factor systems $\mathcal{F}_i$ such that each $\mathcal{F}_i$ is invariant under the action of some iterate of $\mathcal{O}$. Choose $s > 0$ so that each $\mathcal{F}_i$ is $\mathcal{O}^s$-invariant. Let $f : G \to G$ be a relative train track map representing $\mathcal{O}^s$ and let $\emptyset = G_0 \subset G_1 \subset \cdots \subset G_m$ be a filtration for $f : G \to G$ such that each $\mathcal{F}_i$ is realized by a filtration element. In other words, $\mathcal{C}$ is a nested subsequence of $\emptyset \sqsubset \mathcal{F}(G_1) \sqsubset \cdots \sqsubset \mathcal{F}(G_m) = \{[[F_n]]\}$. Since $\mathcal{C}$ is maximal, $f : G \to G$ is reduced. $\qquad\blacksquare$

# 3. Attracting laminations

3.1. *Attracting laminations associated to exponentially-growing strata.* Measured foliations play a central role in Thurston's classification of mapping classes. When working in $\mathrm{Out}(F_n)$, one's method of attack is determined, to a great extent, by how one chooses to generalize measured foliations. In this paper, we adopt a Nielsen-like point of view similar to that of [HT]. In particular, we work with laminations rather than foliations, we make extensive use of the space of ends (in this case $\partial F_n$) and we restrict our considerations to topological rather than measure theoretic properties.

An important feature of our approach is that we work directly in $\mathcal{B}$. Thus an attracting lamination $\Lambda^+$, defined below, is a closed set of lines and not a single point in a space of (measured) laminations. This makes certain arguments longer and perhaps less transparent but it has an essential advantage: it enlarges the basin of attraction for the attracting laminations (see Theorem 6.0.1 and Remark 6.0.2). This is crucial for our application of the 'Tits ping-pong' argument via Corollary 3.4.3.

In this subsection we define attracting laminations in terms of the action of an outer automorphism $\mathcal{O}$ on the space of lines $\mathcal{B}$ and then begin to develop their properties in terms of relative train track maps $f : G \to G$ that represent $\mathcal{O}$. A different approach to laminations may be found in [Lus].

*Definitions* 3.1.1. We say that $\beta' \in \mathcal{B}$ is *weakly attracted to* $\beta \in \mathcal{B}$ *under the action of* $\mathcal{O}$ if $\mathcal{O}^k_\#(\beta') \to \beta$. (We describe the attraction as weak to emphasize that we are working in a non-Hausdorff space that ignores measure.) A subset $U \subset \mathcal{B}$ is an *attracting neighborhood of* $\beta \in \mathcal{B}$ *for the action of* $\mathcal{O}$ if $\mathcal{O}_\#(U) \subset U$ and if $\{\mathcal{O}^k_\#(U) : k \geq 0\}$ is a neighborhood basis for $\beta$ in $\mathcal{B}$. If $U$ is an attracting neighborhood of $\beta$ for the action of $\mathcal{O}$, then $\beta'$ is weakly attracted to $\beta$ under the action of $\mathcal{O}$ if and only if some $\mathcal{O}^k_\#(\beta') \in U$. The reader will easily check that if $\alpha \subset G$ realizes $\beta'$ and $\gamma \subset G$ realizes $\beta$, then $\beta'$ is weakly attracted to $\beta$ if and only if each subpath of $\gamma$ is contained in $f^k_\#(\alpha)$ for all sufficiently large $k$. We sometimes say that $\gamma$ *is weakly attracted to* $\alpha$ *under the action of* $f_\#$.

*Remark* 3.1.2. Since $\mathcal{B}$ is not Hausdorff, limits are not unique. For example, working by analogy in the mapping class group of a closed surface $S$, suppose that $\phi : S \to S$ is a homeomorphism in Thurston normal form and that there are two pseudo-Anosov components $S_1$ and $S_2$. Denote the stable and unstable laminations of $\phi | S_i$ by $\Lambda_i^-$ and $\Lambda_i^+$ respectively; let $\lambda_i^+$ be a leaf of $\Lambda_i^+$. Every complete geodesic $\gamma \subset S$ that intersects the interior of $S_i$ and is not a leaf of $\Lambda_i^-$ is weakly attracted to $\lambda_i^+$. Thus most complete geodesics in $S$ are weakly attracted to both $\lambda_1^+$ and $\lambda_2^+$.

*Definitions* 3.1.3. A bi-infinite path $\sigma$ in a marked graph $G$ is *birecurrent* if every finite subpath of $\sigma$ occurs infinitely often as an unoriented subpath of each end of $\sigma$. A line in $G$ is birecurrent if the path representing it, with either choice of orientation, is birecurrent. An element of $\mathcal{B}$ is birecurrent if some, and hence any (see Lemma 3.1.4), realization in a marked graph is birecurrent.

LEMMA 3.1.4. *If some realization of $\beta \in \mathcal{B}$ in a marked graph is birecurrent then every realization of $\beta \in \mathcal{B}$ in a marked graph is birecurrent. If $\beta$ is birecurrent, then $\Psi_\#(\beta)$ is birecurrent for all $\Psi \in \mathrm{Out}(F_n)$.*

*Proof.* Suppose that $\sigma \subset G$ and $\sigma' \subset G'$ are realizations of $\beta$ and that $\sigma$ is birecurrent. Let $h : G \to G'$ be a homotopy equivalence that respects the markings and let $C$ be the bounded cancellation constant for $h$. Choose lifts $\tilde{\sigma} \subset \tilde{G}$, $\tilde{\sigma}' \subset \tilde{G}'$ and $\tilde{h} : \tilde{G} \to \tilde{G}'$ such that $\tilde{h}_\#(\tilde{\sigma}) = \tilde{\sigma}'$.

Given a finite subpath $\tilde{\sigma}'_0 \subset \tilde{\sigma}'$, extend $\tilde{\sigma}'_0$ to $\tilde{\tau}' \subset \tilde{\sigma}'$ by adding $C$ initial and terminal edges. Choose a finite subpath $\tilde{\tau} \subset \tilde{\sigma}$ such that $\tilde{h}_\#(\tilde{\tau}) \supset \tilde{\tau}'$. Since $\sigma$ is birecurrent, each end of $\tilde{\sigma}$ contains infinitely many copies $\tilde{\tau}_i$ of $\tilde{\tau}$. Define $\tilde{\mu}'_i$ by removing $C$ initial and terminal edges from $\tilde{h}_\#(\tilde{\tau}_i)$. Lemma 2.3.1 implies that $\tilde{\mu}'_i \subset \tilde{\sigma}'$. By construction, each $\tilde{\mu}'_i$ contains a copy of $\tilde{\sigma}'_0$ and we have verified that $\sigma'$ is birecurrent.

For the second part of the lemma, suppose that $\beta$ is realized by $\sigma \subset G$ and that $\sigma$ is birecurrent. Choose a topological representative $f : G \to G$ of $\Psi$. The preceding argument, with $h : G \to G'$ replaced by $f : G \to G$, carries over with no other changes to prove that $f_\#(\sigma)$ is birecurrent. Since $f_\#(\sigma)$ realizes $\Psi_\#(\beta)$, this completes the proof. ◻

*Definitions* 3.1.5. A closed subset $\Lambda^+$ of $\mathcal{B}$ is an *attracting lamination for* $\mathcal{O}$ if it is the closure of a single point $\beta$ that:

(1) is birecurrent;

(2) has an attracting neighborhood for the action of some iterate of $\mathcal{O}$;

(3) is not carried by an $\mathcal{O}$-periodic free factor of rank one.

$\beta$ is said to be *generic for* $\Lambda^+$ or simply $\Lambda^+$*-generic*. We denote the *set of attracting laminations* for $\mathcal{O}$ by $\mathcal{L}(\mathcal{O})$.

LEMMA 3.1.16.   $\mathcal{L}(\mathcal{O})$ *is* $\mathcal{O}$-*invariant.*

*Proof.* Suppose that $\beta$ is generic with respect to $\Lambda^+ \in \mathcal{L}(\mathcal{O})$. Lemma 3.1.4 implies that $\mathcal{O}_{\#}(\beta)$ is birecurrent. If $V$ is an attracting neighborhood for $\beta$ under the action of $\mathcal{O}^s$ then $\mathcal{O}_{\#}(V)$ is an attracting neighborhood for $\mathcal{O}_{\#}(\beta)$ under the action of $\mathcal{O}^s$. If $[[F]]$ is an $\mathcal{O}$-periodic, rank-one free factor that carries $\mathcal{O}_{\#}(\beta)$ and if $\Phi : F_n \to F_n$ represents $\mathcal{O}^{-1}$, then $[[\Phi(F)]]$ is an $\mathcal{O}$-periodic rank-one free factor that carries $\beta$. Thus $\mathcal{O}(\beta)$ is generic with respect to $\mathcal{O}_{\#}(\Lambda^+) \in \mathcal{L}(\mathcal{O})$.                                                                          □

In order to analyze $\mathcal{L}(\mathcal{O})$, we bring relative train track maps into the discussion.

*Definitions* 3.1.7. Assume that $f : G \to G$ and $\emptyset = G_0 \subset G_1 \subset \cdots \subset G_K = G$ are a relative train track map and filtration representing $\mathcal{O}$ and that $H_r$ is an exponentially-growing stratum. For each edge $E$ of $H_r$ and $k \geq 0$, define the *k-tile determined by* $E$ to be the unoriented path determined by $f_{\#}^k(E)$; i.e. $f_{\#}^k(E)$ with either of its orientations. A path in $G$ is called a *tile* if it is a $k$-tile for some $k$. A *k-tiling* of a path in $G_r$ is a decomposition into subpaths that are either $k$-tiles or are contained in $G_{r-1}$. A bi-infinite path $\lambda \subset G$ has an *exhaustion by tiles* if it can be written as the increasing union of tiles; equivalently $\lambda$ has an exhaustion by tiles if each of its subpaths occurs as a subpath of a tile in $\lambda$. If $\lambda$ has an exhaustion by tiles, then condition (3) in the definition of relative train track maps (subsection 2.5) implies that $\lambda \subset G_r$ is $r$-legal. We say that a line in $G$ has a $k$-tiling or has an exhaustion by tiles if the path representing it, with either choice of orientation, has this property.

A nonnegative matrix $M$ is *aperiodic* if it has an iterate $M^k$ that is positive; i.e. if each entry of $M^k$ is positive. Aperiodic matrices are irreducible but the converse is not true. See [Sen] for the precise relationship between aperiodic and irreducible matrices. If $f : G \to G$ is a relative train track map, then we say that *an exponentially-growing stratum* $H_r$ *is aperiodic* if the transition submatrix $M_r$ is aperiodic and that $f : G \to G$ is *eg-aperiodic* if each exponentially-growing stratum is aperiodic.

The following lemma records some elementary but useful observations.

LEMMA 3.1.8.   *Assume that* $H_r$ *is an exponentially-growing stratum.*

1. *Every path in* $G_r$ *has a 0-tiling.*

2. *If* $\sigma$ *is* $r$-legal *and has a* $k$-tiling, *then* $f_{\#}(\sigma)$ *has a* $(k+1)$-tiling.

3. *If* $k < l$, *then each* $l$-tile has a $k$-tiling; *if* $M_r^{k_0}$ *is positive and* $l - k \geq k_0$ *then each* $k$-tile occurs as a subpath of each $l$-tile.

4. *The $ij$<sup>th</sup> coefficient of $M_r^k$ is the number of times that the $i$<sup>th</sup> edge in $H_r$ is crossed, in either direction, by the $k$-tile determined by the $j$<sup>th</sup> edge in $H_r$.*

5. *If $\lambda$ has an exhaustion by tiles then $f_\#(\lambda)$ has an exhaustion by tiles.*

*Proof.* (1) is immediate from the definitions. (2), (4) and (5) follow from Lemma 2.5.1. (3) follows from (1), (2), and (4).     □

The natural way to find attracting laminations is to look at weak limits of $\mathcal{O}^k(\alpha)$ for some circuit $\alpha$. When working with respect to a relative train track map $f : G \to G$, one can look at the weak limits of $f^k(E)$ for some edge $E$. We use this simple approach in the next pair of lemmas.

LEMMA 3.1.9.    *Suppose that $f : G \to G$ and $\emptyset = G_0 \subset G_1 \subset \cdots \subset G_K = G$ are a relative train track map and filtration representing $\mathcal{O}$ and that $H_r$ is an aperiodic exponentially-growing stratum. Then there is an attracting lamination $\Lambda^+$ with generic leaf $\beta$ so that $H_r$ is the highest stratum crossed by the realization $\lambda$ of $\beta$ in $G$.*

*Proof.* Choose an edge $E$ of $H_r$ and $m > 0$ so that $f_\#^m(E) = \alpha E \beta$ for some nontrivial paths $\alpha, \beta \subset G_r$. Let $h = f^m$; choose lifts $\tilde{E}, \tilde{\alpha}, \tilde{\beta}$ and $\tilde{h} : \Gamma \to \Gamma$ so that $\tilde{h}(\tilde{E}) = \tilde{\alpha}\tilde{E}\tilde{\beta}$. Define $\tilde{\tau}_j = \tilde{h}_\#^j(\tilde{E})$ and note that $\tilde{\tau}_j$ is a lift of a $jm$-tile. Then $\tilde{\tau}_0 = \tilde{E}$, $\tilde{\tau}_1 = \tilde{\alpha}\tilde{\tau}_0\tilde{\beta}$ and more generally $\tilde{\tau}_{j+1} = \tilde{\alpha}_j\tilde{\tau}_j\tilde{\beta}_j$, for nontrivial paths $\tilde{\alpha}_j$ and $\tilde{\beta}_j$. The $\tilde{\tau}_j$'s are therefore an increasing sequence of lifts of tiles whose union is a bi-infinite path $\tilde{\lambda} \subset \Gamma$ that is fixed by $\tilde{h}_\#$. We claim that the projection $\lambda \subset G$ realizes an element $\beta \in \mathcal{B}$ that is generic with respect to some element of $\mathcal{L}(\mathcal{O})$.

Since $\tilde{E}$ is mapped over itself by $\tilde{h}$, there is a point $\tilde{x} \in \tilde{E}$ that is fixed by $\tilde{h}$. After replacing $m$ by a multiple if necessary, we may assume that the $h_\#$-image of any edge in $H_r$ contains at least two $H_r$-edges. Define $\tilde{\lambda}_k$ to be the subpath of $\tilde{\lambda}$ that begins with the $k$<sup>th</sup> $\tilde{H}_r$-edge to the left of $\tilde{E}$ and ends with the $k$<sup>th</sup> $\tilde{H}_r$-edge to the right of $\tilde{E}$; define $V_k = N(\tilde{\lambda}_k)$. Lemma 2.5.1 implies that $\tilde{h}_\#(\tilde{\lambda}_k) \supset \tilde{\lambda}_{2k}$. By Lemma 2.3.1(3), $\tilde{h}_\#(V_k) \subset V_{k+1}$ for all sufficiently large $k$. The $V_k$'s are a neighborhood basis for $\lambda$ and so for all sufficiently large $k$, $V_k$ is an attracting neighborhood of $\lambda$ for the action of $\mathcal{O}^m$.

Since the difference between the number of edges in $\tilde{h}_\#(\tilde{\lambda}_k)$ and the number of edges in $\tilde{\lambda}_k$ is unbounded, the $\tilde{\lambda}_k$'s cannot be subpaths of a single $\tilde{h}_\#$-invariant axis. In other words, $\lambda$ is not a circuit and so cannot be carried by any free factor of rank one.

By construction, $\lambda$ has an exhaustion by tiles. We now use this to show that $\lambda$ has a $k$-tiling for all $k \geq 1$. A $k$-tiling of $\lambda$ corresponds to a subdivision of $\tilde{\lambda}$ and so is determined by the vertices of $\tilde{\lambda}$ that are the endpoints of the subdivision pieces. By Lemma 3.1.8(3), we may assume that each tile $\tau_i$ in an exhaustion of $\lambda$ has a $k$-tiling and so defines a finite set $\tilde{V}_i$ of vertices of $\tilde{\lambda}$. After passing to a subsequence, we may assume that a vertex $\tilde{v} \in \tilde{\lambda}$ satisfies either $\tilde{v} \in \tilde{V}_i$ for all large $i$ or $\tilde{v} \notin \tilde{V}_i$ for all large $i$. The set of vertices that satisfy the former condition determines a $k$-tiling of $\lambda$.

Lemma 2.5.1 implies that the first and last edges of any tile are contained in $H_r$. Thus each end of $\lambda$ must contain infinitely many edges in $H_r$. Lemma 3.1.8(3) and the existence of $k$-tilings for all $k$ imply that each tile occurs infinitely often in each end of $\lambda$. Since every finite subpath of $\lambda$ is contained in a tile, $\lambda$ is bicurrent. ☐

Having proved that $\mathcal{L}(\mathcal{O})$ is not empty, we next list some useful properties of generic leaves.

LEMMA 3.1.10.    *Assume that $\beta \in \mathcal{B}$ is a generic line of some $\Lambda^+ \in \mathcal{L}(\mathcal{O})$, that $f : G \to G$ and $\emptyset = G_0 \subset G_1 \subset \cdots \subset G_K = G$ are a relative train track map and filtration representing $\mathcal{O}$ and that $\lambda$ is the realization of $\beta$ in $G$. Then*:

(1) *The highest stratum $H_r$ crossed by $\lambda$ is exponentially growing.*

(2) *$\lambda$ is $r$-legal.*

*Assume that tiles are defined with respect to $H_r$*:

(3) *$\lambda$ has a $k$-tiling for all $k \geq 1$.*

(4) *$\lambda$ has an exhaustion by tiles.*

*Proof.* We argue by induction on the rank $n$ of $F_n$. The $n = 1$ case is vacuous so we may assume that the lemma holds for outer automorphisms of free groups with rank smaller than $n$.

As a first case suppose that $\lambda \subset G_{K-1}$. Let $m$ be the smallest positive integer so that the component $C$ of $G_{K-1}$ that contains $\lambda$ is $f^m$-invariant. The inductive hypothesis, applied to the restriction of $\mathcal{O}^m$ to $\pi_1(C)$ completes the proof.

We now assume that $\lambda$ contains edges of $H_K$. Necessarily $\lambda$ must cross some edges of $H_K$ infinitely many times. Choose $s > 0$ so that $\lambda$ has an attracting neighborhood for the action of $\mathcal{O}^s$ and let $\alpha \neq \lambda \subset G$ be a circuit that is weakly attracted to $\lambda$ under the action of $f_\#^s$. Since $f_\#^{sl}(\alpha)$ weakly converges to $\lambda$ as $l \to \infty$, $f_\#^s$ cannot act periodically on $\alpha$ with period different from one; since $\alpha \neq \gamma$, $f_\#^s$ cannot fix $\alpha$. Thus the number of edges in $f_\#^{sl}(\alpha)$ grows without bound. Since $f_\#^{sl}(\alpha)$ weakly converges to a line with infinitely many $H_K$-edges, the number of $H_K$-edges in $f_\#^{sl}(\alpha)$ grows without bound. It follows that $H_K$ must be an exponentially-growing stratum. This completes the proof of (1).

For (2), let $j$ be the number of illegal turns that $\alpha$ has in $H_r$. The number of illegal turns of $f_{\#}^{sl}(\alpha)$ in $H_r$ is bounded above by $j$. Choose a subpath $\lambda_0$ of $\lambda$. Since $\lambda$ is a weak limit of the $f_{\#}^{sl}(\alpha)$'s, $\lambda_0$ occurs as a subpath of the periodic line determined by $f_{\#}^{sl}(\alpha)$ for all large $l$. Since the length of the circuit $f_{\#}^{sl}(\alpha)$ increases without bound, $\lambda_0$ is covered by two fundamental domains of the line $f_{\#}^{si}(\alpha)$ for all large $l$. It follows that the number of illegal turns of $\lambda_0$ in $H_r$ is bounded above by $2j$. Birecurrence, and the fact that $\lambda_0$ was arbitrary, therefore imply that $\lambda$ is $r$-legal.

For (3), fix $k \geq 0$ and let $\tilde{\lambda} \subset \Gamma$ be a lift of of $\lambda$. A $k$-tiling of $\lambda$ corresponds to a subdivision of $\tilde{\lambda}$ and so is determined by the vertices of $\tilde{\lambda}$ that are the endpoints of the subdivision pieces.

Let $q$ be the number of edges in $\alpha$. For any finite subpath $\lambda_0 \subset \lambda$ let $\lambda_1 \subset \lambda$ be a finite subpath that contains $2q+1$ copies of $\lambda_0$. As in the previous case, if $l$ is sufficiently large, then $\lambda_1$ occurs as a subpath of the periodic line determined by $f_{\#}^{sl}(\alpha)$ that is covered by two fundamental domains. In particular at least one copy of $\lambda_0$ occurs as a subpath of $f_{\#}^{sl}(E)$ for some edge $E$ of $G_r$. We conclude that $\lambda$ is an increasing union of finite subpaths that have $k$-tilings. The $k$-tilings of these subpaths correspond to finite sets $\tilde{V}_i$ of vertices of $\tilde{\lambda}$. After passing to a subsequence, we may assume that a vertex $\tilde{v} \in \tilde{\lambda}$ satisfies either $\tilde{v} \in \tilde{V}_i$ for all large $i$ or $\tilde{v} \notin \tilde{V}_i$ for all large $i$. The set of vertices that satisfy the former condition determines a $k$-tiling of $\lambda$. This proves (3).

To prove (4), choose a finite subpath $\lambda_0 \subset \lambda$. Birecurrence implies that there is a finite subpath $\lambda_1 \subset \lambda$ that contains two disjoint copies of $\lambda_0$. After enlarging $\lambda_0$ if necessary we may assume that $\lambda_0$ contains at least one edge of $H_r$. By (3) $\lambda$ has a $k$-tiling where $k$ is so large that each $k$-tile contains more edges than $\lambda_1$ does. In any $k$-tiling of $\lambda$ there are at most two $k$-tiles that intersect $\lambda_1$; one of these must contain a copy of $\lambda_0$. We have now shown that every finite subpath of $\lambda$ is contained in a tile in $\lambda$. This completes the proof of (4).                                                                 □

COROLLARY 3.1.11.  *Assume that* $f : G \to G$ *and* $\emptyset = G_0 \subset G_1 \subset \cdots \subset G_K = G$ *are a relative train track map and filtration representing* $\mathcal{O}$, *that* $H_r$ *is an aperiodic exponentially-growing stratum and that tiles are defined with respect to* $H_r$. *Assume further that* $\beta \in \mathcal{B}$ *is* $\Lambda^+$-*generic for some* $\Lambda^+ \in \mathcal{L}(\mathcal{O})$ *and that* $H_r$ *is the highest stratum crossed by the realization of* $\beta$ *in* $G$. *Then* $\{N(\tau) : \tau$ *is a tile*$\}$ *is a neighborhood basis in* $\mathcal{B}$ *for* $\beta$. *In particular, all such* $\beta$ *have the same closure.*

*Proof.* Let $\lambda \subset G$ be the realization of $\beta$. Lemma 3.1.10(3) and Lemma 3.1.8(3) imply that $\lambda$ contains all tiles. Conversely, Lemma 3.1.10(4) implies that every subpath of $\lambda$ is contained in a tile in $\lambda$.                                      □

*Definitions* 3.1.12. Lemma 3.1.9, Lemma 3.1.10 and Corollary 3.1.11 imply that for any relative train track map representing $\mathcal{O}$ and for any aperiodic exponentially-growing stratum $H_r$ there is a unique element $\Lambda^+ \in \mathcal{L}(\mathcal{O})$ with the property that $H_r$ is the highest stratum crossed by the realization $\lambda \subset G$ of a $\Lambda^+$-generic line. We say that $H_r$ *is the stratum determined by* $\Lambda^+$ and that $\Lambda^+$ *is the attracting lamination associated to* $H_r$.

LEMMA 3.1.13. $\mathcal{L}(\mathcal{O})$ *is finite*.

*Proof.* Choose a relative train track map $f : G \to G$ and filtration $\emptyset = G_0 \subset G_1 \subset \cdots \subset G_K = G$ representing $\mathcal{O}$. If $f : G \to G$ is eg-aperiodic (i.e. if each exponentially-growing stratum $H_r$ is aperiodic) then there is a one-to-one correspondence between the exponentially-growing strata and the elements of $\mathcal{L}(\mathcal{O})$. Under these conditions, the lemma is clear.

On the other hand, suppose that some $M_r$ is not aperiodic. There is a partition ([Sen]) of the edges of $H_r$ into $s > 1$ sets $P_1, \ldots P_s$ such that for each edge $E \in P_i$, the edge path $f_\#(E)$ only crosses edges in $P_{i+1 (\mathrm{mod}\ s)}$ and in $G_{r-1}$. The matrix $M_r^s$ is not irreducible and so the filtration for $f$ must be enlarged to obtain a filtration for $f^s$. After replacing $f$ by $f^s$, $H_r$ divides into $s$ exponentially-growing strata. If $s$ is maximal, then the transition matrix for each of these $s$ exponentially-growing strata is aperiodic and irreducible. We have shown that some iterate $\mathcal{O}^p$ of $\mathcal{O}$ is represented by an eg-aperiodic relative train track map . Since $\mathcal{L}(\mathcal{O}^p) = \mathcal{L}(\mathcal{O})$, we are reduced to the previous case.  □

It is natural to focus on the case that each element of $\mathcal{L}(\mathcal{O})$ is $\mathcal{O}$-invariant. In the analogy with the mapping class group this corresponds to assuming that the pseudo-Anosov components are fixed rather than permuted. The following lemma relates this hypothesis to relative train track maps.

LEMMA 3.1.14. *The following are equivalent*:

(1) *Each element of* $\mathcal{L}(\mathcal{O})$ *is* $\mathcal{O}$-*invariant*.

(2) *Each element of* $\mathcal{L}(\mathcal{O})$ *has an attracting neighborhood for* $\mathcal{O}_\#$.

(3) *Every relative train track map* $f : G \to G$ *representing* $\mathcal{O}$ *is* eg-*aperiodic*.

(4) *Some relative train track map* $f : G \to G$ *representing* $\mathcal{O}$ *is* eg-*aperiodic*.

*Proof.* It is obvious that $(3) \Rightarrow (4)$.

Suppose that $f : G \to G$ is an eg-aperiodic relative train track map for $\mathcal{O}$, that $\Lambda^+ \in \mathcal{L}(\mathcal{O})$ and that $H_r$ is the exponentially-growing stratum associated to $\Lambda^+$. If $\lambda \subset G_r$ is $\Lambda^+$-generic, then $f_\#(\lambda)$ is $\mathcal{O}(\Lambda^+)$-generic. Since $H_r$ is the highest stratum crossed by $f_\#(\lambda)$, Corollary 3.1.11 implies that $\mathcal{O}(\Lambda^+) = \Lambda^+$. Thus $(4) \Rightarrow (1)$.

Suppose that $f : G \to G$ is a relative train track map representing $\mathcal{O}$, that $H_r$ is an exponentially-growing stratum and that $M_r$ is not aperiodic. There is a partition of the edges of $H_r$ into $s > 1$ sets $P_1, \ldots P_s$ such that for each edge $E \in P_i$, the edge path $f_\#(E)$ only crosses edges in $P_{i+1(\mathrm{mod}\ s)}$ and in $G_{r-1}$. The matrix $M_r^s$ is not irreducible and so the filtration for $f$ must be enlarged to obtain a filtration for $f^s$. After replacing $f$ by $f^s$, $H_r$ divides into $s$ exponentially-growing strata, one for each $P_i$. By Lemma 3.1.9 each of these contributes an element to $\mathcal{L}(\mathcal{O})$ that clearly does not have an attracting neighborhood for $\mathcal{O}$. Thus $(2) \Rightarrow (3)$.

Suppose that $\Lambda^+ \in \mathcal{L}(\mathcal{O})$ is $\mathcal{O}$-invariant, that $\beta$ is a $\Lambda^+$-generic line and that $V$ is an attracting neighborhood for $\beta$ with respect to the action of $\mathcal{O}^s$. Each $\mathcal{O}_\#^i(\beta) \in \mathcal{O}_\#^i(V)$ is generic with respect to $\Lambda^+$. Corollary 3.1.11 implies that $\beta \in \mathcal{O}_\#^i(V)$. Thus $U = V \cap \mathcal{O}_\#(V) \cap \cdots \cap \mathcal{O}_\#^{s-1}(V)$ is a neighborhood of $\beta$ that satisfies $\mathcal{O}_\#(U) \subset U$. Moreover, $\mathcal{O}_\#^s(U) \subset V$. It follows that $U$ is an attracting neighborhood for $\beta$. Thus $(1) \Rightarrow (2)$. $\qquad\square$

We conclude this subsection with a pair of lemmas that are needed for future reference.

Lemma 3.1.15.   *Assume that $H_r$ is an aperiodic exponentially-growing stratum for a train track map $f : G \to G$, that $\Lambda^+ \in \mathcal{L}(\mathcal{O})$ is associated to $H_r$ and that $\delta$ is a line in $\Lambda^+$ that is not entirely contained in $G_{r-1}$. Then the closure of $\delta$ is all of $\Lambda^+$. If $\delta$ is birecurrent then it is $\Lambda^+$-generic.*

*Proof.* Fix $k \geq 0$. By Lemma 3.1.10(3), each $\Lambda^+$-generic line has a $k$-tiling. Since $\delta$ is a weak limit of $\Lambda^+$-generic lines, $\delta$ is an increasing union of finite subpaths that have $k$-tilings. It follows (cf. the proof of Lemma 3.1.10(3)) that $\delta$ has a $k$-tiling. If $\delta \not\subset G_{r-1}$, then $\delta$ must contain at least one $k$-tile. Since $k$ is arbitrary, Corollary 3.1.11 and Lemma 3.1.8(3) imply that the closure of $\delta$ contains each $\Lambda^+$-generic line and so contains $\Lambda^+$. It follows immediately that $\delta$ satisfies condition (3) of Definition 3.1.5. Condition (2) of Definition 3.1.5 follows from the fact that every neighborhood of a generic line is also a neighborhood of $\delta$. If $\delta$ is birecurrent then condition (1) is also satisfied and $\delta$ is $\Lambda^+$-generic. $\qquad\square$

Lemma 3.1.16.   *A generic line of $\Lambda^+ \in \mathcal{L}(\mathcal{O})$ is never a circuit.*

*Proof.* A set in $\mathcal{B}$ consisting of a single circuit is closed in $\mathcal{B}$. If a $\Lambda^+$-generic line $\beta$ is a circuit, then $\Lambda^+ = \beta$. Choose a relative train track map $f : G \to G$ representing $\mathcal{O}$. Since $\mathcal{L}(\mathcal{O})$ is finite and invariant, the realization $\lambda \subset G$ for $\beta$ is invariant for the action of some iterate of $f_\#$. But this contradicts the fact that $\lambda \subset G_r$ is $r$-legal and crosses edges in $H_r$ for some exponentially-growing stratum $H_r$. $\qquad\square$

3.2. *Paired laminations.* In this subsection we define a pairing between $\mathcal{L}(\mathcal{O})$ and $\mathcal{L}(\mathcal{O}^{-1})$ that is analogous to the pairing between the stable and unstable foliations of a pseudo-Anosov homeomorphism. We will need the following lemma and corollary.

LEMMA 3.2.1. *Suppose that $F_n$ has generators $\{a_1, \ldots, a_n\}$ and that $F_{n-1}$ is the subgroup generated by $\{a_1, \ldots, a_{n-1}\}$. If $\Phi : F_n \to F_n$ is an automorphism and $F_{n-1}$ is $\Phi$-invariant, then $\Phi(a_n)$ contains exactly one occurrence of $a_n$ or $\bar{a}_n$.*

*Proof.* The restriction $\Phi|F_{n-1}$ extends by $a_n \mapsto a_n$ to an automorphism $\Phi' : F_n \to F_n$. After replacing $\Phi$ by $\Phi(\Phi')^{-1}$, we may assume that $\Phi|F_{n-1}$ is the identity.

Let $G = R_n$ be the rose with $n$ petals and let $f : G \to G$ be the obvious topological representative of $\Phi$. The edges of $G$ are labeled $e_1, \ldots e_n$ and the restriction of $f$ to the subgraph $G_{n-1}$ consisting of $e_1, \ldots e_{n-1}$ is the identity. If the $f$-image of an initial (respectively terminal) segment of $e_n$ is contained in $G_{n-1}$, then fold that initial (respectively terminal) segment into $G_{n-1}$. We may now assume that $f(e_n)$ begins and ends with $e_n$ or $\bar{e}_n$. If $f$ is an immersion, then it is a homeomorphism and we are done. Assume then that $f$ is not an immersion. The only fold that can take place is between the initial and terminal ends of $e_n$. Let $p : G \to G'$ be the maximal such fold and let $f' : G' \to G$ be the induced map (i.e. $f = f'p$). By construction, $f'$ is an immersion and hence a homeomorphism. But $G'$ has two vertices while $G$ has only one. This contradiction completes the proof. $\square$

COROLLARY 3.2.2. *If $f : G \to G$ is a topological representative and $H_i$ is a stratum that consists of a single edge $E$, then $f(E)$ crosses $E$, in either direction, at most once.*

*Proof.* We may assume without loss of generality that $G_i = G$. Suppose at first that $G_{i-1}$ is connected. Choose a maximal tree $T$ for $G$ that is contained in $G_{i-1}$ and that contains both endpoints of $E$. Let $B = \{a_1, \ldots, a_{n-1}, E\}$ be the basis of $F_n$ determined by the edges of $G \setminus T$, let $F_{n-1}$ be the subgroup generated by $\{a_1, \ldots, a_{n-1}\}$ and let $\Phi : F_n \to F_n$ be the automorphism induced by $f$ and $T$. The corollary follows directly from Lemma 3.2.1.

Suppose now that $G_{i-1}$ is not connected. Denote the components of $G_{i-1}$ by $C_1$ and $C_2$ and suppose at first that $C_1$ and $C_2$ are $f$-invariant. If $C_j$ has a single vertex, let $h_j : G \to G$ be the identity. If $C_j$ has at least two vertices, choose one, $x_j$, that is not an endpoint of $E$ and let $h_j : G \to G$ be a map with support in $C_j$ that is homotopic to the identity and that satisfies $h_j f(x_j) = x_j$. Let $f' : G \to G$ be the topological representative defined by $f'(e) = (h_1 h_2 f)_\#(e)$ for each edge $e$ of $G$. Since $f(E)$ and $f'(E)$ cross $E$ the same number of times, we may replace $f$ with $f'$. In particular, we may assume that $f$ fixes $x_1$ and $x_2$. Add an edge $F$ to $G$ with endpoints $x_1$ and $x_2$, and extend $f$ by $F \mapsto F$. This defines a topological representative (of a different outer automorphism) to which the previous argument applies.

If $f$ permutes $C_1$ and $C_2$ then, arguing as above, we may assume that there are vertices $x_j \in C_j$ that are permuted by $f$. Add $F$ to $G$ as above and extend $f$ by $F \mapsto \bar{F}$. The proof concludes as in the previous case. $\qquad\square$

*Definitions* 3.2.3. Since an attracting lamination $\Lambda^+$ is the closure of a single line $\beta$, any free factor that carries $\beta$ carries every line in $\Lambda^+$. Corollary 2.6.5 therefore implies that there is a unique free factor $F^i$ of minimal rank whose conjugacy class $[[F^i]]$ carries every line in $\Lambda^+$. We denote $[[F^i]]$ by $F(\Lambda^+)$. *The rank of $\Lambda^+$ is defined to be the rank of $F^i$.* For any connected graph $G$, the *rank of $G$ is defined to be the rank of $H_1(G)$.*

We say that the laminations $\Lambda^+$ and $\Lambda^-$ of the following lemma are *paired*. A key point in the proof is that $\mathcal{O}$ and $\mathcal{O}^{-1}$ have the same invariant free factor systems.

LEMMA 3.2.4. *For each $\Lambda^+ \in \mathcal{L}(\mathcal{O})$ there is a unique $\Lambda^- \in \mathcal{L}(\mathcal{O}^{-1})$ such that $F(\Lambda^+) = F(\Lambda^-)$.*

*Proof.* We induct on the rank $n$ of $F_n$. If $n = 1$, then there are no exponentially-growing strata in any relative train track map representing an iterate of $\mathcal{O}$ and so $\mathcal{L}(\mathcal{O})$ is empty. We may therefore assume that the lemma holds for outer automorphisms of $F_k$ where $k < n$.

There is no loss in replacing $\mathcal{O}$ by an iterate. We may therefore assume that each element of $\mathcal{L}(\mathcal{O})$ is $\mathcal{O}$-invariant. Choose $\Lambda^+ \in \mathcal{L}(\mathcal{O})$. Since $F(\Lambda^+)$ is unique, $F(\Lambda^+)$ is $\mathcal{O}$-invariant. If the rank of $\Lambda^+$ is less than $n$, then the inductive hypothesis, applied to the restriction of $\mathcal{O}$ to $F(\Lambda^+)$, implies that there exists a unique $\Lambda^- \in \mathcal{L}(\mathcal{O}^{-1})$ such that $F(\Lambda^+) = F(\Lambda^-)$. We may therefore assume that there is a pairing between the elements of $\mathcal{L}(\mathcal{O})$ with rank less than $n$ and the the elements of $\mathcal{L}(\mathcal{O}^{-1})$ with rank less than $n$.

For any exponentially-growing stratum $H_r$ there are bi-infinite paths in $G_r$ crossing edges in $H_r$. It follows that the rank of each component of $G_{r-1}$ is less than the rank of $G_r$. In particular, there is at most one element of $\mathcal{L}(\mathcal{O})$ or $\mathcal{L}(\mathcal{O}^{-1})$ with rank $n$. It therefore suffices to assume that there is an element $\Lambda^+ \in \mathcal{L}(\mathcal{O})$ with rank $n$ and to prove that there is an element $\Lambda^- \in \mathcal{L}(\mathcal{O}^{-1})$ with rank $n$.

After we replace $\mathcal{O}$ by a further iterate if necessary, there is (Lemma 2.6.7) a reduced relative train track map $f : G \to G$ and filtration $\emptyset = G_0 \subset G_1 \subset \cdots \subset G_K = G$ representing $\mathcal{O}$. Now, $\Lambda^+$ is associated to the top stratum $H_K$ and each element of $\mathcal{L}(\mathcal{O})$ and $\mathcal{L}(\mathcal{O}^{-1})$ with rank less than $n$ is carried by $\mathcal{F} = \mathcal{F}(G_{K-1})$. Corollary 3.2.2 implies that $H_K$ is not a single edge; since subdivision vertices in $H_K$ can be removed without loss, $H_K$ is not an arc. It follows that $\mathcal{F}$ is neither a single (conjugacy class of a) free factor of rank $n-1$ nor a pair of (conjugacy classes of) free factors with rank sum equal to $n$.

Choose (Lemma 2.6.7) a relative train track map $f' : G' \to G'$ representing an iterate of $\mathcal{O}^{-1}$ such that $\mathcal{F}$ is realized by a filtration element. The transition submatrix for a non-exponentially-growing stratum is a permutation and so has an iterate that equals the identity. We may therefore assume, after replacing $f'$ by an iterate and enlarging the filtration if necessary, that each non-exponentially-growing stratum is a single edge. We may also assume that $f : G \to G$ is eg-aperiodic. If $H'_{K'}$ is the topmost stratum then, since $f : G \to G$ is reduced, $\mathcal{F} = \mathcal{F}(G'_{K'-1})$. Since $H'_{K'}$ cannot be a zero stratum, the concluding observation of the preceding paragraph rules out the possibility that it is a single edge. Thus $H'_{K'}$ is exponentially growing. Since the expanding lamination associated to $H'_{K'}$ is not carried by $\mathcal{F}$, it must have rank $n$.                                                                                      □

### 3.3. Expansion factors.

In this subsection we assume that $\Lambda^+$ is an attracting lamination for some element of $\mathrm{Out}(F_n)$. Define the *stabilizer of* $\Lambda^+$ to be $\mathrm{Stab}(\Lambda^+) = \{\Psi \in \mathrm{Out}(F_n) : \Psi_\#(\Lambda^+) = \Lambda^+\}$. The following corollary (of Proposition 3.3.3 below) is the main result of this subsection; it is essential to our reduction of the Tits alternative for $\mathrm{Out}(F_n)$ to the Tits alternative for $UPG(F_n)$.

COROLLARY 3.3.1.   *There is a homomorphism* $PF_{\Lambda^+} : \mathrm{Stab}(\Lambda^+) \to \mathbb{Z}$ *such that* $\Psi \in \mathrm{Ker}(PF_{\Lambda^+})$ *if and only if* $\Lambda^+ \notin \mathcal{L}(\Psi)$ *and* $\Lambda^+ \notin \mathcal{L}(\Psi^{-1})$.

The analogous result for the mapping class group is an immediate corollary of the fact (exposé 12 of [FLP]) that the measured foliations associated to a pseudo-Anosov homeomorphism are uniquely ergodic. Any mapping class that topologically preserves the measured foliation must projectively fix its invariant transverse measure and so multiplies this transverse measure by some scalar factor. The assignment of the logarithm of this scalar factor to the mapping class defines a homomorphism to $\mathbb{R}$. After rescaling this provides a homomorphism to $\mathbb{Z}$ analogous to $PF_{\Lambda^+}$.

Because we are working in $\mathcal{B}$ and not a more structured space that takes measures into account, we cannot measure the attraction factor directly. Instead of an invariant measure defined on the lamination itself, we use a length function on paths in a marked graph. The length function depends on the choice of marked graph but the factor by which an element of $\mathrm{Stab}(\Lambda^+)$ expands this length does not.

*Definition* 3.3.2. Assume that $f : G \to G$ and $\emptyset = G_0 \subset G_1 \subset \cdots \subset G_K = G$ are a relative train track map and filtration for an element of $\mathrm{Stab}(\Lambda^+)$ and that $\Lambda^+$ is the attracting lamination associated to the (necessarily aperiodic) exponentially-growing stratum $H_r$. For any path $\sigma \subset G$ define $\mathrm{EL}_r(\sigma)$ to be the edge length of $\sigma$, counting only the edges of $H_r$ that are entirely contained in $\sigma$. We say that $\Psi \in \mathrm{Stab}(\Lambda^+)$ *asymptotically expands* $\Lambda^+$ *by the factor* $\mu$ if for every such choice of $\emptyset = G_0 \subset G_1 \subset \cdots \subset G_K = G$ and $f : G \to G$, every topological representative $g : G \to G$ of $\Psi$ and for all $\eta > 0$

$$(*) \qquad\qquad \mu - \eta < \frac{EL_r(g_\#(\sigma))}{EL_r(\sigma)} < \mu + \eta$$

whenever $\sigma$ is contained in a $\Lambda^+$-generic line and $EL_r(\sigma)$ is sufficiently large.

*For the remainder of this subsection we assume that* $f : G \to G$, $H_r$, $\emptyset = G_0 \subset G_1 \subset \cdots \subset G_K = G$ *and* $\Lambda^+$ *are as in Definition* 3.3.2.

The following proposition relates asymptotic expansion of $\Lambda^+$ to Perron-Frobenius eigenvalues. In particular, it implies that the Perron-Frobenius eigenvalue associated to an exponentially-growing stratum of a relative train track map $f : G \to G$ depends only on the outer automorphism $\mathcal{O}$ determined by $f$ and on the element of $\mathcal{L}(\mathcal{O})$ that is associated to the stratum.

PROPOSITION 3.3.3. (1) *Every* $\Psi \in \mathrm{Stab}(\Lambda^+)$ *asymptotically expands* $\Lambda^+$ *by some factor* $\mu = \mu(\Psi)$.

(2) $\mu(\Psi\Psi') = \mu(\Psi)\mu(\Psi')$.

(3) $\mu(\Psi) > 1$ *if and only if* $\Lambda^+ \in \mathcal{L}(\Psi)$.

(4) *If* $\Lambda^+ \in \mathcal{L}(\Psi)$ *and* $f' : G' \to G'$ *is a relative train track map for* $\Psi$, *then* $\mu(\Psi) = \mu'_s$ *is the Perron-Frobenius eigenvalue for the transition submatrix* $M'_s$ *of the exponentially-growing stratum* $H'_s$ *associated to* $\Lambda^+$.

Our main result follows easily from Proposition 3.3.3.

*Proof of Corollary* 3.3.1. Define $PF^*_\Lambda(\Psi) = \log(\mu(\Psi))$. Proposition 3.3.3 and the observation that $PF^*_\Lambda(\Psi^{-1}) = -PF^*_\Lambda(\Psi)$ (which follows from Proposition 3.3.3(2)) imply that each $\mu(\Psi)$, other than 1, occurs as the Perron-Frobenius eigenvalue for an irreducible matrix of uniformly bounded size. It follows (cf. page 37 of [BH1]) that the image of $PF^*_\Lambda$ is an infinite discrete subset of $\mathbb{R}$ and is hence isomorphic to $\mathbb{Z}$. Identify the image with $\mathbb{Z}$ and call the resulting homomorphism $PF_\Lambda$. The desired properties follow immediately from Proposition 3.3.3. □

We need a pair of preliminary estimates before beginning the proof of Proposition 3.3.3.

*Definitions* 3.3.4. Let $\{E_i\}$ be the edges of $H_r$ and let $\mu_r$ be the Perron-Frobenius eigenvalue for $M_r$. The Perron-Frobenius theorem ([Sen]) implies that $\mu_r^{-n} M_r^n$ converges to a matrix $M^*$ whose columns are projectively equal. Normalize the column vectors of $M^*$ so that the sum of the entries is one and denote this common *frequency vector* by $A = (a_i)$. Let $\tau_i^k$ be the $k$-tile $f_\#^k(E_i)$.

Every $\Lambda^+$-generic line $\lambda$ has a unique 0-tiling. Pushing this forward by $f_\#$, as in Lemma 3.1.8(2), produces a 1-tiling of $f_\#(\lambda)$ that we call the standard 1-tiling of $f_\#(\lambda)$. Continue this to define the *standard $k$-tiling* of $f_\#^k(\lambda)$. Since $\Lambda^+$ is $\mathcal{O}$-invariant, $\lambda = f_\#^k(\gamma_k)$ for some $\Lambda^+$-generic line $\gamma_k$. In this way every $\Lambda^+$-generic line has a standard $k$-tiling for all $k \geq 0$.

The following lemma states that the $k$-tiles are 'evenly distributed' in $\Lambda$-generic lines.

LEMMA 3.3.5. *Fix $\varepsilon > 0$ and $k \geq 0$. Suppose that $\sigma$ is a finite subpath of a $\Lambda^+$-generic line $\lambda$. Among all $k$-tiles in the standard $k$-tiling of $\lambda$ that are entirely contained in $\sigma$, denote the fraction that equals $\tau_i^k$ by $\alpha_{ik}(\sigma)$. If $EL_r(\sigma)$ is sufficiently large, then $a_i - \varepsilon < \alpha_{ik}(\sigma) < a_i + \varepsilon$.*

*Proof.* Choose $l > 0$. In the standard $l$-tiling of $\lambda$, there are at most two $l$-tiles that intersect $\sigma$ but are not entirely contained in $\sigma$. If $EL_r(\sigma)$ is sufficiently large, then there is no loss in ignoring the $H_r$ edges of $\sigma$ in these two $l$-tiles and we may assume that $\sigma$ is a union of $l$-tiles and paths in $G_{r-1}$. It therefore suffices to assume that $\sigma$ is an $l$-tile for some arbitrarily large $l$. By Lemma 2.5.1 and the definitions, it suffices to assume that $k = 0$. This case is an immediate consequence of Lemma 3.1.8(4). $\qquad\square$

LEMMA 3.3.6. *Assume that $g : G \to G$ is a topological representative and that $\Lambda^+$ is $g_\#$-invariant. There is a constant $C_1 = C_1(g)$ satisfying $EL_r(g_\#(\delta)) < C_1$ for any subpath $\delta \subset G_{r-1}$ of a $\Lambda^+$-generic line.*

*Proof.* If the lemma fails, then there exist $\Lambda$-generic lines $\lambda_j$ and finite subpaths $\delta_j \subset G_{r-1}$ of $\lambda_j$ such that the central segment of $g_\#(\delta_j)$ obtained by removing the first and last $j$ edges contains at least one edge in $H_r$. After passing to a subsequence, we may assume that the $\delta_j$'s are an increasing sequence whose union is a line $\delta^* \subset G_{r-1}$ with the property that $g_\#(\delta^*) \not\subset G_{r-1}$. Since $\delta^*$ is a line in $\Lambda^+$ whose closure is not all of $\Lambda^+$ and since $g_\#$ is a homeomorphism that preserves $\Lambda^+$, $g_\#(\delta^*)$ is a line in $\Lambda^+$ whose closure is not all of $\Lambda^+$. This contradicts Lemma 3.1.15 and so completes the proof. $\qquad\square$

*Proof of Proposition* 3.3.3. Assume that $f : G \to G$, $\emptyset = G_0 \subset G_1 \subset \cdots \subset G_K = G$, $H_r$ and $g : G \to G$ are as in Definition 3.3.2, that $a_i$ and $\tau_i^k$ are as in Definitions 3.3.4 and that $\sigma$ is a finite subpath of a $\Lambda^+$-generic line. Define

$$\mu_k = \frac{\sum_i a_i EL_r(g_\#(\tau_i^k))}{\sum_i a_i EL_r(\tau_i^k)}.$$

We will show that for all $\varepsilon > 0$ : if $k$ is sufficiently large (relative to $\varepsilon$) and if $EL_r(\sigma)$ is sufficiently large (relative to $\varepsilon$ and $k$) then

$$(**) \qquad (1-\varepsilon)\mu_k \leq \frac{EL_r(g_\#(\sigma))}{EL_r(\sigma)} \leq (1+\varepsilon)\mu_k.$$

It will follow that $\mu = \lim \mu_k$ exists and that for all $\eta > 0$, $(*)$ holds whenever $EL_r(\sigma)$ is sufficiently large (relative to $\eta$).

We will verify $(**)$ relative to the choices of $f : G \to G$, $\emptyset = G_0 \subset G_1 \subset \cdots \subset G_K = G$, $H_r$, and $g : G \to G$ and then check that $\mu$ does not depend on these choices.

Let $\mathrm{bcc}(g)$ be a bounded cancellation constant for $g : G \to G$ and let $C_1$ be the constant of Lemma 3.3.6. In the following list $\sim_\varepsilon$ means that the error of approximation is small relative to $\varepsilon$. Choose $k$ so large that for all $i$:

(A) $C_1/EL_r(\tau_i^k) \sim_\varepsilon 0$,

(B) $\mathrm{bcc}(g)/EL_r(\tau_i^k) \sim_\varepsilon 0$.

Now restrict attention to $\sigma$ where $EL_r(\sigma)$ is so large that:

(C) $\alpha_{ik}(\sigma) \sim_\varepsilon a_i$,

(D) $EL_r(\tau_i^k)/EL_r(\sigma) \sim_\varepsilon 0$.

We can now verify $(**)$. In order to approximate $EL_r(g_\#(\sigma))$ we allow ourselves to make assumptions that result in errors that are a small percentage of the total. There are four such assumptions. The first is that $\sigma$ begins and ends with a $k$-tile in the standard $k$-tiling of $\lambda$. Thus $\sigma$ is a concatenation of $k$-tiles $\gamma_j$ and maximal subpaths $\delta_l \subset G_{r-1}$; let $N(\sigma)$ be the number of $k$-tiles in this decomposition. The error that this assumption contributes is at most twice the number of $H_r$-edges in a $k$-tile and so is controlled by (D). The second assumption is that the approximation in (C) is exact. The third is that each $EL_r(g_\#(\delta_l)) = 0$; this error is controlled by (A). The final assumption is that $g_\#(\sigma) \subset g_\#(\lambda)$ is a concatenation of the $g_\#(\gamma_j)$'s and $g_\#(\delta_l)$'s (with no cancellation at the junctures). This produces an error that is bounded by $2\mathrm{bcc}(g)N(\sigma)$ and so is controlled by (B).

Once these assumptions are made, $EL_r(g_\#(\sigma)) = \sum_i N(\sigma)a_i EL_r(g_\#(\tau_i^k))$. When we apply this again with $g = $ identity, $EL_r(\sigma) = \sum_i N(\sigma)a_i EL_r(\tau_i^k)$. Thus $\frac{EL_r(g_\#(\sigma))}{EL_r(\sigma)} = \mu_k$ and we have verified $(**)$.

If $g^* : G \to G$ is another topological representative of $\Psi$, then there are lifts $\tilde{g} : \Gamma \to \Gamma$ and $\tilde{g}^* : \Gamma \to \Gamma$ such that the distance between $\tilde{g}(\tilde{x})$ and $\tilde{g}^*(\tilde{x})$ is bounded independently of $\tilde{x}$. It follows that $EL_r(g_\#(\sigma)) - EL_r(g_\#^*(\sigma))$ is bounded independently of $\sigma$ and hence that $\mu$ does not depend on the choice of $g$.

Suppose next that $\hat{f} : \hat{G} \to \hat{G}$ and $\emptyset = \hat{G}_0 \subset \hat{G}_1 \subset \cdots \subset \hat{G}_{\hat{K}} = \hat{G}$ is another relative train track map and filtration representing $\mathcal{O}$, that $\hat{\Lambda}^+$ is the attracting lamination associated to the stratum $\hat{H}_s$ and that $\widehat{EL}_s$ is the edge length function that counts edges of $\hat{H}_s$ in $\hat{G}$.

Choose a homotopy equivalence $h : G \to \hat{G}$ that respects the markings and that restricts to an immersion on each edge. Arguing exactly as above we conclude that there is a positive constant $\nu$ so that for all $\varepsilon > 0$

$$(***) \qquad (1-\varepsilon)\nu < \frac{\widehat{EL}_s(h_\#(\sigma))}{EL_r(\sigma)} < (1-\varepsilon)\nu$$

whenever $\sigma$ is contained in a $\Lambda^+$-generic line and $EL_r(\sigma)$ is sufficiently large. The details of this modification are left to the reader.

Suppose now that $\hat{g} : \hat{G} \to \hat{G}$ represents $\Psi$. For any finite path $\sigma \subset G_r$, $h_\# g_\#(\sigma)$ and $\hat{g}_\# h_\#(\sigma)$ differ by initial and terminal segments of uniformly bounded size. It follows, by (***), that

$$\frac{\widehat{EL}_s(\hat{g}_\# h_\#(\sigma))}{\widehat{EL}_s(h_\#(\sigma))} \sim \frac{\widehat{EL}_s(h_\# g_\#(\sigma))}{\widehat{EL}_s(h_\#(\sigma))} \sim \frac{EL_r(g_\#(\sigma))}{EL_r(\sigma)}$$

where the error of approximation goes to 0 as $EL_r(\sigma) \to \infty$ or equivalently as $\widehat{EL}_s(h_\#(\sigma)) \to \infty$. We conclude that $\mu$ is independent of the choice of $f : G \to G$ and $\emptyset = G_0 \subset G_1 \subset \cdots \subset G_K = G$. This completes the proof of part (1) of Proposition 3.3.3.

Suppose that $g : G \to G$ and $g' : G \to G$ are topological representatives for $\Psi$ and $\Psi'$ respectively. If $\sigma \subset G$ is contained in a $\Lambda^+$-generic line, then there is a subpath $\sigma' \subset G$ of $g_\#(\sigma)$ that is contained in a $\Lambda^+$-generic line and that differs from $g_\#(\sigma)$ only in an initial and terminal segment of uniformly bounded length. Similarly, $g'_\# g_\#(\sigma) \subset G$ differs from $g'_\#(\sigma')$ only in an initial and terminal segment of uniformly bounded length. Thus $\mu(\Psi\Psi') = \mu(\Psi)\mu(\Psi')$ and we have proved part (2).

Suppose now that $\mu(\Psi) > 1$ and that $\sigma_0$ is a finite subpath of a $\Lambda^+$-generic line $\lambda$. Let $\sigma_1$ be the subpath of $g_\#(\lambda)$ obtained from $g_\#(\sigma_0)$ by removing the initial and terminal subpaths of length $\mathrm{bcc}(g)$. Lemma 2.3.1(3) implies that $g_\#(N(\sigma_0)) \subset N(\sigma_1)$. If $EL_r(\sigma_0)$ is sufficiently large, then $EL_r(\sigma_1) > EL_r(\sigma_0)$ and we may iterate the argument to produce $\sigma_k$ with increasing $EL_r$-length such that $g_\#(N(\sigma_{k-1})) \subset N(\sigma_k)$ and hence $g_\#^k(N(\sigma_0)) \subset N(\sigma_k)$. Since $g_\#^k(\lambda)$ is $\Lambda^+$-generic, Corollary 3.1.11 implies that $\lambda \in N(\sigma_k)$ for all $k$.

By Lemma 3.1.10(4) $\sigma_0$ is contained in some tile, say an $l$-tile, in $\lambda$. By Lemma 3.3.5 there exists $k$ so that $\sigma_k$ contains every $l$-tile. In particular, $N(\sigma_k) \subset N(\sigma_0)$. Lemma 3.3.5 and Corollary 3.1.11 imply that the $N(\sigma_j)$'s are a neighborhood basis for $\lambda$ and hence that $N(\sigma_0)$ is an attracting neighborhood of $\lambda$ for the action of $\Psi_\#^k$. It follows that $\Lambda^+ \in \mathcal{L}(\Psi)$ and we have verified half of (3).

If $\Lambda^+ \in \mathcal{L}(\Psi)$, we may assume that the relative train track map $f : G \to G$ used to compute $\mu$ represents $\Psi$, that $g = f$ and that $\sigma$ is a $k$ tile $\tau_i^k$ for some large $k$. With these assumptions $EL_r(g_\#(\sigma))$ is the $i^{\mathrm{th}}$ column sum of $M_r^{k+1}$ and $EL_r(\sigma)$ is the $i^{\mathrm{th}}$ column sum of $M_r^k$. The Perron-Frobenius theorem (cf. the definition of $M^*$ above) implies that $EL_r(g_\#(\sigma))/EL_r(\sigma) \to \mu_r$ as $k \to \infty$. This completes the proof of parts (3) and (4). $\qquad\square$

3.4. *Detecting $F_2$ via laminations.* We now show how expanding laminations can be used to prove that a group of outer automorphisms contains a free subgroup of rank two. Our criterion is based on a technique of Tits (Proposition 1.1 of [Tit]), a version of which appears in the next lemma.

LEMMA 3.4.1. *Suppose that a group $\mathcal{H}$ acts on a space $X$, that there are subsets $U^+, U^-, V^+$ and $V^-$ of $X$, a point $x \in X$ and elements $f, g \in \mathcal{H}$ such that*

1. $x \notin (U^+ \cup U^- \cup V^+ \cup V^-)$,

2. $f(\{x\} \cup U^+ \cup V^+ \cup V^-) \subset U^+$,

3. $f^{-1}(\{x\} \cup U^- \cup V^+ \cup V^-) \subset U^-$,

4. $g(\{x\} \cup U^+ \cup U^- \cup V^+) \subset V^+$,

5. $g^{-1}(\{x\} \cup U^+ \cup U^- \cup V^-) \subset V^-$.

*Then the subgroup of $\mathcal{H}$ generated by $f$ and $g$ is isomorphic to $F_2$.*

In our case, $\mathcal{H}$ will be a subgroup of Out($F_n$) and $X$ will be the space $\mathcal{B}$. The 'ping pong' method of Tits can be recast as follows.

LEMMA 3.4.2. *Suppose that:*

- $\Lambda^+ \in \mathcal{L}(\mathcal{O})$ *and* $\Lambda^- \in \mathcal{L}(\mathcal{O}^{-1})$ *are paired and $\mathcal{O}$-invariant.*

- $\Gamma^+ \in \mathcal{L}(\Psi)$ *and* $\Gamma^- \in \mathcal{L}(\Psi^{-1})$ *are paired and $\Psi$-invariant.*

- *Generic lines in $\Lambda^+$ and $\Lambda^-$ are weakly attracted to $\Gamma^+$ (respectively $\Gamma^-$) under the action of $\Psi$ (respectively $\Psi^{-1}$).*

- *Generic lines in $\Gamma^+$ and $\Gamma^-$ are weakly attracted to $\Lambda^+$ (respectively $\Lambda^-$) under the action of $\mathcal{O}$ (respectively $\mathcal{O}^{-1}$).*

*Then $\mathcal{O}^N$ and $\Psi^N$ generate a free subgroup of rank two for all sufficiently large $N$.*

*Proof.* Let $\lambda^\pm$ and $\gamma^\pm$ be generic lines for $\Lambda^\pm$ and $\Gamma^\pm$; let $U^\pm$ and $V^\pm$ be attracting neighborhoods of $\lambda^\pm$ and $\gamma^\pm$ respectively.

There exists $k \geq 0$ such that $\mathcal{O}_\#^k(\gamma^+) \subset U^+$. Since $\{\Psi_\#^l(V^+)\}$ is a neighborhood basis for $\gamma^+$, there exists $l \geq 0$ such that $\mathcal{O}_\#^k(\Psi_\#^l(V^+)) \subset U^+$. This inclusion remains valid if $k$ and/or $l$ are increased.

Repeating this argument on various combinations of $\mathcal{O}^{\pm 1}$ and $\Psi^{\pm 1}$, we see that for all sufficiently large $K$ and $L$,

- $\mathcal{O}_\#^K \Psi_\#^L(V^+), \mathcal{O}_\#^K \Psi_\#^{-L}(V^-) \subset U^+$,

- $\mathcal{O}_\#^{-K} \Psi_\#^L(V^+), \mathcal{O}_\#^{-K} \Psi_\#^{-L}(V^-) \subset U^-$,

- $\Psi_\#^K \mathcal{O}_\#^L(U^+), \Psi_\#^K \mathcal{O}_\#^{-L}(U^-) \subset V^+$,

- $\Psi_\#^{-K} \mathcal{O}_\#^L(U^+), \Psi_\#^{-K} \mathcal{O}_\#^{-L}(U^-) \subset V^-$.

Defining $U_L^+$ (respectively $U_L^-, V_L^+, V_L^-$) by $\mathcal{O}_\#^L(U^+)$ (respectively $\mathcal{O}_\#^{-L}(U^-)$, $\Psi_\#^L(V^+)$, $\Psi_\#^{-L}(V^-)$) and defining $N = K + L$, we have

(1) $\mathcal{O}_\#^N(V_L^+), \mathcal{O}_\#^N(V_L^-) \subset U_L^+$,

(2) $\mathcal{O}_\#^{-N}(V_L^+), \mathcal{O}_\#^{-N}(V_L^-) \subset U_L^-$,

(3) $\Psi_\#^N(U_L^+), \Psi_\#^N(U_L^-) \subset V_L^+$,

(4) $\Psi_\#^{-N}(U_L^+), \Psi_\#^{-N}(U_L^-) \subset V_L^-$.

Since $U^\pm$ and $V^\pm$ are attracting neighborhoods, we also have

(5) $\mathcal{O}_\#^N(U_L^+) \subset U_L^+; \mathcal{O}_\#^{-N}(U_L^-) \subset U_L^-; \Psi_\#^N(V_L^+) \subset V_L^+;$ and $\Psi_\#^{-N}(V_L^-) \subset V_L^-$.

Choose a circuit $\sigma \subset G$ that is weakly attracted to $\Lambda^+$ under the action of $\mathcal{O}$ and to $\Lambda^-$ under the action of $\mathcal{O}^{-1}$. By (3) and (4), $\mathcal{O}_\#^m(\sigma)$ is weakly attracted to $\Gamma^+$ (respectively $\Gamma^-$) under the action of $\Psi$ (respectively $\Psi^{-1}$) for all sufficiently large $m$, say $m \geq M$. Let $x = \mathcal{O}_\#^M(\sigma) \in \mathcal{B}$. Since (Lemma 3.1.16) generic lines in attracting laminations cannot be circuits, we can choose $L$ so that $x \notin (U_L^+ \cup U_L^- \cup V_L^+ \cup V_L^-)$. For large $K$ and hence large $N$, $\mathcal{O}_\#^N(x) \in U_L^+$; $\mathcal{O}_\#^{-N}(x) \in U_L^-$; $\Psi_\#^N(x) \in V_L^+$; and $\Psi_\#^{-N}(x) \in V_L^-$.

We have now verified all the hypotheses of Lemma 3.4.1. $\qquad \square$

Now we state the specific application of Lemma 3.4.2 to be used.

COROLLARY 3.4.3. *Suppose that $\Lambda^+ \in \mathcal{L}(\mathcal{O})$ and $\Lambda^- \in \mathcal{L}(\mathcal{O}^{-1})$ are paired and $\mathcal{O}$-invariant, that $\mathcal{H}$ is a subgroup of $\mathrm{Out}(F_n)$ containing $\mathcal{O}$ and that there is an element $\psi \in \mathcal{H}$ such that generic lines of the four laminations $\psi^{\pm 1}(\Lambda^\pm)$ are weakly attracted to $\Lambda^+$ under the action of $\mathcal{O}$ and to $\Lambda^-$ under the action of $\mathcal{O}^{-1}$. Then $\mathcal{H}$ contains a free subgroup of rank two.*

*Proof.* Define $\Psi = \psi \mathcal{O} \psi^{-1}$ and note that $\Gamma^+ = \psi_{\#}(\Lambda^+) \in \mathcal{L}(\Psi^+)$ and $\Gamma^- = \psi_{\#}(\Lambda^-) \in \mathcal{L}(\Psi^{-1})$ are paired and $\Psi$-invariant. The condition on $\psi_{\#}^{\pm 1}(\Lambda^{\pm})$ can be restated as:

1. $\Gamma^{\pm}$ is weakly attracted to $\Lambda^+$ (respectively $\Lambda^-$) under the action of $\mathcal{O}$ (respectively $\mathcal{O}^{-1}$).

2. $\psi_{\#}^{-1} \Lambda^{\pm}$ is weakly attracted to $\psi_{\#}^{-1}(\Gamma^+)$ (respectively $\psi_{\#}^{-1}(\Gamma^-)$) under the action of $\psi^{-1} \Psi \psi$ (respectively $\psi^{-1} \Psi^{-1} \psi$).

This last condition is equivalent to:

3. $\Lambda^{\pm}$ is weakly attracted to $\Gamma^+$ (respectively $\Gamma^-$) under the action of $\Psi$ (respectively $\Psi^{-1}$).

The corollary therefore follows from Lemma 3.4.2. $\qquad\square$

## 4. Splittings

4.1. *Preliminaries and non-exponentially-growing strata.* To understand the action of an outer automorphism $\mathcal{O}$ on the space $\mathcal{B}$ of abstract lines, we choose a relative train track map $f : G \to G$ representing $\mathcal{O}$ and study the induced action $f_{\#}$ on the space $\mathcal{B}(G)$ of lines in $G$. One advantage of working in $\mathcal{B}(G)$ is that it is possible to subdivide a bi-infinite path in $G$ into subpaths in $G$. We will be interested in subdivisions in which the action of $f_{\#}$ on the whole is the 'sum' of the action of $f_{\#}$ on the parts. We make this precise as follows.

Suppose that $\sigma = \ldots \sigma_{l-1} \sigma_l \ldots$ is a decomposition of a path or circuit $\sigma \subset G$ into nontrivial subpaths. If $\sigma$ is a finite path or a circuit, then the decomposition is assumed to be finite but infinite paths may have infinite decompositions. Decompositions of a circuit $\sigma = \sigma_1 \ldots \sigma_l$ are always assumed to be cyclic; in particular, $\sigma_1$ and $\bar{\sigma}_l$ do not have a common initial segment. If $\sigma$ is a path then we assume that there are at least two subpaths in the decomposition but if $\sigma$ is a circuit then we allow $\sigma = \sigma_1$.

We say that $\sigma = \ldots \sigma_{l-1} \sigma_l \ldots$ is a *k-splitting* if

$$f_{\#}^k(\sigma) = \ldots f_{\#}^k(\sigma_{l-1}) f_{\#}^k(\sigma_l) \ldots$$

is a decomposition into subpaths and is a *splitting* if it is a k- splitting for all $k > 0$.

If $\tilde{f} : \Gamma \to \Gamma$ and $\tilde{\sigma} \subset \Gamma$ are lifts, then a decomposition $\tilde{\sigma} = \ldots \tilde{\sigma}_{l-1} \tilde{\sigma}_l \ldots$ into subpaths (with at least two pieces and with finitely many pieces if $\tilde{\sigma}$ is finite) is called a *k-splitting* if $\tilde{f}_{\#}^k(\tilde{\sigma}) = \ldots \tilde{f}_{\#}^k(\tilde{\sigma}_{l-1}) \tilde{f}_{\#}^k(\tilde{\sigma}_l) \ldots$ is a decomposition into subpaths, and is called a *splitting* if it is a k-splitting for all $k \geq 0$. If $\sigma$ is a path in $G$ then every $k$-splitting (splitting) of $\tilde{\sigma}$ projects to a $k$-splitting (splitting) of $\sigma$. If $\sigma$ is a circuit and $\tilde{\sigma}$ is the axis of the covering translation $T$, then $T$-invariant $k$-splittings (splittings) of $\tilde{\sigma}$ project to $k$-splittings (splittings) of $\sigma$.

Decompositions of $\tilde{\sigma}$ are determined by the juncture points $\tilde{J}$ of the subpaths. If the decomposition determined by $\tilde{J}$ is a $k$-splitting (splitting) then we say that $\tilde{\sigma}$ *can be $k$-split (split) at $\tilde{J}$*. If $\tilde{J}$ contains a single point $\tilde{x}$, then we say that $\tilde{\sigma}$ *can be $k$-split (split) at $\tilde{x}$*. As a matter of notation, we will only use $\cdot$ to separate subpaths if the separation is a splitting.

We first record some elementary properties of $k$-splittings and splittings.

LEMMA 4.1.1. (1) *If $\sigma$ is a circuit, $\sigma = \sigma_1'$ is a splitting and $\sigma_1' = \sigma_1 \cdot \ldots \cdot \sigma_l$, then $\sigma = \sigma_1 \cdot \ldots \cdot \sigma_l$. In other words, splittings of the path $\sigma_1$ determine splittings of the circuit $\sigma$.*

*If $\sigma \subset G$ is a path or circuit then*:

(2) *If $\sigma = \sigma_1 \cdot \sigma_2$ and $\sigma_1 = \sigma_1' \cdot \sigma_2'$ then $\sigma = \sigma_1' \cdot \sigma_2' \cdot \sigma_2$. The analogous result with the roles of $\sigma_1$ and $\sigma_2$ reversed also holds.*

(3) *$\tilde{\sigma}$ can be $k$-split at $\tilde{x}$ if and only if $\tilde{f}^k(\tilde{x}) \in \tilde{f}_{\#}^k(\tilde{\sigma})$.*

(4) *$\{\tilde{x} \in \tilde{\sigma} : \tilde{\sigma}$ can be $k$-split at $\tilde{x}\}$ is closed.*

(5) *If $\tilde{f}_{\#}(\tilde{\sigma})$ can be split at $\tilde{y}$ and if $\tilde{x} \in \tilde{\sigma}$ satisfies $\tilde{f}(\tilde{x}) = \tilde{y}$, then $\tilde{\sigma}$ can be split at $\tilde{x}$.*

*If $\sigma \subset G$ is a path then*

(6) *Assume that $\alpha = \alpha_1 \alpha_2$ is a $k$-splitting, that $\sigma = \alpha\beta$ is a decomposition into subpaths and that not all of $f_{\#}^k(\alpha_2)$ is canceled when $f_{\#}^k(\alpha)f_{\#}^k(\beta)$ is tightened to $f_{\#}^k(\sigma)$. Let $\beta' = \alpha_2\beta$. Then $\sigma = \alpha_1\beta'$ is a $k$-splitting.*

*Proof.* Parts (1), (2), and (3) follow immediately from the definitions. Parts (4), (5), and (6) follow from (3). □

The next lemma complements Lemma 4.1.1(1).

LEMMA 4.1.2. *Every circuit $\sigma \subset G$ has a splitting $\sigma = \sigma_1$.*

*Proof.* Choose lifts $\tilde{f} : \Gamma \to \Gamma$ and $\tilde{\sigma} \subset \Gamma$ and let $T : \Gamma \to \Gamma$ be the covering translation with axis $\tilde{\sigma}$. The set $\tilde{S}_k = \{\tilde{x} \in \tilde{\sigma} : \tilde{f}^k(\tilde{x}) \in \tilde{f}_{\#}^k(\tilde{\sigma})\}$ is closed. An easy induction argument shows that $f^N$ maps $\cap_{k=1}^N \tilde{S}_k$ onto $\tilde{f}_{\#}^N(\tilde{\sigma})$ for all $N \geq 1$. Since $\cap_{k=1}^N \tilde{S}_k$ is $T$-invariant and nonempty, it must intersect each fundamental domain of $\tilde{\sigma}$. Thus the $T$-invariant set $\cap_{k=1}^{\infty} \tilde{S}_k$ is nonempty. Choose a $T$-orbit $\tilde{J} \subset \cap_{k=1}^{\infty} \tilde{S}_k$. If $\tilde{x}, \tilde{y} \in \tilde{J}$ and $\tilde{x} < \tilde{y}$ in the ordering induced from $\tilde{\sigma}$ then for all $k > 0$, $f^k(\tilde{x}), f^k(\tilde{y}) \in f^k_\#(\tilde{\sigma})$ and $f^k(\tilde{x}) < f^k(\tilde{y})$ in the ordering induced from $f^k_\#(\tilde{\sigma})$. It follows immediately that $\tilde{\sigma}$ can be $k$-split at $\tilde{J}$ for all $k > 0$, and hence can be split at $\tilde{J}$. Since $\tilde{J}$ is an orbit of $T$, there is an induced splitting $\sigma = \sigma_1$. $\qquad\square$

If $\sigma$ intersects the interior of some edge of $H_i$ then we say that $\sigma$ *intersects* $H_i$ *nontrivially*. We organize splittings of $\sigma$ according to the growth rate of the highest stratum that $\sigma$ intersects nontrivially. After the next definition, we begin with the basic splitting lemma for the non-exponentially-growing case.

*Definition* 4.1.3. Suppose that $f : G \to G$ is a topological representative, that the non-exponentially-growing stratum $H_i$ consists of a single edge $E_i$ and that $f(E_i) = E_i u_i$ for some path $u_i \subset G_{i-1}$. (It is part of the definition of an improved relative train track map that each nonexponentially-growing stratum consists of a single edge.) We say that paths of the form $E_i \gamma \bar{E}_i$, $E_i \gamma$, and $\gamma \bar{E}_i$, where $\gamma \subset G_{i-1}$, are *basic paths of height* $i$.

The restriction on the endpoints of $\sigma$ in the following lemma reduces the number of special cases that we must consider.

LEMMA 4.1.4. *Suppose that $f : G \to G$ and $E_i$ are as in Definition* 4.1.3. *Suppose further that $\sigma \subset G_i$ is a path or circuit that intersects $H_i$ nontrivially and that the endpoints of $\sigma$, if any, are not contained in the interior of $E_i$. Then $\sigma$ has a splitting whose pieces are either basic paths of height $i$ or are contained in $G_{i-1}$.*

*Proof.* Suppose at first that $\sigma$ is a path. Choose lifts $\tilde{f} : \Gamma \to \Gamma$ and $\tilde{\sigma} \subset \Gamma$. Fix $k > 0$. There is an initial segment $E_i^k$ of $E_i$ such that $f^k(E_i^k) = E_i$. No other points in $G_i$ are mapped by $f^k$ into the interior of $E_i$. If a copy of $E_i$ cancels with a copy of $\bar{E}_i$ when $f^k(\sigma)$ is tightened to $f^k_\#(\sigma)$, then there is a subpath $\mu$ in $\sigma$ connecting a copy of $E_i^k$ to a copy of $\bar{E}_i^k$ such that $f^k_\#(\mu) = *$. But $\mu$ is a closed path and $f$ is a homotopy equivalence so this is impossible. We conclude that no such cancellation occurs and hence (Lemma 4.1.1(3)) that $\tilde{\sigma}$ can be $k$-split at any point in the interior of a lift of $E_i^k$ or $\bar{E}_i^k$. Lemma 4.1.1(4) implies that $\tilde{\sigma}$ can be $k$-split at the initial vertex of any lift of $E_i$ and at the terminal vertex of any lift of $\bar{E}_i$. Since $k$ is arbitrary $\tilde{\sigma}$ can be split at these points. Lemma 4.1.1(2) and induction allow us to split at all such points simultaneously. The induced splitting of $\sigma$ satisfies the conclusions of the lemma.

If $\sigma$ is a circuit, first apply Lemma 4.1.2 to obtain a splitting $\sigma = \sigma_1$. If the basepoint of $\sigma_1$ is not contained in the interior of $E_i$ or $\bar{E}_i$, then by our previous argument, $\sigma_1$ has a splitting of the right type. Lemma 4.1.1(1) produces the desired splitting of $\sigma$. Suppose then that the basepoint of $\sigma_1$ is contained in the interior of $E_i$ or $\bar{E}_i$. Arguing as in the previous case, we see that $\tilde{\sigma}_1$ can be split at the initial vertex of any lift of $E_i$ and at the terminal vertex of any lift of $\bar{E}_i$. Let $\sigma_1 = \tau_1 \cdot \ldots \cdot \tau_m$ be the resulting splitting of $\sigma_1$ and let $\tau_1' = \tau_m \tau_1$. Then $\sigma = \tau_1' \cdot \tau_2 \ldots \tau_{m-1}$ is the desired splitting.  □

4.2. *Exponentially-growing strata.* We return to the non-exponentially-growing case in Section 5.5. In this subsection we focus on the exponentially growing case. If $H_r$ is an exponentially-growing stratum, then denote the maximal invariant set $\{x \in H_r : f^k(x) \in H_r, \forall k \geq 0\}$ by $I_r$. The pre-image of $I_r$ in $\Gamma$ is denoted $\tilde{I}_r$. The train track property gives the following basic splitting property.

LEMMA 4.2.1. *Suppose that $f : G \to G$ is a relative train track map, that $H_r$ is an exponentially-growing stratum and that $\sigma \subset G_r$ is an $r$-legal path. If $\tilde{x} \in \tilde{\sigma} \cap \tilde{I}_r$ and if every neighborhood of $\tilde{x}$ in $\tilde{\sigma}$ intersects $\tilde{H}_r$ nontrivially, then $\tilde{\sigma}$ can be split at $\tilde{x}$.*

*Proof.* Let $\tilde{\sigma} = \tilde{\sigma}_1 \tilde{\sigma}_2$ be the decomposition determined by subdividing at $\tilde{x}$. Lemma 4.1.1(3), Lemma 2.5.1 and induction on $k$ imply that $\tilde{\sigma} = \tilde{\sigma}_1 \tilde{\sigma}_2$ is a $k$-splitting for all $k$.  □

Lemma 4.2.2 below is a local version of Lemma 4.2.1. Its proof exploits the fact that exponential growth dominates bounded loss if the initial length is large enough.

Suppose that $f : G \to G$ is a relative train track map, that $H_r$ is an exponentially-growing stratum, that $\tau$ is a path in $G$ and that $\alpha \subset G_r$ is a subpath of $\tau$ with endpoints at vertices. If there are $k$ $H_r$-edges to the left and to the right of $\alpha$ in $\tau$, define $W_k(\alpha)$ to be the subpath of $\tau$ that begins with the $k^{\text{th}}$ $H_r$-edge to the left of $\alpha$ and ends with the $k^{\text{th}}$ $H_r$-edge to the right of $\alpha$. We say that $\alpha$ is *$k$-protected in $\tau$* if its first and last edges are in $H_r$, if $W_k(\alpha) \subset G_r$ and if $W_k(\alpha)$ is $r$-legal.

LEMMA 4.2.2. *Assume that $f : G \to G$ is a relative train track map and that $H_r$ is an exponentially-growing stratum. There is a constant $K$ so that if $\tau$ is a path in $G$ and if $\alpha \subset G_r$ is a $K$-protected subpath of $\tau$, then $\tau$ can be split at the endpoints of $\alpha$.*

*Proof.* Choose $l$ so that the $f^l$-image of an edge in $H_r$ contains at least two edges in $H_r$. Let $K = 2lC$ where $C$ is a bounded cancellation constant for $f : G \to G$ .

We will show that if $\alpha$ is $K$-protected in $\tau$, then $\tau$ can be $i$-split at the endpoints of $\alpha$ for $1 \leq i \leq l$. Moreover, $f_\#^l(\alpha)$ is $K$-protected in $f_\#^l(\tau)$. Iteration of this argument proves that $\tau$ can be $i$-split at the endpoints of $\alpha$ for all $i$.

Lemma 4.1.1(6), Lemma 2.5.1 and the bounded cancellation lemma imply that if $\alpha$ is $k$-protected for $k > C$, then $\tau$ can be 1-split at the endpoints of $\alpha$ and that $f_\#(\alpha)$ is $k - C$ protected. Thus if $\alpha$ is $K$ protected, then $\tau$ can be $i$-split at the endpoints of $\alpha$ and also at the endpoints of $W_{K-lC}(\alpha) = W_{lC}(\alpha)$ for $1 \leq i \leq l$. Since $f_\#^l(W_{lC}(\alpha)) \supset W_{2lC}(f_\#^l(\alpha)) = W_K(f_\#^l(\alpha))$, $f_\#^l(\alpha)$ is $K$-protected in $f_\#^l(\tau)$.                                                                   □

In view of the fact that $\mathcal{B}(G)$ is homeomorphic to $\mathcal{B}$, it is clear what it means for a bi-infinite path in $G$ to be weakly attracted to a generic line of an element of $\mathcal{L}(\mathcal{O})$. We now extend this (and adapt the notation) so that it applies to arbitrary paths in $G$.

*Definition 4.2.3.* Suppose that $f : G \to G$ is a relative train track map, that $H_r$ is an aperiodic exponentially-growing stratum, that $\Lambda^+ \in \mathcal{L}(\mathcal{O})$ is associated to $H_r$ and that $\sigma \subset G_r$ is a path or circuit. Then $\sigma$ *is weakly attracted to* $\Lambda^+$ if each finite subpath of the realization of some (and hence any) generic line of $\Lambda^+$ occurs as an unoriented subpath of $f_\#^k(\sigma)$ for all sufficiently large $k$.

COROLLARY 4.2.4. *Assume that $f : G \to G$ is a relative train track map, that $\Lambda^+$ is the expanding lamination associated to an aperiodic exponentially-growing stratum $H_r$ and that $\sigma \subset G$ is a path or circuit. Then the following are equivalent*:

1. *$\sigma$ is weakly attracted to $\Lambda^+$.*

2. *Some $f_\#^k(\sigma)$ splits into subpaths, at least one of which is weakly attracted to $\Lambda^+$.*

3. *Some $f_\#^k(\sigma)$ splits into subpaths, at least one of which is an edge of $H_r$.*

*Proof.* If $\sigma$ is a circuit then choose $\sigma_1$ as in Lemma 4.1.2. For every finite subpath $\lambda_0$ of a generic line $\lambda$, there is a subpath $\lambda_1$ of $\lambda$ that contains two disjoint copies of $\lambda_0$. If $\lambda_1$ occurs as an unoriented subpath of $f_\#^k \sigma$ then $\lambda_0$ occurs as an unoriented subpath of $f_\#^k \sigma_1$. It follows that $\sigma$ is weakly attracted to $\Lambda^+$ if and only if $\sigma_1$ is weakly attracted to $\Lambda^+$. We may therefore assume that $\sigma$ is a path.

Corollary 3.1.11 and Lemma 3.1.8(3) imply that any edge of $H_r$ is weakly attracted to $\Lambda^+$. Thus (3) implies (2). It is immediate from the definitions of splitting and weak attraction that (2) implies (1). If $\sigma$ is weakly attracted to $\Lambda^+$, then some $f_\#^k(\sigma)$ contains $r$-legal subpaths in $G_r$ that contain arbitrarily many edges in $H_r$. Lemma 4.2.2 implies that $f_\#^k(\sigma)$ splits into subpaths, at least one of which is $r$-legal, is contained in $G_r$ and contains edges of $H_r$. Condition (3) now follows from Lemma 4.2.1.                    □

Lemma 4.2.2 suggests that if $H_r$ is an exponentially-growing stratum, then $\sigma \subset G_r$ can be split into pieces that are either $r$-legal or are neighborhoods of an illegal turn in $H_r$. We make that precise in Lemma 4.2.6 below. First we choose the appropriate neighborhoods of an illegal turn in $H_r$.

If $f : G \to G$ is a relative train track map and $H_r$ is an exponentially-growing stratum, then define $P_r$ to be the set of paths $\rho \subset G_r$ such that:

(i) Each $f_\#^k(\rho)$ contains exactly one illegal turn in $H_r$.

(ii) The initial and terminal (possibly partial) edges of each $f_\#^k(\rho)$ are contained in $H_r$.

(iii) The number of $H_r$-edges in $f_\#^k(\rho)$ is bounded independently of $k$.

The following lemma is essentially contained in the proof of Lemma 5.11 of [BH1]. We give a proof here for the convenience of the reader. Recall (subsection 2.2) that we do not distinguish between a path $\sigma \subset G$ and its associated edge path $E_0' E_1 \ldots E_k'$. In the following proof it is necessary to keep track of those edges in the edge path that are contained in $H_r$. We write $\sigma \cap H_r$ for the ordered sequence of edges and partial edges of $E_0' E_1 \ldots E_k'$ that are contained in $H_r$.

LEMMA 4.2.5.   $P_r$ is a finite $f_\#$-invariant set.

*Proof.* The $f_\#$-invariance of $P_r$ is immediate from the definition. Decompose $\rho \in P_r$ as a concatenation $\rho = \alpha\beta$ where $\alpha$ and $\beta$ are $r$-legal.

When $f^k(\alpha)$ and $f^k(\beta)$ are tightened to $f_\#^k(\alpha)$ and $f_\#^k(\beta)$, no $H_r$-edges are canceled. When $f_\#^k(\alpha)f_\#^k(\beta)$ is tightened to $f_\#^k(\rho)$, an initial segment of $f_\#^k(\bar{\alpha})$ is canceled with an initial segment of $f_\#^k(\beta)$. Since $f_\#^k(\rho)$ has an illegal turn in $H_r$, the first noncanceled edges in $f_\#^k(\bar{\alpha})$ and $f_\#^k(\beta)$ are contained in $H_r$. The cancellation between $f_\#^k(\bar{\alpha})$ and $f_\#^k(\beta)$ is therefore determined by $f_\#^k(\bar{\alpha}) \cap H_r$ and $f_\#^k(\beta) \cap H_r$ : the two paths cancel until the first distinct elements of $f_\#^k(\bar{\alpha}) \cap H_r$ and $f_\#^k(\beta) \cap H_r$ are encountered. Note that $f_\#^k(\bar{\alpha}) \cap H_r$ and $f_\#^k(\beta) \cap H_r$ are determined by $\alpha \cap H_r$ and $\beta \cap H_r$.

We claim that as $\rho$ varies over all elements of $P_r$, $\rho \cap H_r$ takes on only finitely many values. If $\rho$ has a splitting then at least one of the resulting pieces is $r$-legal and intersects $H_r$ nontrivially. This contradicts Lemma 2.5.1 and condition(iii). Lemma 4.2.2 therefore implies that the number of $H_r$ edges in $\rho$ is bounded independently of $\rho$. Let $\alpha_0$ and $\bar{\beta}_0$ be the initial (possibly partial) edges of $\alpha \cap H_r$ and $\bar{\beta} \cap H_r$ respectively. The only possible difficulty is that there might exist $\rho = \alpha\beta$ and $\rho' = \alpha'\beta'$ where the only difference between $\rho \cap H_r$ and $\rho' \cap H_r$ is that the length of $\alpha_0$ and $\alpha'_0$ and the length of $\beta_0$ and $\beta'_0$ may differ. Suppose for concreteness that $\alpha'_0$ is a proper subset of $\alpha_0$ and that $A$ is their difference. Property (ii) implies that the number of edges in $f^k_\#(A) \cap H_r$ grows without bound. Property (iii) therefore implies that for sufficiently large $k$, edges in $f^k_\#(A) \cap H_r$ must be canceled with edges of $f^k_\#(\beta) \cap H_r$. This implies that all of $f^k_\#(\bar{\alpha}') \cap H_r$ is canceled with a proper initial segment $X$ of $f^k_\#(\beta) \cap H_r$. But $X$, like every initial segment of $f^k_\#(\beta) \cap H_r$, either contains or is contained in $f^k_\#(\beta') \cap H_r$. In the former case, all of $f^k_\#(\beta')$ is canceled with part of $f^k_\#(\bar{\alpha}')$; in the latter case all of $f^k_\#(\alpha')$ is canceled with part of $f^k_\#(\bar{\beta}')$. In either case, $f^k_\#(\rho')$ is $r$-legal in contradiction to condition (i). We have now verified the claim.

Property(ii), the fact that $\alpha_0$ and $\beta_0$ take on only finitely many values and the fact that the number of $H_r$-edges in $\alpha$ and in $\beta$ are bounded independently of $\rho$ imply that there exists $k > 0$, independent of $\rho$, such that $f^k_\#(\rho)$ is obtained from $f^k_\#(\alpha_0)$ and $f^k_\#(\beta_0)$ by concatenating and by cancelling at the juncture. This implies that $f^k_\#(\rho)$, and hence $\rho$, takes on only finitely many values.    $\square$

LEMMA 4.2.6.   *Suppose that $f : G \to G$ is a relative train track map, that $H_r$ is an exponentially-growing stratum, that $\sigma \subset G_r$ is a path or circuit and that each $f^k_\#(\sigma)$ has the same finite number of illegal turns in $H_r$. Then $\sigma$ can be split into subpaths that are either $r$-legal or elements of $P_r$.*

*Proof.* We may assume by Lemma 4.1.2 and Lemma 4.1.1(1) that $\sigma$ is a path. We induct on the number $m$ of illegal turns that $\sigma$ has in $H_r$. If $m = 0$ then $\sigma$ is $r$-legal and there is nothing to prove. Suppose that $m = 1$. Write $\sigma = \alpha\beta$ where $\alpha$ and $\beta$ are $r$-legal subpaths. After possibly splitting off an initial segment of $\alpha$ and a terminal segment of $\beta$ according to Lemma 4.2.2, we may assume that $\alpha$ and $\beta$ contain only finitely many $H_r$-edges.

It is convenient to work with lifts $\tilde{f} : \Gamma \to \Gamma$, $\tilde{H}_r$ and $\tilde{\sigma} = \tilde{\alpha}\tilde{\beta}$ to the universal cover. We first check that if $\tilde{\alpha}$ is infinite then there is at least one point in $\tilde{\alpha}$ at which $\tilde{\sigma}$ can be split. For each $k > 0$, $\tilde{f}^k_\#(\tilde{\sigma})$ is obtained from $\tilde{f}^k_\#(\tilde{\alpha})$ and $\tilde{f}^k_\#(\tilde{\beta})$ by concatenating and by cancelling terminal edges of $\tilde{f}^k_\#(\tilde{\alpha})$ with initial edges of $\tilde{f}^k_\#(\tilde{\beta})$. Since each $\tilde{f}^k_\#(\tilde{\sigma})$ has an illegal turn in $\tilde{H}_r$, not all of the $\tilde{H}_r$-edges of $\tilde{f}^k_\#(\tilde{\alpha})$ are canceled during this process. Let $\tilde{x}$ be the initial vertex of the first $\tilde{H}_r$-edge in $\tilde{\sigma}$. Lemma 4.2.1 and Lemma 4.1.1(6) imply that $\tilde{\sigma}$ can be $k$-split at $\tilde{x}$. Since $k$ is arbitrary $\tilde{\sigma}$ can be split at $\tilde{x}$.

The set of points at which $\tilde{\sigma}$ can be split is closed. If $\tilde{\sigma}$ can be split at a point in $\tilde{\alpha}$, choose the splitting point that is closest to the terminal end of $\tilde{\alpha}$. After splitting $\tilde{\sigma}$ at this point we may assume that there are no splitting points for $\tilde{\sigma}$ in $\tilde{\alpha}$. By a completely similar argument we may also assume that there are no splitting points for $\tilde{\sigma}$ in $\tilde{\beta}$. It remains to show that if $\sigma$ has no splittings, then $\sigma$ is an element of $P_r$.

Condition (i) of $P_r$ follows from the hypothesis of the lemma. If condition (iii) is violated, then Lemma 4.2.2 implies that some $f_{\#}^k(\sigma)$ has a splitting. Lemma 4.1.1(5) then implies that $\sigma$ has a splitting, which is a contradiction. The first (possibly partial) edge of $\sigma$ must be contained in $H_r$; otherwise, as above, $\sigma$ can be split at the initial vertex of the first edge of $\alpha$ in $H_r$. The same argument shows that the terminal (possibly partial) edge of $\beta$ is contained in $H_r$. This implies (ii) and completes the proof in the $m = 1$ case.

Suppose now that $m > 1$. Decompose $\tilde{\sigma} = \tilde{\sigma}_1 \ldots \tilde{\sigma}_{m+1}$ so that each juncture is an illegal turn in $\tilde{H}_r$ and each $\tilde{\sigma}_i$ is $r$-legal. Each $\tilde{f}_{\#}^k(\tilde{\sigma})$ has a decomposition $\tilde{f}_{\#}^k(\tilde{\sigma}) = \tilde{\tau}_1^k \ldots \tilde{\tau}_{m+1}^k$ into maximal $r$-legal subpaths. The set $\tilde{S}_k^2 = \{\tilde{x} \in \tilde{\sigma}_2 : \tilde{f}^k(\tilde{x}) \in \tilde{f}_{\#}^k(\tilde{\sigma})\}$ is closed. An easy induction argument shows that $f^N$ maps $\cap_{k=1}^N \tilde{S}_N^2$ onto $\tilde{\tau}_2^N$ for all $N \geq 1$. It follows that $\cap_{k=1}^\infty \tilde{S}_k^2$ is nonempty and that it is therefore possible (Lemma 4.1.1(3)) to split $\tilde{\sigma}$ at a point in $\tilde{\sigma}_2$. This splits $\sigma$ into subpaths that have fewer than $m$ illegal turns in $H_r$ and induction on $m$ completes the proof. □

## 5. Improved relative train track maps

Corollary 3.4.3 connects the Tits alternative to the action of $\mathcal{O}$ on bi-infinite paths and in particular to the basins of attraction for an expanding lamination pair $\Lambda^\pm$. In the next section we state and prove our weak attraction theorem which characterizes these basins of attraction for 'topmost' laminations. In this section, we lay the groundwork for that analysis.

This is the most technical section in the paper. All future references to results in this section will be made to subsection 5.1. It is therefore possible to skim the proof of Theorem 5.1.5 and follow the proof of Theorem 7.0.1.

### 5.1. *Statements.*

*Definition* 5.1.1. A nontrivial path $\sigma \subset G$ is a *periodic Nielsen path* for $f : G \to G$ if $f_{\#}^k(\sigma) = \sigma$ for some $k \geq 1$; if $k = 1$, then we sometimes simply say that $\sigma$ is a *Nielsen path*. We say that the periodic Nielsen path $\sigma$ is *indivisible* if it cannot be written as a concatenation of nontrivial periodic Nielsen paths.

The relative train track maps of [BH1] must be modified to suit our present needs. There, it was sufficient to control the Nielsen paths; in this paper, we must also control the periodic Nielsen paths. We do this by replacing $\mathcal{O}$ by an iterate in which every periodic Nielsen path is a Nielsen path and then applying the techniques of [BH1]. This is carried out in subsection 5.2. [BH1] contains a characterization of the irreducible outer automorphisms that arise as pseudo-Anosov mapping classes. In subsection 5.3, we generalize this by characterizing geometric strata (defined below).

*Remark* 5.1.2. Lemma 4.2.6 implies that if $f : G \to G$ is a relative train track map and $H_r$ is an exponentially-growing stratum, then the indivisible periodic Nielsen paths in $G_r$ that intersect $H_r$ nontrivially are precisely the elements of $P_r$ that have periodic orbit under the action of $f_\#$.

*Definition* 5.1.3. Suppose that $H_i$ is a single edge $E_i$ and that $f(E_i) = E_i \tau^l$ for some closed indivisible Nielsen path $\tau \subset G_{i-1}$ and some $l > 0$. The *exceptional paths of height* $i$ are those paths of the form $E_i \tau^k \bar{E}_j$ or $E_i \bar{\tau}^k \bar{E}_j$ where $k \geq 0$, $j \leq i$, $H_j$ is a single edge $E_j$ and $f(E_j) = E_j \tau^m$ for some $m > 0$. The set of exceptional paths of height $i$ is $f_\#$-invariant. It is an easy consequence of Lemma 4.1.1(3) that an exceptional path of height $i$ has no splittings. If $\sigma$ is an exceptional path of some unspecified height, then we sometimes simply say that $\sigma$ *is an exceptional path*.

*Definition* 5.1.4. Suppose that:

- $g : Q \to Q$ is a homotopy equivalence of a graph $Q$, all of whose components are noncontractible.

- $\alpha_1, \ldots, \alpha_m$ are circuits in $Q$ that are permuted by $g_\#$.

- $S$ is a compact surface with $m+1$ boundary components, $\alpha_1^*, \ldots, \alpha_m^*$ and $\rho^*$.

- $\phi : S \to S$ is a pseudo-Anosov homeomorphism that permutes the $\alpha_i^*$'s in the same way that $g_\#$ permutes the $\alpha_i$'s.

Let $\mathcal{A}$ be the union of $m$ annuli $A_1, \ldots, A_m$. Define $Y$ to be the space obtained from $Q \cup \mathcal{A} \cup S$ by attaching one end of $A_i$ to $\alpha_i$ and the other end to $\alpha_i^*$. Extend $g \cup \phi$ to a homotopy equivalence $h : Y \to Y$ by interpolating between $g(\alpha_i)$ and $\phi(\alpha_i^*)$ on $\mathcal{A}$. We say that $h : Y \to Y$ is *a geometric extension* of $g : Q \to Q$.

Suppose that $f : G \to G$ is a topological representative and that $H_i$ is an exponentially-growing stratum. We say that $H_i$ *is a geometric stratum* if there exists $h : Y \to Y$ as above and a homotopy equivalence $\Phi : (Y, Q) \to$ (noncontractible components of $G_i$, noncontractible components of $G_{i-1}$) such that $f\Phi \simeq \Phi h$. In particular, $\Phi_\#$ identifies the outer automorphism induced by $h$ with the outer automorphism induced by restricting $f$ to the noncontractible components of $G_i$.

Relative train track maps are topological representatives whose exponentially-growing strata satisfy certain extra conditions. We will also add conditions on the zero strata and on the non-exponentially-growing strata. By passing to an iterate and subdividing if necessary (or by adding homological restrictions as in [BFH2]), we may assume that each non-exponentially-growing stratum $H_i$ is a single edge $E_i$ and that $f(E_i) = E_i u_i$ for some path $u_i \subset G_{i-1}$. We introduce a move, called sliding, to arrange that $u_i$ be a closed path and that $f(E_i) = E_i \cdot u_i$. This is carried out in subsection 5.4 and analyzed further in subsection 5.5.

A relative train track map that satisfies the conclusions of the following theorem is said to be an *improved relative train track map*. The proof of Theorem 7.0.1 is given in subsection 5.6.

THEOREM 5.1.5. *For every outer automorphism $\mathcal{O}$ and $\mathcal{O}$-invariant free factor system $\mathcal{F}$ there is an* eg-aperiodic *relative train track map $f : G \to G$ and filtration $\emptyset = G_0 \subset G_1 \subset \cdots \subset G_K = G$ representing an iterate of $\mathcal{O}$ with the following properties.*

- $\mathcal{F} = \mathcal{F}(G_r)$ *for some filtration element $G_r$.*

- $f$ *is reduced (Definition 2.6.6).*

- *Every periodic Nielsen path has period one.*

- *For every vertex $v \in G$, $f(v)$ is a fixed point. If $v$ is an endpoint of an edge in a non-exponentially-growing stratum then $v$ is a fixed point. If $v$ is the endpoint of an edge in an exponentially-growing stratum $H_i$ and if $v$ is also contained in a noncontractible component of $G_{i-1}$, then $v$ is a fixed point.*

- $H_i$ *is a zero stratum if and only if it is the union of the contractible components of $G_i$.*

- *If $H_i$ is a* zero stratum, *then*

  z-(i)  $H_{i+1}$ *is an exponentially-growing stratum.*

  z-(ii)  $f|H_i$ *is an immersion.*

  z-(iii)  *Each vertex in $H_i$ that has valence less than three in $G_{i+1}$ is the endpoint of an edge of $H_{i+1}$.*

- *If $H_i$ is a* non-exponentially-growing stratum, *then*

  ne-(i)  $H_i$ *is a single edge $E_i$.*

  ne-(ii)  $f(E_i) = E_i \cdot u_i$ *for some closed path $u_i \subset G_{i-1}$ whose basepoint is fixed by $f$.*

ne-(iii) *If $\sigma \subset G_i$ is a basic path of height $i$ (Definition 4.1.3) that does not split as a concatenation of two basic paths of height $i$ or as a concatenation of a basic path of height $i$ with a path contained in $G_{i-1}$, then either: (i) some $f_\#^k(\sigma)$ splits into pieces, one of which equals $E_i$ or $\bar{E}_i$; or (ii) $u_i$ is a Nielsen path and some $f_\#^k(\sigma)$ is an exceptional path of height $i$.*

- *If $H_i$ is an* exponentially-growing *stratum then*

eg-(i) *There is at most one indivisible Nielsen path $\rho_i \subset G_i$ that intersects $H_i$ nontrivially. The initial edges of $\rho_i$ and $\bar{\rho}_i$ are distinct (possibly partial) edges in $H_i$.*

eg-(ii) *If $\rho_i \subset G_i$ is an indivisible Nielsen path that intersects $H_i$ nontrivially and if $H_i$ is not geometric, then there is an edge $E$ of $H_i$ that $\rho_i$ crosses exactly once. (See also Lemma 5.1.7 below.)*

eg-(iii) *If $H_i$ is geometric then there is an indivisible Nielsen path $\rho_i \subset G_i$ that intersects $H_i$ nontrivially and satisfies the following properties: (i) $\rho_i$ is a closed path with basepoint in the interior of $H_i$; (ii) the circuit determined by $\rho_i$ corresponds to the unattached peripheral curve $\rho^*$ of $S$; and (iii) the surface $S$ is connected.*

Lemma 5.1.7 below is used to analyze nongeometric exponentially-growing strata.

*Definition* 5.1.6. For any subgraph $X$ of $G$ and finite path $\rho \subset G$, define $\langle X, \rho \rangle$ to be the groupoid of paths that can be decomposed into a concatenation of subpaths that are either entirely contained in $X$ or are equal to $\rho$ or $\bar{\rho}$.

LEMMA 5.1.7. *Suppose that $f : G \to G$ is reduced, that $H_r$ is an aperiodic exponentially-growing stratum, that $\rho_r \subset G_r$ is a Nielsen path that crosses some edge $E$ of $H_r$ exactly once and that the first and last (possibly partial) edges of $\rho_r$ are contained in $H_r$. Then the endpoints of $\rho_r$ are distinct and if both endpoints are contained in $G_{r-1}$, then at least one of them is contained in a contractible component of $G_{r-1}$. If $X$ is a subgraph of $G$ that does not contain any edges of $H_r$, then there is a free factor system that carries the same bi-infinite paths as $\langle X, \rho_r \rangle$.*

*Proof.* We may assume, after subdividing if necessary, that the endpoints of $\rho_r$ are vertices. Let $\hat{G}$ be the graph obtained from $G$ by removing the edge $E$ and adding a new edge $\hat{E}$ with endpoints equal to the initial and terminal endpoints of $\rho$. Decompose $\rho$ into the concatenation of subpaths $\rho = \alpha E \beta$ where $\alpha$ and $\beta$ are disjoint from $E$. There is a homotopy equivalence $h : G \to \hat{G}$ that is the 'identity' on all edges other than $E$ and that satisfies $h(E) = \bar{\alpha}\hat{E}\bar{\beta}$. (The homotopy inverse sends $\hat{E}$ to $\rho_r$.) The map $h_\#$ induces a bijection between the bi-infinite paths in $\langle X, \rho_r \rangle$ and the bi-infinite paths in the subgraph $X \cup \hat{E}$ of $\hat{G}$. In particular, $\langle X, \rho_r \rangle$ carries exactly the same bi-infinite paths as the free factor system $\mathcal{F}(X \cup \hat{E})$.

Suppose now that $X = G_{r-1}$. If the endpoints of $\rho_r$ are equal or are both contained in noncontractible components of $G_{r-1}$, then $\mathcal{F}(X \cup \hat{E})$ is strictly larger than $\mathcal{F}(X) = \mathcal{F}(G_{r-1})$. An $r$-legal bi-infinite path in $\langle X, \rho_r \rangle$ must lie entirely in $X$. Thus $\mathcal{F}(X \cup \hat{E})$ does not carry any line that is generic for the element of $\mathcal{L}(\mathcal{O})$ associated to $H_r$ and so is strictly smaller than $\mathcal{F}(G_r)$. Since $\langle X, \rho_r \rangle$ is $f_\#$-invariant, this contradicts the assumption that $f : G \to G$ is reduced and completes the proof.     □

An outer automorphism $\mathcal{O}$ is said to have *polynomial growth* if some (and hence every) relative train track map  representing $\mathcal{O}$ has no exponentially-growing strata.  The set of all polynomial growth outer automorphisms is denoted $\mathrm{PG}(F_n)$. An element of $\mathrm{GL}(n, \mathbb{Z})$ is *unipotent* if it is conjugate to an upper triangular matrix with ones on the diagonal. We say that an element of $\mathrm{PG}(F_n)$ is *unipotent* if its image in $\mathrm{GL}(n, \mathbb{Z})$ is unipotent. The set of unipotent elements of $\mathrm{PG}(F_n)$ is denoted $\mathrm{UPG}(F_n)$ and plays a central role in [BFH2] and [BFH1]. We conclude this subsection with a strengthening of Theorem 7.0.1 for elements of $\mathrm{UPG}(F_n)$.

THEOREM 5.1.8. *Suppose that $\mathcal{O} \in \mathrm{UPG}(F_n)$ and that $\mathcal{F}$ is an $\mathcal{O}$-invariant free factor system. Then there is a relative train track map $f : G \to G$ and filtration $\emptyset = G_0 \subset G_1 \subset \cdots \subset G_K = G$ representing $\mathcal{O}$ with the following properties.*

1. *$\mathcal{F} = \mathcal{F}(G_r)$ for some filtration element $G_r$.*

2. *Each $H_i$ is a single edge $E_i$ satisfying $f(E_i) = E_i \cdot u_i$ for some closed path $u_i \subset G_{i-1}$.*

3. *Every vertex of $G$ is fixed by $f$.*

4. *Every periodic Nielsen path has period one.*

5. *If $\sigma$ is any path with endpoints at vertices, then there exists $M = M(\sigma)$ so that for each $m \geq M$, $f_\#^m(\sigma)$ splits into subpaths that are either single edges or are exceptional.*

6. *$M(\sigma)$ is a bounded multiple of the edge length of $\sigma$.*

The proof of Theorem 5.1.8 is given in subsection 5.7.

5.2. *Nielsen paths in exponentially-growing strata.* This subsection and the one to follow focus on exponentially-growing strata.

*Definition* 5.2.1. Suppose that $\mathcal{O}$ fixes each element of $\mathcal{L}(\mathcal{O}) = \{\Lambda_1^+, \ldots, \Lambda_l^+\}$, that $\mu_i$ is the expansion factor for the action of $\mathcal{O}$ on $\Lambda_i^+$ and that $\mathcal{F}$ is an $\mathcal{O}$-invariant free factor system. We say that a topological representative $f : G \to G$ of $\mathcal{O}$ and filtration $\emptyset = G_0 \subset G_1 \subset \cdots \subset G_K = G$ are $\mathcal{F}$-*Nielsen minimized* if:

(1) $f : G \to G$ is reduced.

(2) $\mathcal{F} = \mathcal{F}(G_s)$ for some filtration element $G_s$.

(3) There are exactly $l$ exponentially-growing strata and the Perron-Frobenius eigenvalues of their transition submatrices equal $\mu_1, \ldots, \mu_l$.

(4) If $H_r$ is an exponentially-growing stratum, then every indivisible periodic Nielsen path $\rho \subset G_r$ that intersects $H_r$ nontrivially has period one.

(5) If $C$ is a contractible component of some $G_i$, then $f^j(C) \subset G_{i-1}$ for some $j > 0$.

(6) For each exponentially-growing stratum $H_r$, let $N_r(f)$ be the number of indivisible Nielsen paths $\rho \subset G_r$ that intersect $H_r$ nontrivially. Then $N(f) = \sum_r N_r(f)$ is as small as possible subject to conditions (1)–(5).

*Remark* 5.2.2. Proposition 3.3.3(4) implies that condition (3) is satisfied by every relative train track map. In the course of proving Lemma 5.2.5 below, we must consider topological representatives that are not relative train track maps. That is why the definition does not include the hypothesis that $f : G \to G$ is a relative train track map and why (4) and (6) do not refer to $P_r$.

One might expect to leave (1), (2) and (5) out of the definition and add them later as separate conditions. In that case the minimization of $N(f)$ in (6) would take place over a larger collection of topological representatives and it is an *a priori* possibility that the absolute minimum $N(f)$ would not occur for a relative train track map satisfying (1), (2) and (5).

If we eliminate conditions (1), (2), (4) and (5) from the definition then we recover the definition of stable topological representative from page 42 of [BH1].

LEMMA 5.2.3. *For any $\mathcal{O}$ and $\mathcal{F}$, there exists an $\mathcal{F}$-Nielsen minimized relative train track map $f : G \to G$ and filtration representing $\mathcal{O}^k$ for some $k > 0$.*

*Proof.* After we replace $\mathcal{O}$ by an iterate if necessary, $\mathcal{O}$ fixes each element of $\mathcal{L}(\mathcal{O})$. Lemma 2.6.7 produces a relative train track map and filtration that represents some $\mathcal{O}^k$ and that satisfies (1), (2) and (5) in the definition of $\mathcal{F}$-Nielsen minimized. Lemma 3.1.14 implies that $f : G \to G$ is eg-aperiodic. Proposition 3.3.3(4) implies that $f : G \to G$ satisfies (3). Since (Remark 5.1.2) every indivisible periodic Nielsen path is an element of the finite set $P_r$, we may assume, by increasing $k$, that (4) is satisfied. Fix $k$, and choose, among all topological representatives for $\mathcal{O}^k$ that satisfy (1)–(5), one, say $f' : G' \to G'$, that minimizes $N(f)$ and is hence $\mathcal{F}$-Nielsen minimizing. If $f' : G' \to G'$ is a relative train track map, we are done. If not, modify $f' : G' \to G'$ by performing 'core subdivisions' and 'collapsing of inessential paths' to construct a relative train track map $f : G \to G$. This is fully described in the proof of Lemma 5.13, Lemma 5.14 and Theorem 5.12 of [BH1]. Since $f : G \to G$ is a relative train track map, (3) is satisfied. The construction involves only subdivision, folding and the collapse of pre-trivial forests. As discussed in the proof of Lemma 2.6.7, (1) and (2) are satisfied. These operations preserve (5) and do not change the period of any indivisible periodic Nielsen path nor any $N_r(f)$. Thus $f : G \to G$ is still $\mathcal{F}$-Nielsen minimized. $\qquad\square$

Suppose that $f : G \to G$ is a relative train track map and that $H_r$ is an exponentially-growing stratum. If $\rho_r \subset G_r$ is an indivisible Nielsen path that intersects $H_r$ nontrivially, then $\rho_r = \alpha\beta$ where $\alpha$ and $\beta$ are $r$-legal and the turn at the juncture of $\bar{\alpha}$ and $\beta$ is an illegal turn in $H_r$. We say that the *fold at the illegal turn of $\rho_r$ in $H_r$ is a full fold* if either all of the initial (possibly partial) edge $E_1$ of $\bar{\alpha}$ can be folded with all or part of the initial (possibly partial) edge $E_2$ of $\beta$ or all of $E_2$ can be folded with all or part of $E_1$.

LEMMA 5.2.4. *Suppose that $f : G \to G$ is an $\mathcal{F}$-Nielsen minimized relative train track map, that $H_r$ is an exponentially-growing stratum and that $\rho_r \subset G_r$ is an indivisible Nielsen path that intersects $H_r$ nontrivially. Then the fold at the illegal turn of $\rho_r$ in $H_r$ is a full fold.*

*Proof.* This is a slight modification of Lemma 5.17 of [BH1]. There is a decomposition $\rho_r = \alpha\beta$ into $r$-legal paths in $G_r$ and there is a path $\tau \subset G_r$ such that $f_\#(\alpha) = \alpha\tau$ and $f_\#(\bar{\beta}) = \bar{\beta}\tau$. We may assume, after subdividing if necessary, that the endpoints of $\rho_r$ are vertices. It suffices (page 25 of [BH1]) to show that $\alpha$ and $\beta$ cannot both be single edges.

Suppose to the contrary that both $\alpha$ and $\beta$ are single edges. Let $G'$ be the graph obtained from $G$ by identifying $\alpha$ and $\bar{\beta}$ and let $q : G \to G'$ be the quotient map. Since $f(\alpha) = \alpha\tau$ and $f(\bar{\beta}) = \bar{\beta}\tau$, there is an induced map $f' : G' \to G'$ defined by $f'(q(E)) = q_\# f(E)$ for each edge $E$ of $G$. If there are edges with trivial $f'$-image, then they form a tree and we collapse each component of the tree. After repeating this tighten and collapse procedure finitely many times (cf. subsection 2.4) we arrive at a topological representative we continue to call $f' : G' \to G'$ and a quotient map we continue to call $q : G \to G'$. An edge in $G$ is collapsed if and only if some iterate of $f_\#$ maps it to $\rho_r$ or $\bar{\rho}_r$. In particular, only edges in zero strata can be collapsed. We claim that conditions (1)–(5) in the definition of $\mathcal{F}$-Nielsen minimized are satisfied by $f' : G' \to G'$ and the filtration with elements of the form $G_i' = q(G_i)$. (If each component of $H_j$ is collapsed to a point then $q(G_j)$ is not added to the filtration.)

Lemma 5.1.7 implies that the endpoints of $\rho_r$ are distinct and that at least one of them is not contained in a noncontractible component of $G_{r-1}$. Thus $q$ is a homotopy equivalence (cf. subsection 2.4) and $\mathcal{F}(G_i) = \mathcal{F}(G_i')$ for all $i$. This implies conditions (1) and (2).

It is easy to check that condition (5) is stable under these collapse and tighten operations and we leave this to the reader.

If $E$ is an edge of $G_r$, then $q(E)$ is an edge of $G'$ and $f'(q(E)) = qf(E)$. (If $E$ belongs to $H_r$, then this uses the fact that an $r$-legal path in $G_r$ does not cross the turn $(\bar{\alpha}, \beta)$.) If $E$ is an edge of $G \setminus G_r$ that is not collapsed by $q$, then $f'q(E)$ is obtained from $qf(E)$ by cancelling edges in $H_r$. Thus $H_i'= q(H_i)$ is exponentially-growing if and only if $H_i$ is exponentially-growing and the Perron-Frobenius eigenvalues for $M_i$ and $M_i'$ are equal. This implies that condition (3) holds.

For every path $\sigma \subset G$ and $k > 0$, $(f')_\#^k(q_\#(\sigma)) = q_\# f_\#^k(\sigma)$. In particular, if $\sigma$ is a periodic Nielsen path for $f$, then $\sigma' = q_\#(\sigma)$ is a periodic Nielsen path for $f'$ and the period of $\sigma'$ is at most the period of $\sigma$. If $\sigma \neq \rho_r$, then $\sigma'$ is not trivial.

Conversely, suppose that $\sigma' \subset G_i'$ is a periodic Nielsen path for $f' : G' \to G'$. We choose a path $\sigma \subset G_i$ satisfying $q_\#(\sigma) = \sigma'$ as follows. If the endpoints of $\sigma'$ do not lie in the edge $q(\alpha) = q(\bar{\beta})$, then there is a unique path $\sigma \subset G$ satisfying $q_\#(\sigma) = \sigma'$. If an endpoint of $\sigma'$ lies in $q(\alpha) = q(\bar{\beta})$ but is distinct from the initial endpoint of $q(\alpha) = q(\bar{\beta})$, then there is a unique path $\sigma \subset G$ that has periodic endpoints and that satisfies $q_\#(\sigma) = \sigma'$. Finally, if $\sigma'$ begins or ends at the initial endpoint of $q(\alpha) = q(\bar{\beta})$, then there is a unique path $\sigma \subset G$ that does not begin or end with $\rho_r$ or $\bar{\rho}_r$ and that satisfies $q_\#(\sigma) = \sigma'$. In all cases, $\sigma$ is a periodic Nielsen path and the period of $\sigma'$ equals the period of $\sigma$; moreover, if $\sigma'$ is indivisible, then $\sigma$ is indivisible. Condition (4) for $f'$ now follows from condition (4) for $f$.

Since $q_\#(\rho_r)$ is trivial, we have decreased $N(f)$. This contradiction completes the proof.     □

The following lemma is the main result of this subsection.

LEMMA 5.2.5.   *If $f : G \to G$ is an $\mathcal{F}$-Nielsen minimized relative train track map and $H_r$ is an exponentially-growing stratum, then there exists at most one indivisible Nielsen path $\rho_r \subset G_r$ that intersects $H_r$ nontrivially. Moreover, if there is such an indivisible Nielsen path $\rho_r$, then:*

- *The first and last (possibly partial) edges of $\rho_r$ are contained in $H_r$.*

- *The illegal turn of $\rho_r$ in $H_r$ is the only illegal turn in $H_r$.*

- *$\rho_r$ crosses every edge in $H_r$ at least once.*

- *Either $\rho_r$ crosses every edge of $H_r$ exactly twice or $\rho_r$ crosses some edge of $H_r$ exactly once.*

*Proof.*  This is Theorem 5.15 of [BH1] with the word stable replaced by $\mathcal{F}$-Nielsen minimized. Having proved Lemma 5.2.4, we see that the proof given in [BH1] carries over to this context without change.    □

5.3. *Geometric strata.*  The following proposition generalizes Theorem 4.1 of [BH1]. The entire subsection is devoted to its proof.

PROPOSITION 5.3.1.   *Suppose that $f : G \to G$ is an $\mathcal{F}$-Nielsen minimized relative train track map with filtration $\emptyset = G_0 \subset G_1 \subset \cdots \subset G_K = G$, that $H_r$ is an exponentially-growing stratum and that $\rho_r \subset G_r$ is an indivisible Nielsen path that crosses every edge of $H_r$ exactly twice. Then:*

- *The endpoints of $\rho_r$ are equal and are not contained in $G_{r-1}$.*

- *The initial (possibly partial) edges of $\rho_r$ and $\bar{\rho}_r$ are distinct.*

- *$H_r$ is a geometric stratum. Moreover, the associated surface $S$ is connected and its unattached peripheral curve $\rho^*$ corresponds to the circuit determined by $\rho_r$.*

We assume throughout this subsection that $K = r$ and that the endpoints of $\rho_r$ are vertices. This causes no loss of generality in our proof of Proposition 5.3.1.

We begin by recalling a pair of definitions from page 46 of [BH1].

*Definition* 5.3.2 (folding $\rho_r$: the proper case). Suppose that $f : G \to G$, $H_r$ and $\rho_r$ are as in the hypotheses of Proposition 5.3.1. Decompose $\rho_r = \alpha\beta$ into a concatenation of maximal $r$-legal subpaths and let $E_1 \subset H_r$ and $E_2 \subset H_r$ be the initial edges of $\bar{\alpha}$ and $\beta$ respectively. Lemma 5.2.4 implies that one of the edge paths $f(E_i)$, $i = 1$ or 2, is an initial subpath of the other. For concreteness, suppose that $f(E_1)$ is an initial subpath of $f(E_2)$. Assume that

$f(E_1)$ is a proper subpath of $f(E_2)$. (The case that $f(E_1) = f(E_2)$ is handled in Definition 5.3.4.)

Let $b$ be the (possibly trivial) maximal subpath of $G_{r-1}$ that follows $E_1$ in $\bar{\alpha}$. Lemma 2.5.1 implies that $f(E_1)f_\#(b)$ is an initial segment of $f_\#(\bar{\alpha})$ that is followed in $f_\#(\bar{\alpha})$ by an edge in $H_r$. Since $f : G \to G$ is a relative train track map , $f_\#(b)$ is nontrivial whenever $b$ is nontrivial. The last edge of $f(E_2)$ and the first edge of $f_\#(\bar{\alpha})$ that does not cancel with an edge of $f_\#(\beta)$ are contained in $H_r$. Thus $f_\#(E_1b)$ is a proper initial segment of $f(E_2)$. Define $F : G \to G'$ to be the generalized fold (see subsection 2.4) of $E_1b$ with the corresponding proper initial subpath of $E_2$. There is a map $g : G' \to G$ such that $gF \simeq f$ rel $\mathcal{V}$. We refer to $F : G \to G'$ as the *extended fold* (*determined by $\rho_r$*) and to $g : G' \to G$ as a *map induced by the extended fold*. (Thus the extended fold determined by $\rho_r$ is a particular type of generalized fold; not every generalized fold is the extended fold for some $\rho_r$.)

Since $f$ is a relative train track map and $H_r$ is the highest stratum, $g(E')$ does not cross the turn $(\bar{E}_1, E_2)$ for any edge $E'$ of $G'$. Thus $f' = Fg : G' \to G'$ is a topological representative. The filtration for $f' : G' \to G'$ is defined by $H'_i = H_i$ for $i < r$ and $H'_r = (H_r \setminus E_2) \cup E'_2$. We say that $f' : G' \to G'$ is obtained from $f : G \to G$ by *folding $\rho_r$* and that $\rho'_r = F_\#(\rho_r)$ is *the indivisible Nielsen path determined by $\rho_r$*.

The following lemma states, among other things, that $f' : G' \to G'$ is a relative train track map. This would fail if we simply folded $E_1$ with an initial segment of $E_2$.

LEMMA 5.3.3.   *With notation as in Definition* 5.3.2,

- *$f' : G' \to G'$, $H'_r$, $\rho'_r$ and $G'_{r-1}$ satisfy the hypotheses of Proposition* 5.3.1.

- *If $f' : G' \to G'$, $H'_r$, $\rho'_r$ and $G'_{r-1}$ satisfy the conclusions of Proposition* 5.3.1 *then $f : G \to G$, $H_r$, $\rho_r$ and $G_{r-1}$ satisfy the conclusions of Proposition* 5.3.1.

- *$H_r$ and $H'_r$ have the same number of edges.*

*Proof.* A proof that $f' : G' \to G'$ is a relative train track map and that $H'_r$ is an exponentially-growing stratum is contained in the proof of Theorem 5.15 of [BH1]. As in the proof of Lemma 5.2.4, $f' : G' \to G'$ is $\mathcal{F}$-Nielsen minimizing.

The extended fold $F : G \to G'$ can be described as follows. $E_2$ is subdivided into two pieces if $b$ is trivial and three pieces if $b$ is nontrivial. The first piece is labeled $E_1$ and then identified with $E_1 \subset H_r$; the middle piece, if it exists, is labeled $b$ and then identified with $b \subset G_{r-1}$; and the last piece is labeled $E'_2$ and is the new edge of $H'_r$. To construct $\rho'_r$, subdivide and relabel

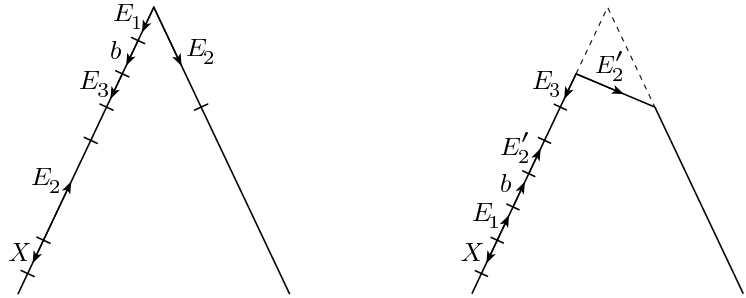

*both* copies of $E_2$ in $\rho_r$. Remove $E_1$ and $b$ (if it is nontrivial) from $\alpha$ to form $\alpha'$ and remove the first and middle segment (if it exists) of the subdivided $E_2$ that is the first edge of $\beta$ to form $\beta'$; $\rho'_r = \alpha'\beta'$. The key point is that the second copy of $E_2$ in $\rho_r$ contributes to $\rho'_r$ a copy of the edges that are removed from $\alpha$ to form $\alpha'$. The losses and gains exactly balance so that $\rho'_r$ crosses every $H'_r$ edge exactly twice. The completes the proof of the first part of the lemma.

Since $F$ : (noncontractible components of $G_r$, noncontractible components of $G_{r-1}$) → (noncontractible components of $G'_r$, noncontractible components of $G'_{r-1}$) is a homotopy equivalence, $H_r$ is a geometric stratum if and only if $H'_r$ is a geometric stratum. By construction, $F_\#$ maps the circuit determined by $\rho_r$ to the circuit determined by $\rho'$. If the initial endpoint $v_1$ of $\rho_r$ is not contained in $G_{r-1}$, then each edge in the link of $v_i$ is in $H_r$. Each time that the interior of $\rho_r$ passes through $v_1$, it crosses two edges of $H_r$. The total number of times that $\rho_r$ crosses the $H_r$-edges in the link of $v_1$ is even. Since $\rho_r$ starts at $v_1$, it must also end at $v_1$. Thus the conclusion that the endpoints of $\rho_r$ are equal is a consequence of the other conclusions of Proposition 5.3.1. If an endpoint of $\rho_r$ is in $G_{r-1}$ or if the initial edges of $\rho_r$ and $\bar{\rho}_r$ are equal then the same is true for $F_\#(\rho_r) = \rho'_r$. This completes the second part of the lemma.

It is clear from the definitions that $H_r$ and $H'_r$ have the same number of edges.  □

*Definition* 5.3.4 (folding $\rho_r$: the improper case). Suppose that $f : G \to G$, $H_r$ and $\rho_r$ are as in Proposition 5.3.1. Decompose $\rho_r = \alpha\beta$ into a concatenation of maximal $r$-legal subpaths and let $E_1 \subset H_r$ and $E_2 \subset H_r$ be the initial edges of $\bar{\alpha}$ and $\beta$ respectively. Assume that $f(E_1) = f(E_2)$.

Define $F : G \to G'$ to be the fold of $E_1$ with $E_2$. There is an induced map $g : G' \to G$ satisfying $gF = f$. As in the previous case, $f' = Fg : G' \to G'$ is a topological representative. We may think of $G_{r-1}$ as a subgraph of $G'$. The filtration for $f' : G' \to G'$ is defined by $G'_i = G_i$ for $i < r$ and $G'_r = G'$. If $f' : G' \to G'$ is a relative train track map then let $\rho'_r = F_\#(\rho_r)$. If $f' : G' \to G'$ is not a relative train track map, then restore the relative train track property by collapsing inessential connecting paths and by performing core subdivisions in $G_{r-1}$ as in the proof of Lemma 5.2.3. This process may change the combinatorial type of $G'_{r-1}$ but an edge of $H'_r = F(H_r)$ is only affected by being shortened according to a core subdivision. We abuse notation by denoting the resulting relative train track map by $f' : G' \to G'$ and the top stratum by $H'_r$ even though $G'$ may have changed. $F_\#(\rho_r)$ determines an indivisible Nielsen path $\rho'_r$ that intersects $H'_r$ nontrivially.

LEMMA 5.3.5.    *With notation as in Definition* 5.3.4,

- $f' : G' \to G'$, $H'_r$, $\rho'_r$ *and* $G'_{r-1}$ *satisfy the hypotheses of Proposition* 5.3.1.

- *If* $f' : G' \to G'$,   $H'_r$, $\rho'_r$ *and* $G'_{r-1}$ *satisfy the conclusions of Proposition* 5.3.1 *then* $f : G \to G$, $H_r$, $\rho_r$ *and* $G_{r-1}$ *satisfy the conclusions of Proposition* 5.3.1.

- $H'_r$ *has one less edge than* $H_r$ *does.*

*Proof.* By construction $f' : G' \to G'$ is a relative train track map. As in the proof of Lemma 5.2.3, $f' : G' \to G'$ is $\mathcal{F}$-Nielsen minimizing. Collapsing inessential connecting paths and performing core subdivisions in $G_{r-1}$ has no effect on the way that $\rho'_r$ crosses edges in $H'_r$. Thus the argument used in Lemma 5.3.3 to prove that $\rho'_r$ crosses each edge of $H'_r$ exactly twice applies in this context as well. This completes the proof of the first part of the lemma.

The second part is proved exactly as it was in Lemma 5.3.3. The third part is immediate from the construction. $\qquad\square$

While proving Proposition 5.3.1, there is no loss in replacing $f : G \to G$ by $f' : G' \to G'$ produced by either a proper fold of $\rho_r$ (Definition 5.3.2) or by an improper fold of $\rho_r$ (Definition 5.3.4). Moreover, this process can be repeated by folding $\rho'_r$ and so on. We refer to this as *repeatedly folding the indivisible Nielsen path* . If at some point the fold is improper, then the number of edges in the $r$-stratum is reduced by the folding process. Since the number of edges in the $r$-stratum never increases, improper folds occur only finitely many times. *We may therefore assume that as we repeatedly fold the indivisible Nielsen path, each fold is proper.*

Every topological representative of $G$ factors (subsection 2.4) as a sequence of folds followed by a homeomorphism. In general, there is no preferred way to choose the folds (although there is an *a priori* bound on the number of folds that occur) but in our context it is natural to begin at the illegal turn of $\rho_r$ in $H_r$. This is the key to Lemma 5.3.6 below. To make its statement precise, note that for any extended fold $F : G \to G'$, we may think of $G_{r-1}$ as a subgraph of both $G$ and $G'$ and with respect to this notation, $F|G_{r-1} =$ identity.

LEMMA 5.3.6.    *If* $f : G \to G$ *and* $\rho_r \subset G_r$ *satisfy the hypotheses of Proposition* 5.3.1, *then there exist*:

- *a composition of extended folds $f_r : G \to G^1$,*

- *a composition $f_{r-1} : G^1 \to G^2$ of folds involving edges in $G_{r-1}$,*

- *a homeomorphism $\theta : G^2 \to G$,*

*such that $f \simeq \theta f_{r-1} f_r$ rel $\mathcal{V}$.*

*Proof.* We use the notation of Definition 5.3.2; we also let $\rho'_r = \alpha' \beta' \subset G'_r$ be the decomposition of $\rho'_r$ into maximal $r$-legal subpaths. Let $E'_1$ and $E'_2$ be the initial edges of $\bar{\alpha}'$ and $\beta'$ respectively and let $F' : G' \to G''$ be the extended fold of $\rho'_r$ with respect to $f' : G' \to G'$. We claim that either $g$ cannot be folded at $(E'_1, E'_2)$ or $F' : G' \to G''$ is a generalized fold for $g : G' \to G$. (By construction, $F' : G' \to G''$ is a generalized fold for $f' : G' \to G'$.) More precisely, write $E'_2$ as a concatenation of subpaths $E'_2 = \mu'_1 \mu'_2 \mu'_3$ satisfying $f'(\mu'_1) = f'(E'_1)$ and $f'(\mu'_2) = f'_\#(b')$ where $b'$ is the maximal subpath of $G'_{r-1}$ following $E'_1$ in $\rho'_r$. (It is possible that both $\mu'_2$ and $b'$ are trivial.) We will show that if $g(E'_1)$ and $g(E'_2)$ have a nontrivial common initial segment, then $g(\mu'_1) = g(E'_1)$ and $g'(\mu'_2) = g'_\#(b')$.

As a first step toward verifying the claim, suppose that $g(E'_1)$ and $g(E'_2)$ have a common initial segment but that the maximal $g$-fold of $E'_1$ and $E'_2$ does not use all of $\mu'_1$. In other words, suppose that the $g$-fold is not full. Let $\hat{F} : G' \to \hat{G}$ be the maximal $g$-fold of $E'_1$ and $E'_2$, let $\hat{g} : \hat{G} \to G$ be the induced map satisfying $\hat{g}\hat{F} = g$ and let $\hat{G}_i = \hat{F}(G'_i)$ for $1 \le i \le r$. Since the fold is not full, $\hat{H}_r = \hat{F}(H'_r)$ has one more edge than do $H'_r$ and $H_r$.

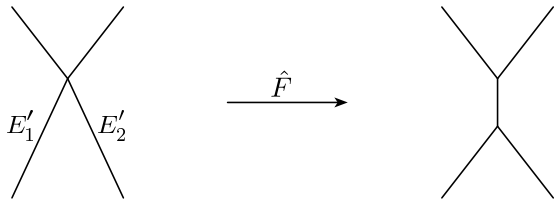

Since $\hat{F}$ is maximal, $\hat{g}$ cannot be folded at the newly created vertex. By Lemma 5.2.5, $(E'_1, E'_2)$ is the only illegal turn for $f'$ that involves an edge in $H'_r$. Since every fold for $g$ is also a fold for $f' = Fg$, $F'$ is the only fold for $g$ that involves an edge in $H'_r$. It follows that all folds for $\hat{g}$ have both edges in $\hat{G}_{r-1}$. If $\hat{\gamma} \subset \hat{G}_{r-1}$ is any nontrivial path with endpoints in $\hat{H}_r \cap \hat{G}_{r-1} = H'_r \cap G'_{r-1} = H_r \cap G_{r-1}$, then there exists a nontrivial path $\gamma \subset G_r$ with endpoints equal to those of $\hat{\gamma}$ such that $\hat{F}_\# F_\#(\gamma) = \hat{\gamma}$; in particular, $\hat{g}_\#(\hat{\gamma}) = f_\#(\gamma)$ is nontrivial. It follows that folding edges in $\hat{G}_{r-1}$ according to $\hat{g}|\hat{G}_{r-1}$ will not cause previously distinct vertices in $\hat{H}_r \cap \hat{G}_{r-1}$ to become identified. Thus no new illegal turns involving the edges of $\hat{H}_r$ are created during such folds. Consecutively fold edges in $\hat{G}_{r-1}$ according to $\hat{g}|\hat{G}_{r-1}$ until no more folds are possible; call this composition of folds $f_{r-1} : \hat{G} \to G^*$. There is an induced immersion $\theta : G^* \to G$ satisfying $f = \theta f_{r-1} \hat{F} F$. Since $G$ has no valence-one vertices, $\theta$ is a homeomorphism. Now $\theta(f_{r-1}(\hat{G}_{r-1})) = f(G_{r-1}) = G_{r-1}$ so that $\theta(f_{r-1}(\hat{H}_r)) = H_r$. But $\theta f_{r-1} | \hat{H}_r$ induces a bijection on edges and so this contradicts the fact that $\hat{H}_r$ and $H_r$ do not have the same number of edges. We conclude that if $g$ can be folded at $(E'_1, E'_2)$, then all of $E'_1$ can be folded with a proper initial segment of $E'_2$.

Let $E$ be the first edge of $f_\#(\bar{\beta})$ that is not part of the maximum common initial segment of $f_\#(\bar{\alpha})$ and $f_\#(\beta)$. Then the initial edge $E'$ in $F(E) \subset G'$ is the first edge of $g_\#(\beta')$ that is not part of the maximum common initial segment of $g_\#(\bar{\alpha}')$ and $g_\#(\beta')$. Since $E$ is contained in $H_r$, $E'$ is contained in $H'_r$. In particular, $E'$ is not contained in $b'$ and we have verified the claim.

If $g$ cannot be folded at $(\bar{E}'_1, E'_2)$, then define $f_r = F$ and construct $f_{r-1}$ and $\theta$ exactly as above. Otherwise, let $F' : G' \to G''$ be the extended fold of $\rho'_r$ with respect to $f' : G' \to G'$ and let $g' : G'' \to G$ be the induced map satisfying $g'F' \simeq g$ rel $\mathcal{V}'$. If $g'$ cannot be folded at the illegal turn of $\rho''_r$, then define $f_r = F'F$ and construct $f_{r-1}$ and $\theta$ exactly as above. Otherwise repeat the argument of the claim to conclude that the extended fold of $\rho''_r$ is a generalized fold for $g'$.

Continue in this manner until we arrive at the desired factorization. The process must terminate because there is a bound to the number of folds that occur in any factorization of $f : G \to G$ .      $\square$

We associate a surface $S$ and a graph $K$ to $f : G \to G$ and $\rho_r$ as follows. We think of $\rho_r$ as a map with domain $I \times \{0\} \subset I \times [0, 1]$ and subdivide $I \times \{0\}$ into subintervals that map either to individual edges of $H_r$ or into maximal subpaths $\{b_l\}$ of $G_{r-1}$. The edges in the subdivision of $I \times \{0\}$ are labeled according to the oriented images in $G$. For each edge $E_i$ of $H_r$, identify the two edges of $I \times \{0\}$ that are labeled $E_i$. The quotient of $I \times \{0\}$ by this identification rule is a graph $K$. The quotient of all of $I \times [0, 1]$ by the identification rule is a surface $S$ that deformation retracts to $K$. The edges of $K$ inherit a labeling that defines a map $h : K \to G$. For each edge $E_i \subset H_r$ there is a unique edge in $K$ labeled $E_i$; for notational simplicity, we refer to this edge as $E_i$. The other edges of $K$ form a subgraph $K_0$ consisting of edges that are mapped to $G_{r-1}$. Note that $K_0 \subset \partial S$. Let $\mathcal{V}_K$ be the vertex set of $K$.

The next lemma states that an extended fold $F : G \to G'$ can be lifted to a homotopy equivalence between the graphs $K$ and $K'$ associated to the indivisible Nielsen paths $\rho_r$ and $\rho'_r$ respectively.

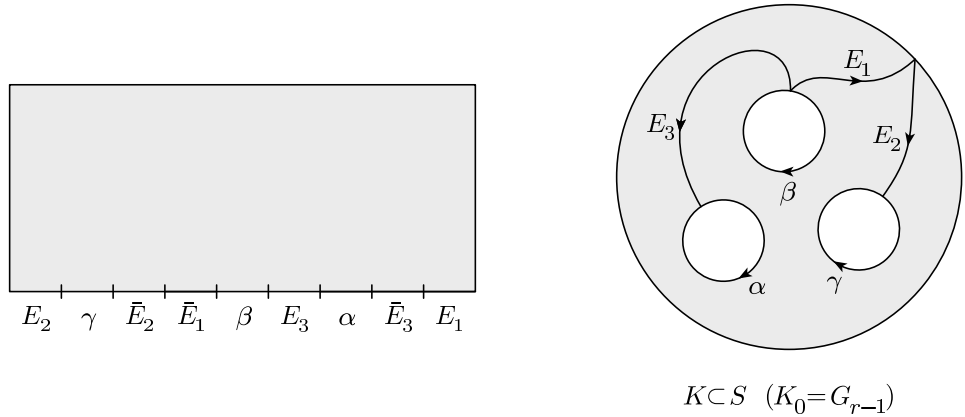

$K \subset S$   $(K_0 = G_{r-1})$

LEMMA 5.3.7.   *Suppose that $f : G \to G$, $f' : G' \to G'$, $\rho_r$ and $\rho'_r$ are as in Definition 5.3.2. Suppose further that $K$ and $K'$ are the graphs associated to $\rho_r$ and $\rho'_r$ respectively and that $h : K \to G$ and $h' : K' \to G'$ are the associated labeling maps. Then there is a homotopy equivalence $F_K : K \to K'$ such that:*

- *$h' F_K = Fh$.*

- *$F_K$ induces a bijection between the vertices of $K$ and the vertices of $K'$.*

- *$F_K | K_0 : K_0 \to K'_0$ is a homeomorphism.*

*Proof.* We use the notation of Definition 5.3.2. As described in the proof of Lemma 5.3.3, there is a one-to-one correspondence between the edges of $K$ and the edges of $K'$; this correspondence preserves labels with the single exception that $E_2$ corresponds to $E'_2$.

To study the link structures of $K$ and $K'$, we think of $\rho_r = e_1 \cdot \ldots \cdot e_m$ as a concatenation of subpaths where each $e_j$ is labeled by either an edge $E_i \subset H_r$ or a maximal subpath $b_l \subset H_r$. Thus each $e_j$ is identified with an edge of $K$. The link of a vertex in $K$ is the equivalence class of oriented edges in $K$ generated by $\bar{e}_j \sim e_{j+1}$ for $1 \leq j \leq m-1$. Let $\sim'$ be the analogous equivalence relation on oriented edges in $K'$.

If $b$ is trivial, let $A = E_1$ (thought of as an edge of $K$); if $b$ is not trivial, let $A = b_K$, the edge of $K$ that corresponds to $b$. Let $E_3$ be the edge that follows $b$ in $\bar{\alpha}$. The generating relations for $\sim'$ are obtained from the generating relations for $\sim$ as follows. Erase $E_1 \sim E_2$ and add $E_3 \sim' E'_2$. If $A = b_K$ or if the second occurrence of $E_1$ in $\rho_r$ is not followed by the second occurrence of $E_3$ in $\rho_r$, then erase $\bar{A} \sim E_3$. (This accounts for the changes caused by the shortening of $\alpha$ and $\beta$ to $\alpha'$ and $\beta'$.) If $\rho_r$ does not begin with $E_2$ or end with $\bar{E}_2$, then $\rho_r$ passes through a turn $(X, E_2)$ with $X \neq E_1$. Erase $X \sim E_2$, add $X \sim' E_1$ and add $\bar{A} \sim' E'_2$. (This accounts for the subdivision of $E_2$ into the edge path $AE'_2$ or $E_1 AE'_2$.) Finally, change $\sim$ to $\sim'$ in the remaining generating relations.

The combinatorial types of $K$ and $K'$ differ in at most two links. The initial endpoint of $E_2$ is moved from the vertex containing $X$ and $E_1$ and is added to the vertex containing $E_3$ and $\bar{A}$. (These vertices need not be distinct.)

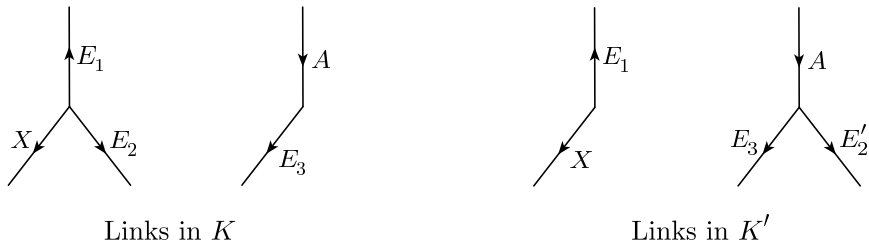

Links in $K$                                    Links in $K'$

We may think of $K \setminus E_2$ as a subgraph of both $K$ and $K'$; thus $K$ is formed by adding $E_2$ and $K'$ is formed by adding $E_2'$. The maps $h$ and $h'$ agree on $K \setminus E_2$; $h(E_2) = E_2$ and $h'(E_2') = E_2'$. Define $F_K$ to be the identity on $K \setminus E_2$ and $F_K(E_2) = E_1 b_K E_2'$ where we allow the possibility that $b_K$ is trivial. The desired properties of $F_K$ follow immediately. $\qquad\square$

The next step in the proof of Proposition 5.3.1 is to show that $f : G \to G$ can be lifted to $K$.

COROLLARY 5.3.8. *There is a homotopy equivalence $f_K : K \to K$ such that*:

(i) $hf_K \simeq fh$ rel $\mathcal{V}_K$.

(ii) $f_K$ *permutes the elements of* $\mathcal{V}_K$.

(iii) $f_K | K_0$ *is a homeomorphism.*

*Proof.* Let $\theta, f_{r-1}$ and $f_r$ be as in Lemma 5.3.6. There is a relative train track map $f^1 : G^1 \to G^1$ that is obtained by iteratively folding indivisible Nielsen paths, starting with $\rho_r$ and continuing through the extended folds that make up $f_r$. Let $h^1 : K^1 \to G^1$ be the associated labeling map. Lemma 5.3.7 implies that there exists $F_1 : K \to K^1$ such that $h^1 F_1 = f_r h$.

Let $\rho_r^1 = (f_r)_\#(\rho_r) \subset G^1$ be the corresponding indivisible Nielsen path and let $\{b_l^1\}$ be the maximal subpaths of $\rho_r^1$ in $G_{r-1}^1$. Since $f^1$ is a relative train track map and $f^1 | G_{r-1} = f | G_{r-1} = \theta f_{r-1} | G_{r-1}$, each $\theta_\#(f_{r-1})_\#(b_l^1)$ is a nontrivial path in $G_{r-1}$. By construction $\theta_\# f_{r-1} | H_r^1 : H_r^1 \to H_r^1$ is injective and $\theta_\#(f_{r-1})_\#(\rho_r^1) = f_\#(\rho_r) = \rho_r$. It follows that $\rho_r$ is obtained from $\rho_r^1$ by replacing each $b_l^1$ with $\theta_\#(f_{r-1})_\#(b_l^1)$ and relabeling an edge $E^1$ of $H_r^1$ by the edge $\theta_\#(f_{r-1})_\#(E^1)$; thus $K$ is obtained from $K^1$ by changing the edge labels $b_l^1$ to $\theta_\#(f_{r-1})_\#(b_l^1)$ and the edge labels $E^1$ to $\theta_\#(f_{r-1})_\#(E^1)$. This induces a homeomorphism $F_2 : K^1 \to K$ such that $hF_2 \simeq \theta f_{r-1} h^1$ rel vertices.

Define $f_K = F_2 F_1$. Then $h f_K = h F_2 F_1 \simeq \theta f_{r-1} h^1 F_1 = \theta f_{r-1} f_r h = f h$ rel $\mathcal{V}_K$. Since $F_1$ induces a bijection of vertices and $F_2$ is a homeomorphism, $f_K$ permutes the vertices of $K$. Finally, since $F_1 | K_0$ and $F_2$ are homeomorphisms, $f_K | K_0$ is a homeomorphism.                                                                 □

The next lemma exploits the fact that $f : G \to G$ is reduced. If $\hat{v}$ is vertex of $K$, we denote its *link*, thought of as the oriented edges with initial vertex $\hat{v}$, by $\mathrm{Lk}(K, \hat{v})$. The link of $v$ in $G$ is denoted $\mathrm{Lk}(G, v)$.

LEMMA 5.3.9.    1. *If* $h(\mathrm{Lk}(K, \hat{v})) \subset H_r$, *then* $h(\mathrm{Lk}(K, \hat{v})) = \mathrm{Lk}(G, h(\hat{v}))$.

2. *The endpoints of* $\rho_r$ *are equal and are not contained in* $G_{r-1}$; *the initial edge of* $\rho_r$ *is distinct from the initial edge of* $\bar{\rho}_r$.

3. *For each component* $C$ *of* $K_0$, $h_\#(C)$ *is a nontrivial circuit in* $G_{r-1}$.

*Proof.* Since $f_K | K_0$ is a homeomorphism and $f_K$ permutes the vertices of $K$, we may assume, after replacing $f$ by an iterate if necessary, that $f_K$ fixes every vertex of $K$ and restricts to the identity on $K_0$.

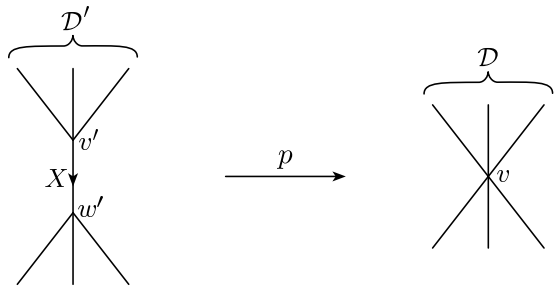

Suppose that $\mathcal{D} = h(\mathrm{Lk}(K, \hat{v})) \subset H_r$ is not all of $\mathrm{Lk}(G, h(\hat{v}))$. Define a new graph $G'$ from $G$ as follows. Replace $v = h(\hat{v})$ by a pair of vertices $v'$ and $w'$; reattach the edges of $\mathcal{D}$ to $v'$ ( where they are labeled $\mathcal{D}'$) and the remaining edges in $\mathrm{Lk}(G, v)$ to $w'$; add an edge $X$ connecting $v'$ and $w'$. Define $p : G' \to G_r$ to be the homotopy equivalence that collapses $X$ to $v$. Let $G' = G \setminus X$ and $G'_{r-1} = p^{-1}(G_{r-1}) \setminus X$. Note that $p | G'_{r-1} : G'_{r-1} \to G_{r-1}$ is a homeomorphism so that we may think of $G_{r-1}$ as a subgraph of both $G$ and $G'$. With this notation $p | G_{r-1} = $ identity. By construction, there is a map $h^* : K \to G'_r$ satisfying $p h^* = h$. If $E$ is an edge of $H_r$, and so also an edge of $K$, then for all $k > 0$, $h^*_\# f_K^k(E) \subset G'_r$ is a lift of $f_\#^k(E)$. Thus $G'_r$ carries the lamination $\Lambda^+ \in \mathcal{L}(\mathcal{O})$ associated to $H_r$ and $\mathcal{F}(G'_r) \neq \mathcal{F}(G_{r-1})$.

We next show that $\mathcal{F}(G'_r)$ is $\mathcal{O}$-invariant by defining a topological representative $f' : G' \to G'$ such that $p f' \simeq f p$ rel vertices and such that $f'$ restricts to a self map of $G'_r$. On the edges of $G'_{r-1} = G_{r-1}$, $f' = f$. If $E'$ is an edge of $H'_r = G'_r \setminus G'_{r-1}$, then $p(E')$ is an edge of $H_r$ and so $h^{-1}p(E')$ is a well-defined edge of $K$; define $f'(E') = h^* f_K h^{-1} p(E')$. Finally define $f'(X) = X$.

Lemma 5.3.8 implies that $pf' \simeq fp$ rel vertices on the edges of $H'_r$ and hence that $pf' \simeq fp$ rel vertices on all edges. Since $f_K$ fixes all vertices, so does $f'$. It follows that $f'$ is continuous. By definition, $f'$ restricts to a self map of $G'_r$.

To complete the proof of (1) we will show that $\mathcal{F}(G'_r) \neq \mathcal{F}(G') = \mathcal{F}(G)$ in contradiction to the assumption that $f : G \to G$ is reduced and our previous observation that $\mathcal{F}(G'_r) \neq \mathcal{F}(G_{r-1})$. If $\mathcal{F}(G'_r) = \mathcal{F}(G')$, then $G'_r$ must have two components, one of which, $Y'$, is contractible. Since $G'_r$ carries $\Lambda^+$, $Y' \subset G'_{r-1}$. But then $Y'$ is a contractible component of $G_{r-1}$ that contains the vertex $v = h(\hat{v})$ and is hence mapped to itself by $f$. This contradicts condition (5) in the definition of $\mathcal{F}$-Nielsen minimizing.

We now turn to the proof of (2). Let $v_1$ be the initial vertex of $\rho_r = \alpha\beta$. If $v_1 \notin G_{r-1}$, then each time the interior of $\rho_r$ passes through $v_1$ it crosses two edges of $H_r$. The total number of times that $\rho_r$ crosses the $H_r$-edges in $\mathrm{Lk}(G, v)$ is even. Since $\rho_r$ starts at $v_1$, it must also end at $v_1$. Suppose that the initial edge $E$ of $\rho_r$ equals the initial edge of $\bar{\rho}_r$. For large $k$, there are initial subpaths $\alpha_0$ and $\beta_0$ of $E$ such that $\alpha = f^k_{\#}(\alpha_0)$ and $\bar{\beta} = f^k_{\#}(\beta_0)$. We may assume without loss that $\alpha_0 \subset \beta_0$ and hence that $\alpha$ is an initial subpath of $\bar{\beta}$. Since the initial edge of $\beta$ is in $H_r$, the difference between the number of edges in $f^m_{\#}(\alpha)$ and $f^m_{\#}(\beta)$ grows without bound, in contradiction to the fact that for all $m \geq 0$, $\rho_r = f^m_{\#}(\rho_r)$ is obtained from the concatenation of $f^m_{\#}(\alpha)$ and $f^m_{\#}(\beta)$ by cancelling at the juncture.

To complete the proof of (2), we assume that $v_1 \in G_{r-1}$ and argue to a contradiction. Let $\hat{v}_1$ be the initial endpoint of the lift to $K$ of the first edge in $\rho_r$. Part (1) implies that $\hat{v}_1$ is contained in a component $C$ of $K_0$. Since $C$ is completed to a circle in $\partial S$ by $\partial I \times [0, 1] \cup I \times \{1\}$, $C$ is an arc and $\hat{v}_1$ is one of its endpoints. For notational concreteness, we give the argument when $C$ has two edges: $A$ with endpoints $\hat{v}_1$ and $\hat{v}_2$ and label $b_1$; and $B$ with endpoints $\hat{v}_2$ and $\hat{v}_3$ and label $b_2$. We also assume that $v_1 = v_3 \neq v_2$ where $v_i = h(\hat{v}_i)$. The arguments in other cases require only straightforward modifications.

The subset of $\mathrm{Lk}(K, \hat{v}_i)$ that projects to $H_r$ is denoted $\hat{\mathcal{D}}_i$; its image in $\mathrm{Lk}(G, v_i)$ is denoted $\mathcal{D}_i$. Note that $\mathrm{Lk}(K, \hat{v}_1) = \hat{\mathcal{D}}_1 \cup A$, $\mathrm{Lk}(K, \hat{v}_2) = \hat{\mathcal{D}}_2 \cup \bar{A} \cup B$ and $\mathrm{Lk}(K, \hat{v}_3) = \hat{\mathcal{D}}_3 \cup \bar{B}$. Define $G'$ to be the graph obtained from $G_r$ as follows. Replace $v_1 = v_3$ by three vertices $v'_1, v'_3$ and $w'_1$; replace $v_2$ by two vertices $v'_2$ and $w'_2$. The edges of $\mathcal{D}_i$ are reattached to $v'_i$ (where they are labeled $\mathcal{D}'_i$). The remaining edges in the link of $v_1$ (respectively $v_2$) are reattached to $w'_1$ (respectively $w'_2$). Add edges $X$ connecting $v'_1$ to $w'_1$, $A'$ connecting $v'_1$ to $v'_2$ and $B'$ connecting $v'_2$ to $v'_3$. Define $p : G' \to G$ to be the homotopy equivalence that collapses $X$ to $v_1$, sends $A'$ to $b_1$ and sends $B'$ to $b_2$. As in the previous case, there is a map $h^* : K \to G' \setminus X$ satisfying $ph^* = h$ and there is an induced topological representative $f' : G' \to G'$ : that satisfies $pf' = fp$; that restricts to the identity on $X \cup A' \cup B'$; and that maps $G' \setminus X$ into itself. Define $G'_{r-1} = p^{-1}(G_{r-1}) \setminus (X \cup A' \cup B')$ and $G'_r = G' \setminus X$. The proof now concludes as in the previous case.

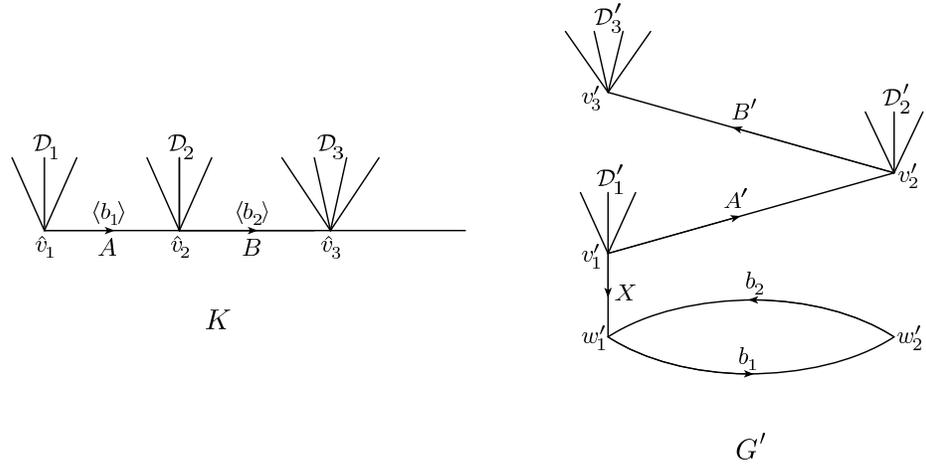

For part (3), we may assume by parts (1) and (2) that $C$ is a loop. We will give the argument when the boundary component $C$ consists of three edges: $A$ with endpoints $\hat{v}_1$ and $\hat{v}_2$ and label $b_1$; $B$ with endpoints $\hat{v}_2$ and $\hat{v}_3$ and label $b_2$; and $C$ with endpoints $\hat{v}_3$ and $\hat{v}_1$ and label $b_3$. As in part (2) we assume that $v_1 = v_3 \neq v_2$. Define $G'$ and $p : G' \to G$ exactly as in part (2). If $h_\#(C)$ is trivial, then $b_3 \simeq \bar{b}_2\bar{b}_1$ rel endpoints and there exists $h^* : K \to G'$ such that $ph^* \simeq h$ rel vertices; one simply defines $h^*(A) = A'$, $h^*(B) = B'$ and $h^*(C) = \bar{B}'\bar{A}'$. The proof now concludes as in the previous cases. □

*Proof of Proposition* 5.3.1. We have already defined $S$. Parts (1) and (2) of Lemma 5.3.9 imply that $\partial I \times \{0,1\} \cup I \times \{1\}$ projects to a component $\rho^*$ of $\partial S$. The components of $\partial S \setminus \rho^* = K_0$ are denoted $\gamma_1^*, \ldots, \gamma_m^*$. Part (3) of Lemma 5.3.9 implies that each $\gamma_i = h_\#(\gamma_i^*)$ is a nontrivial circuit. Since $hf_K \simeq fh$ and $f_K$ permutes the $\gamma_i^*$'s, $f_\#$ permutes the $\gamma_i$'s and the induced permutations of $\{1, \ldots, m\}$ agree.

Define $\hat{G} = G_{r-1} \cup_{h|K_0} K$. The identity map on $G_{r-1}$ and $h$ fit together to give a continuous map $\hat{h} : \hat{G} \to G$. Part (1) of Lemma 5.3.9 implies that $\hat{h}$ induces a bijection between the vertices of $K \setminus K_0$ and the vertices of $G \setminus G_{r-1}$. Thus $\hat{h}$ induces a bijection between the vertices of $\hat{G}$ and the vertices of $G$. Since $\hat{h}$ also induces a bijection on edges, $\hat{h}$ is a homeomorphism.

Let $\mathcal{A}$ be the union of $m$ annuli $A_i, \ldots, A_m$. Define $Y$ to be the space obtained from $G_{r-1} \cup \mathcal{A} \cup S$ by attaching one end of $A_i$ to $\gamma_i$ and the other end to $\gamma_i^*$. By construction, $S = M(\rho_r^K)$ is the mapping cylinder of the quotient map $\rho_r^K : I \times \{0\} \to K$. Since $\hat{h} : \hat{G} \to G$ is a homeomorphism, $Y$ is homeomorphic to the mapping cylinder $M(\rho_r)$, where $\rho_r$ is thought of as a map of the interval into $G$. The natural deformation retraction of $M(\rho_r)$ to $G$ defines $\Phi : (Y, G_{r-1}) \to (G, G_{r-1})$. It remains to show that if $r : S \to K$ is the deformation retraction given by collapsing mapping cylinder lines, then the homotopy equivalence $f_K r : S \to S$ is homotopic to a pseudo-Anosov homeomorphism $\phi$.

Since $r|K_0 = $ identity and $f_K|K_0$ is a homeomorphism, $f_K r$ permutes the components of $K_0 \subset \partial S$. The component $\rho^*$ is freely homotopic in $S$ to the circuit determined by $\rho_r$ and so is also fixed by $(f^K r)_\#$. It follows (Theorem 3.1 of [Hem]) that $f^K r$ is homotopic to a homeomorphism $\phi$. To prove that the mapping class determined by $\phi$ is pseudo-Anosov, it suffices to show that the only periodic conjugacy classes are the peripheral ones.

A nontrivial circuit $\sigma_K \subset K$ determines a nontrivial circuit $\hat{\sigma} \subset \hat{G}$ and hence a nontrivial circuit $\sigma \subset G$. If $\sigma_K$ is periodic under the action of $f_K$ and $\sigma_K \not\subset K_0$, then $\sigma$ is periodic under the action of $f$ and $\sigma \not\subset G_{r-1}$. Lemma 4.2.6 and Lemma 5.2.5 imply that $\sigma$ splits into subpaths that are either entirely contained in $G_{r-1}$ or equal to $\rho_r$ or $\bar{\rho}_r$. Part (2) of Lemma 5.3.9 implies that $\sigma$ is a multiple of $\rho_r$ or $\bar{\rho}_r$.                                             □

5.4. *Sliding.* In this subsection we introduce and study the technique used to arrange condition (ne-ii) of Theorem 7.0.1. We assume throughout this subsection that $f : G \to G$ is a relative train track map, that $H_i$ is a non-exponentially-growing stratum and that each non-exponentially-growing stratum $H_j$ is a single edge $E_i$ satisfying $f(E_j) = E_j u_j$ for some path $u_j \subset G_{j-1}$.

For any path $\alpha \subset G_{i-1}$ with initial endpoint equal to the terminal endpoint of $E_i$ and with terminal endpoint at a vertex of $G_{i-1}$, define a new graph $G'$ by replacing $E_i$ with an edge $E'_i$ that has the same initial endpoint as $E_i$ and the same terminal endpoint as $\alpha$. Every edge of $G \setminus E_i$ is naturally identified with an edge of $G' \setminus E'_i$; we use the same name for the edge in both graphs. Similarly, a path $\beta \subset G$ that does not cross $E_i$ is identified with a path, also called $\beta$, in $G'$.

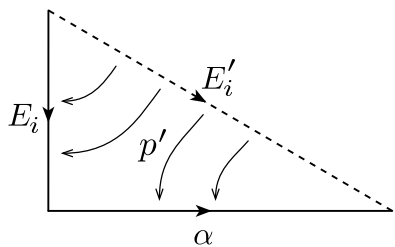

There are homotopy equivalences $p : G \to G'$ and $p' : G' \to G$ that equal the 'identity' on the common edges of $G$ and $G'$ and that satisfy $p(E_i) = E'_i \bar{\alpha}$ and $p'(E_i') = E_i \alpha$ respectively. Define $f' : G' \to G'$ by tightening $pfp' : G' \to G'$; in other words, on each edge of $G'$, $f' = (pfp')_\#$. We say that $f' : G' \to G'$ is obtained from $f : G \to G$ by *sliding* $E_i$ *along* $\alpha$. For each $G_j$ define $G_j' = p(G_j)$.

The basic properties of sliding are listed in the following lemma.

LEMMA 5.4.1. *Suppose that* $f' : G' \to G'$ *is obtained from* $f : G \to G$ *by sliding* $E_i$ *along* $\alpha$. *Then* $f'(E_i') = E_i' u_i'$ *where* $u_i' = [\bar{\alpha} u f(\alpha)] \subset G_{i-1}'$. *Moreover,* $H_j'$ *is exponentially growing (respectively non-exponentially growing) if and only if* $H_j$ *is exponentially growing (respectively non-exponentially growing). If* $f : G \to G$ *is* $\mathcal{F}$-*Nielsen minimizing, then so is* $f' : G' \to G'$ .

*Proof.* We have $f'(E_i') = (pfp')_\#(E_i') = (pf)_\#(E_i \alpha) = p_\#([E_i u f(\alpha)]) = [E_i' \bar{\alpha} u f(\alpha)] = E_i'[\bar{\alpha} u f(\alpha)]$. Sliding has no effect on any $\mathcal{F}(G_j)$ so $f' : G' \to G'$ satisfies conditions (1) and (2) in the definition of $\mathcal{F}$-Nielsen minimizing if $f : G \to G$ does.

Since $f|G_{i-1}$ agrees with $f'|G_{i-1}'$, we may restrict our attention to strata $H_j$ with $j > i$. If $H_j$ is a zero stratum and $E_0$ is an edge in $H_j$, then $f(E_0) \subset G_{j-1}$. Thus $f'(E_0) = (pf)_\#(E_0) \subset G_{j-1}'$ and $H_j'$ is also a zero stratum. If $H_j$ is non-exponentially growing, then, according to our standing assumption, $H_j$ is a single edge $E_j$ and there is a path $u_j \subset G_{j-1}$ such that $f(E_j) = E_j u_j$. Thus $f'(E_j) = (pf)_\#(E_j) = p_\#(E_j u_j) = E_j p_\#(u_j)$ where $p_\#(u_j) \subset G_{j-1}'$. Thus $H_j'$ is non-exponentially growing.

Suppose now that $H_r$ is exponentially growing and that $E$ is an edge of $H_r$. For any nontrivial paths $\beta \subset G$ and $\gamma' \subset G'$ with endpoints at vertices, $(p'p)_\#(\beta) = \beta$ and $(pp')_\#(\gamma') = \gamma'$. In particular, $p_\#(\beta)$ and $p_\#'(\gamma')$ are nontrivial. The train track property implies that $f(E) = a_1 b_1 a_2 \ldots b_l a_{l+1}$ where $a_i \subset H_r$ and $b_i \subset G_{r-1}$ are nontrivial paths. Thus $f'(E) = (pf)_\#(E) = p_\#(a_1 b_1 a_2 \ldots b_l a_{l+1}) = a_1 p_\#(b_1) a_2 \ldots p_\#(b_l) a_{l+1}$. This implies that $H_r$ is exponentially growing and that the transition submatrices, $M_r$ and $M_r'$, and hence the Perron-Frobenius eigenvalues $\mu_r$ and $\mu_r'$, are equal. If $\beta' \subset G_{r-1}'$ is a nontrivial path with endpoints in $H_r' \cap G_{r-1}'$, then $p_\#'(\beta')$ is a nontrivial path with endpoints in $H_r \cap G_{r-1}$ and hence $(fp')_\#(\beta')$ is nontrivial. Thus $f_\#'(\beta') = (pfp')_\#(\beta')$ is nontrivial. We have now verified that $f' : G' \to G'$ is a relative train track map.

If $\sigma' \subset G_r'$ satisfies $(f')_\#^k(\sigma') = \sigma'$, then $p_\#'(\sigma') = p_\#'(pfp')_\#^k(\sigma') = (f)_\#^k p_\#'(\sigma')$. Similarly, if $\sigma \subset G_r$ and $f_\#^k(\sigma) = \sigma$, then $p_\#(\sigma) = (pf^k)_\#(\sigma) = (f')_\#^k p_\#(\sigma)$. In other words, $p$ induces a period-preserving bijection between indivisible periodic Nielsen paths in $G_r$ and indivisible periodic Nielsen paths in $G_r'$. Thus $f' : G' \to G'$ is $\mathcal{F}$-Nielsen minimized if $f : G \to G$ is. $\square$

In order to find good paths along which to slide, we consider a restricted lift of $f$ defined as follows. Choose a lift $\tilde{E}_i$ in the universal cover $\Gamma$ of $G$. Let $\tilde{f} : \Gamma \to \Gamma$ be the lift of $f : G \to G$ that fixes the initial endpoint of $\tilde{E}_i$, let $\tilde{u}_i$ be the lift of $u_i$ satisfying $\tilde{f}(\tilde{E}_i) = \tilde{E}_i \tilde{u}_i$ and let $\Gamma_{i-1} \subset \Gamma$ be the component of the full pre-image of $G_{i-1}$ that contains $\tilde{u}_i$. Since $\tilde{f}(\tilde{E}_i) = \tilde{E}_i \tilde{u}_i$, $\Gamma_{i-1}$ is $\tilde{f}$-invariant. Denote $\tilde{f}|\Gamma_{i-1}$ by $h : \Gamma_{i-1} \to \Gamma_{i-1}$ and note that if $\tilde{q}$ is the initial vertex of $\tilde{u}_i$, then $h(\tilde{q})$ is its terminal vertex.

The sliding operation can be lifted to $\Gamma$ by replacing each lift of $E_i$ with a lift of $E'_i$. Call the resulting tree $\Gamma'$. Let $\tilde{f}' : \Gamma' \to \Gamma'$ be the lift of $f' : G' \to G'$ that fixes the initial endpoint of the lift $\tilde{E}'_i$ that corresponds to $\tilde{E}_i$. ($\tilde{E}_i$ and $\tilde{E}'_i$ have a 'common' initial endpoint.) Lemma 5.4.1 implies that $\Gamma'_{i-1} = \Gamma_{i-1}$ and that $\tilde{f}'|\Gamma'_{i-1} = \tilde{f}|\Gamma_{i-1}$. In this sense, $h : \Gamma_{i-1} \to \Gamma_{i-1}$ is unchanged by the sliding operation.

For any $\tilde{x}, \tilde{y} \in \Gamma_{i-1}$, denote the path connecting $\tilde{x}$ to $\tilde{y}$ by $[\tilde{x}, \tilde{y}]$ and its image under the covering projection $\mathrm{pr} : \Gamma \to G$ by $\mathrm{pr}([\tilde{x}, \tilde{y}]) \subset G_{i-1}$. Paths in $G_{i-1}$ that have their initial endpoint at $q = \mathrm{pr}(\tilde{q})$ and their terminal endpoint at a vertex $v$ are in one-to-one correspondence with paths in $\Gamma_{i-1}$ of the form $[\tilde{q}, \tilde{v}]$ where $\mathrm{pr}(\tilde{v}) = v$ and hence are in one-to-one correspondence with the set of vertices $\tilde{v}$ of $\Gamma_{i-1}$. Thus we may speak of *sliding along the path determined by the vertex $\tilde{v}$*.

LEMMA 5.4.2. *If $f' : G' \to G'$ is obtained from $f : G \to G$ by sliding along the path corresponding to a vertex $\tilde{v}$, then $f'(E'_i) = E'_i u'_i$ where $u'_i = \mathrm{pr}([\tilde{v}, h(\tilde{v})])$.*

*Proof.* This follows immediately from Lemma 5.4.1 and the definition of $h$. □

The following proposition is the main result of this subsection. Conditions (1) and (2) record the fact that we have taken $u$ to be as simple as $h$ will allow. Condition (3) is a strengthening of the assertion that $f'(E'_i) = E'_i \cdot u'_i$ and is used in Lemma 5.5.1.

PROPOSITION 5.4.3. *After subdivision at a periodic orbit of $f : G \to G$ if necessary, there is a vertex $\tilde{v}$ of $\Gamma$ that projects to a periodic point of $f$ so that if $f' : G' \to G'$ is obtained from $f : G \to G$ by sliding along the path determined by $\tilde{v}$, then $f'(E'_i) = E'_i \cdot u'_i$ where:*

1. *$u'_i$ is trivial if and only if $h$ has a fixed point.*

2. *If $u'_i$ is nontrivial, then $u'_i$ is periodic under the action of $f_\#$ if and only if $h$ commutes with a covering translation $T$ of $\Gamma_{i-1}$; in this case, the infinite ray $\tilde{R}' = \tilde{u}'_i h_\#(\tilde{u}'_i) h_\#^2(\tilde{u}'_i) \ldots$ is contained in the axis of $T$.*

3. *If $u'_i$ is nontrivial, then $E'_i \cdot w'_i$ is a splitting for every initial segment $w'_i$ of $u'_i$.*

*Proof.* If Fix($h$) $\neq \emptyset$, then, after subdividing if necessary, we may choose $\tilde{v} \in$ Fix($h$). Lemma 5.4.2 implies that $u'$ is the trivial path. We assume now that $h$ is fixed-point free. Lemma 5.4.2 implies that $u'$ cannot be trivial for any choice of $\tilde{v}$. This verifies condition (1).

Let $\tilde{X} = \{\tilde{x} : \{\tilde{x}, h(\tilde{x}), h^2(\tilde{x}), \ldots\}$ is an ordered subset of an infinite path in $\Gamma_{i-1}\}$. Note that $h(\tilde{X}) \subset \tilde{X}$ and that if $\tilde{x} \in \tilde{X}$ and $[\tilde{x}, h(\tilde{x})] = [\tilde{x}, \tilde{y}] \cdot [\tilde{y}, h(\tilde{x})]$, then $\tilde{y} \in \tilde{X}$. The first step in the proof is to show that $\tilde{X} \neq \emptyset$.

We say that the initial edge of $[\tilde{v}, h(\tilde{v})]$ is *preferred* by the vertex $\tilde{v} \in \Gamma_{i-1}$. If $E$ is preferred by both (respectively neither) of its endpoints, then $h(E)$, thought of as an edge path, contains $\bar{E}$ (respectively $E$). But then some subinterval of $E$ maps to all of $\bar{E}$ (respectively $E$) and so must contain a fixed point. This contradiction implies that $E$ is preferred by exactly one of its endpoints. For each $\tilde{v} \in \Gamma_{i-1}$, let $L(\tilde{v}) = E^0 \cdot E^1 \cdots$ be the infinite path defined by choosing $E^0$ to be the preferred edge for $\tilde{v}$, and by inductively choosing $E^{i+1}$ to be the preferred edge for the terminal endpoint of $E^i$.

Given $\tilde{v}, \tilde{w} \in \Gamma_{i-1}$, denote $[\tilde{v}, \tilde{w}]$ by $\tilde{\gamma}$. An easy induction argument on edge length shows that either the initial edge of $\tilde{\gamma}$ is preferred by the initial vertex of $\tilde{\gamma}$ or the terminal edge of $\tilde{\gamma}$ is preferred by the terminal vertex of $\tilde{\gamma}$ or both . In other words, either the initial edge of $\tilde{\gamma}$ is the initial edge of $L(\tilde{v})$ or the initial edge of the inverse of $\tilde{\gamma}$ is the initial edge of $L(\tilde{w})$ or both. It follows that $L(\tilde{v})$ and $L(\tilde{w})$ have a common infinite end. Define $I(\tilde{v}, \tilde{w})$ to be the initial vertex of $L(\tilde{v}) \cap L(\tilde{w})$; thus $L(I(\tilde{v}, \tilde{w})) = L(\tilde{v}) \cap L(\tilde{w})$.

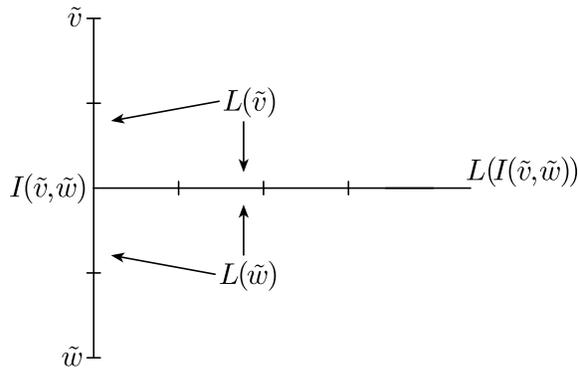

Choose a vertex $\tilde{v}_0$ and inductively define $\tilde{v}_{i+1} = I(\tilde{v}_i, h(\tilde{v}_i))$. Then $L(\tilde{v}_0) = [\tilde{v}_0, \tilde{v}_1][\tilde{v}_1, \tilde{v}_2][\tilde{v}_2, \tilde{v}_3] \ldots$ and $h_\#([\tilde{v}_i, \tilde{v}_{i+1}]) = [h(\tilde{v}_i), h(\tilde{v}_{i+1})]$ contains $[\tilde{v}_{i+1}, \tilde{v}_{i+2}]$. Define $\tilde{Y}_m = \{\tilde{y} \in [\tilde{v}_0, \tilde{v}_1] : h^i(\tilde{y}) \in [\tilde{v}_i, \tilde{v}_{i+1}] \ \forall \ 0 \leq i \leq m\}$. A straightforward induction argument shows that $h(\tilde{Y}_m) = [\tilde{v}_m, \tilde{v}_{m+1}]$ and in particular that $\tilde{Y}_m \neq \emptyset$. The $\tilde{Y}_m$'s are a nested sequence of closed subsets of $[\tilde{v}_0, \tilde{v}_1]$ so that $\cap \tilde{Y}_m \neq \emptyset$. By construction, $\cap \tilde{Y}_m \subset \tilde{X}$ so that $\tilde{X} \neq \emptyset$.

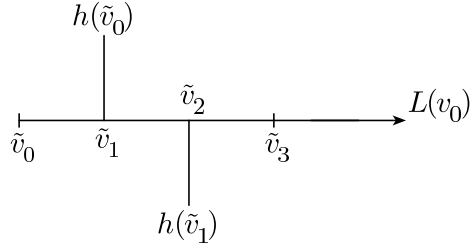

The next step is to choose a vertex $\tilde{v} \in \tilde{X}$. Let $s$ be the smallest positive integer for which there exists $\tilde{x} \in \tilde{X}$ satisfying $\mathrm{pr}[\tilde{x}, h(\tilde{x})] \subset G_s$. Choose such an $\tilde{x}$. If $H_s$ is non-exponentially growing, then $H_s$ is a single edge $E_s$ and $f(E_s) = E_s u_s$ for some path $u_s \subset G_{s-1}$. After replacing $\tilde{x}$ by some $h^k(\tilde{x})$ if necessary, we may assume that $[\tilde{x}, h(\tilde{x})]$ contains at least one entire edge $\tilde{e}$ whose projected image $e$ equals either $E_s$ or $\bar{E}_s$ and that $\mathrm{pr}(h(\tilde{x}))$ is not contained in the interior of $E_s$. If $e = E_s$ let $\tilde{v}$ be the initial edge of $\tilde{e}$; if $e = \bar{E}_s$, let $\tilde{v}$ be the terminal edge of $\tilde{e}$. Lemma 4.1.4 implies that $[\tilde{x}, h(\tilde{x})]$ can be split at $\tilde{v}$ and hence that $\tilde{v} \in \tilde{X}$.

If $H_s$ is exponentially growing, then after replacing $\tilde{x}$ by some $h^k(\tilde{x})$ if necessary, we may assume that each $[h^i(\tilde{x}), h^{i+1}(\tilde{x})]$ has the same number of illegal turns in $H_s$. Lemma 4.2.6 produces a splitting of $[\tilde{x}, h(\tilde{x})]$. If one of the resulting pieces is a lift $\tilde{\rho}$ of some $\rho \in P_s$, let $\tilde{v}$ be its initial endpoint. Replacing $\tilde{x}$ by some $h^k(\tilde{x})$ if necessary, we may assume by Lemma 4.2.5 that $v$ is $f$-periodic. After subdividing at the orbit of $v$, we may assume that $\tilde{v} \in \tilde{X}$ is a vertex. If there are no $\tilde{\rho}$ pieces, then $[\tilde{x}, h(\tilde{x})]$ is $s$-legal. After replacing $\tilde{x}$ by some $h^k(\tilde{x})$ if necessary, we may assume that $[\tilde{x}, h(\tilde{x})]$ contains an entire edge of $H_s$. Let $\tilde{v}$ be an endpoint of such an edge. Lemma 4.2.1 implies that $[\tilde{x}, h(\tilde{x})]$ splits at $\tilde{v}$ and hence that $\tilde{v} \in \tilde{X}$. Replacing $\tilde{v}$ by some $h^i(\tilde{v})$, we may assume that $v$ is a periodic point.

We assume from now on that $\tilde{u}_i' = [\tilde{v}, h(\tilde{v})]$, where $\tilde{v}$ is chosen as above, and that $u_i' \subset G_s$ is its projected image. If $(f')_\#^k(u_i') = u_i'$ for some $k > 0$, then the infinite ray $\tilde{R}' = \tilde{u}_i' h_\#(\tilde{u}_i') h_\#^2(\tilde{u}_i') \dots$ is contained in the axis of the covering translation $T : \Gamma_{i-1} \to \Gamma_{i-1}$ that satisfies $T(\tilde{u}_i') = h_\#^k(\tilde{u}_i')$. Since $h_\#$ preserves the axis of $T$, $h$ commutes with $T$.

Conversely, suppose that there is a covering translation $T$ of $\Gamma_{i-1}$ that commutes with $h$. Then $T([\tilde{x}, h(\tilde{x})]) = [T\tilde{x}, h(T\tilde{x})]$ for all $\tilde{x}$ and so $L(T(\tilde{v})) = T(L(\tilde{v}))$. This implies that $L(\tilde{v})$ and $T(L(\tilde{v}))$ have a common infinite end and hence that $L(\tilde{v})$ and the axis of $T$ have a common infinite end. In particular, $h^l(\tilde{v})$ is contained in the axis of $T$ for all sufficiently large $l$. This implies that there is a uniform bound to the edge length of $[h^l(\tilde{v}), h^{l+1}(\tilde{v})]$, and hence that $(f')_\#^l(u_i')$ takes only finitely many values. After replacing $\tilde{v}$ by some $h^l(\tilde{v})$ if necessary, we may assume that $u_i'$ is periodic under the action of $f'_\#$ and that $\tilde{R}'$ is contained in the axis of $T$. We have now verified (2).

Recall that $\tilde{R}' = \tilde{u}'_i \cdot h_\#(\tilde{u}'_i) \cdot h^2_\#(\tilde{u}'_i) \ldots$ is the infinite ray starting at $\tilde{v}$ and containing $\{h^i(\tilde{v}) : i \geq 0\}$. Let $\tilde{u}'_i = \tilde{\sigma}_1 \cdot \tilde{\sigma}_2 \cdot \ldots \cdot \tilde{\sigma}_{n_0}$ be the splitting provided by Lemma 4.2.1 if $H_s$ is non-exponentially growing and by Lemma 4.2.6 if $H_s$ is exponentially growing. Then $h^i_\#(\tilde{u}'_i) = h^i_\#(\tilde{\sigma}_1) \cdot h^i_\#(\tilde{\sigma}_2) \cdot \ldots \cdot h^i_\#(\sigma_{n_0})$ and $\tilde{R}'$ has an infinite splitting $\tilde{R}' = \tilde{\sigma}_1 \cdot \tilde{\sigma}_2 \cdot \ldots$ where $\tilde{\sigma}_{in_0+j} = h^i_\#(\tilde{\sigma}_j)$.

To verify condition (3), we must show that $[\tilde{v}, h^i(\tilde{v})]$ is contained in $[\tilde{v}, h^i(\tilde{y})]$ for all $\tilde{y} \in \tilde{u}'_i$. It suffices to show that $h^i(\tilde{u}'_i)$ intersects $h^{i-1}_\#(\tilde{u}'_i)$ trivially. We will prove the slightly stronger statement that $h^i(\tilde{\sigma}_j)$, which tightens to $\tilde{\sigma}_{in_0+j}$, intersects $\tilde{\sigma}_{in_0+j-1}$ trivially for $1 \leq j \leq n_0$ and all $i \geq 0$.

Suppose at first that $H_s$ is non-exponentially growing. If the initial edge of $\tilde{\sigma}_j$ is a lift of $E_s$, then Lemma 4.1.4 implies that $h^i(\tilde{\sigma}_j)$ is a lift of $E_s$, possibly followed by a sequence of edges in $\Gamma_{i-1}$ and possibly terminating in a lift of $\bar{E}_s$; if $h^i(\tilde{\sigma}_j)$ does terminate in a lift of $\bar{E}_s$ then the sequence of edges in $\Gamma_{i-1}$ tightens to a nontrivial path. The initial edge of $h^i(\tilde{\sigma}_j)$ is disjoint from the rest of $h^i(\tilde{\sigma}_j)$ and acts as a barrier to keep $h^i(\tilde{\sigma}_j)$ from intersecting $\tilde{\sigma}_{in_0+j-1}$. If the initial edge of $\tilde{\sigma}_j$ is not a lift of $E_s$, then the terminal end of $\tilde{\sigma}_{in_0+j-1}$ is a lift of $\bar{E}_s$ and $h^i(\tilde{\sigma}_j)$ is a sequence of edges in $\Gamma_{i-1}$, possibly terminating in a lift of $\bar{E}_s$; if the last edge of $h^i(\tilde{\sigma}_j)$ is a lift of $\bar{E}_s$, then the sequence of edges in $\Gamma_{i-1}$ tightens to a nontrivial path. The edges in $\Gamma_{i-1}$ cannot cross the terminal edge of $\tilde{\sigma}_{in_0+j-1}$ and the terminating lift of $\bar{E}_s$, if it exists, cannot either because the path between it and $\tilde{\sigma}_{in_0+j-1}$ is nontrivial.

Assume now that $H_s$ is exponentially growing. If $\tilde{\sigma}_j = \tilde{\rho}$ for some $\rho \in P_s$, let $\rho = \alpha\beta$ be the unique decomposition of $\rho$ into $s$-legal subpaths. Lemma 2.5.1 implies that $h^i(\tilde{\alpha})$ (and also $h^i(\tilde{\beta})$) decomposes into a path that projects to $H_s$, followed by a sequence of edges that project to $G_{s-1}$ and tighten to a nontrivial path, followed by a path that projects to $H_s$ and so on. The terminal end of $h^i(\tilde{\alpha})$ and the initial end of $h^i(\tilde{\beta})$ agree up to a point, but there is always an initial subpath of $h^i(\tilde{\alpha})$ that is disjoint from the rest of $h^i(\tilde{\sigma}_j)$. This prevents $h^i(\tilde{\sigma}_j)$ from crossing back to $\tilde{\sigma}_{in_0+j-1}$.

If $\tilde{\sigma}_j \neq \tilde{\rho}$, then $\tilde{\sigma}_j$ is $s$-legal. If the initial edge of $\tilde{\sigma}_j$ projects to $H_s$, then Lemma 2.5.1 implies that the initial edge of $h^i(\tilde{\sigma}_j)$ is disjoint from the rest of $h^i(\tilde{\sigma}_j)$ and prevents $h^i(\tilde{\sigma}_j)$ from crossing back to $\tilde{\sigma}_{in_0+j-1}$. If the initial edge of $\tilde{\sigma}_j$ does not project to $H_s$, then the terminal edge of $\sigma_{in_0+j-1}$ projects to $H_s$ and prevents $h^i(\tilde{\sigma}_j)$ from crossing back to $\tilde{\sigma}_{in_0+j-1}$: edges that project to $G_{s-1}$ cannot cross this barrier and edges that project to $H_s$ do not because Lemma 2.5.1 implies that they are part of $h^i_\#(\tilde{\sigma}_j)$. $\qquad\square$

5.5. *Splitting basic paths.* We proved in Lemma 4.1.4 that if $f : G \to G$ is a relative train track map and $H_i$ is a single edge $E_i$ satisfying $f(E_i) = E_i u_i$ for some path $u_i \subset G_{i-1}$, then a path $\sigma \subset G_i$ splits into subpaths that are either entirely contained in $G_{i-1}$ or are basic paths of height $i$; i.e. are one of the following three types: $E_i\gamma, E_i\gamma\bar{E}_i, \gamma\bar{E}_i$, where $\gamma \subset G_{i-1}$. The path $\gamma\bar{E}_i$ is the inverse of the path $E_i\bar{\gamma}$ so it suffices to consider $E_i\gamma$ and $E_i\gamma\bar{E}_i$. In this subsection we consider further splittings of the paths $E_i\gamma$ and $E_i\gamma\bar{E}_i$.

LEMMA 5.5.1. *Assume that*:

- *$f : G \to G$ is a relative train track map;*

- *each exponentially-growing stratum $H_i$ satisfies conditions eg-(i), eg-(ii), and eg-(iii) of Theorem 7.0.1;*

- *each non-exponentially-growing stratum $H_i$ is a single edge $E_i$ and satisfies the conclusions of Proposition 5.4.3;*

- *if $H_i$ is an exponentially-growing stratum then every indivisible periodic Nielsen path $\rho \subset G_i$ that intersects $H_i$ nontrivially has period one.*

*Then every periodic Nielsen path has period one. If $H_i$ is a non-exponentially-growing stratum and if $\gamma \subset G_{i-1}$ is a nontrivial path then the following are satisfied*:

(1) *If $E_i\gamma$ (respectively $E_i\gamma\bar{E}_i$) can be split at a point in the interior of $E_i$, then $f_\#^m(E_i\gamma) = E_i \cdot \gamma_1$ (respectively $f_\#^m(E_i\gamma\bar{E}_i) = E_i \cdot \gamma_1\bar{E}_i$) for some $m \geq 0$ and $\gamma_1 \subset G_{i-1}$.*

(2) *If $E_i\gamma$ has no splittings, then some $f_\#^m(E_i\gamma)$ is an exceptional path of height $i$ (Definition 5.1.3).*

(3) *If $E_i\gamma\bar{E}_i$ has no splittings, then $E_i\gamma\bar{E}_i$ is an exceptional path of height $i$.*

*Proof.* The proof is by induction. If $G_1 = H_1$ is exponentially growing, then conditions (1)–(3) for $f|G_1$ are vacuous and the condition on periodic Nielsen paths follows directly from our hypotheses. If $G_1 = H_1$ is non-exponentially growing, then $f$ pointwise fixes the single edge $E_1$ in $H_1$. Conditions (1)–(3) are therefore vacuous and the condition on periodic Nielsen paths follows from the fact that $f|G_1$ is the identity.

We assume now that the lemma holds for $f|G_{i-1}$ and prove it for $f|G_i$. If $H_i$ is exponentially growing, then the the condition on periodic Nielsen paths follows from Lemma 4.2.6 and the inductive hypothesis; conditions (1)–(3) are vacuous. We may therefore assume that $H_i$ is a single non-exponentially-growing edge $E_i$ and that $f(E_i) = E_i \cdot u_i$. Let $s$ be the smallest positive integer for which $u_i \subset G_s$.

Proposition 5.4.3 implies that for any nontrivial initial segment $\sigma_1$ of $E_i$, some $f_\#^m(\sigma_1) = E_i \cdot \gamma'$ where $\gamma' \subset G_{i-1}$. Thus if a path $\sigma$ splits as $\sigma = \sigma_1 \cdot \sigma_2$ where $\sigma_1$ is a nontrivial initial segment of $E_i$ then some $f_\#^m(\sigma)$ has a splitting of the form $E_i \cdot \sigma'$. Part (1) of the lemma now follows from the fact that each $f_\#^m(E_i\gamma)$ is of the form $E_i\gamma_1$ and each $f_\#^m(E_i\gamma\bar{E}_i)$ has the form $E_i\gamma_1\bar{E}_i$ where $\gamma_1 \subset G_{i-1}$.

In order to treat (2) and (3) simultaneously, let $\sigma = E_i\gamma$ or $E_i\gamma\bar{E}_i$. Assume that $\sigma$ has no splittings; in particular, $u_i$ is nontrivial.

*Step* 1 (cancelling large middle segments). As a first step in the proof of (2) and (3), we use the absence of splittings to show that if $\sigma = \sigma_1'\sigma_2'\sigma_3'$ is any decomposition into nontrivial subpaths, then there exist $M > 0$ and an $M$-splitting $\sigma = \sigma_1\sigma_2\sigma_3$ such that $\sigma_1$ is an initial subpath of $\sigma_1'$, $\bar{\sigma}_3$ is an initial subpath of $\bar{\sigma}_3'$ and $f_\#^M(\sigma) = f_\#^M(\sigma_1)f_\#^M(\sigma_3)$ where the indicated juncture point is a vertex. There is no loss in assuming that $\sigma_1'$ is contained in the initial edge of $\sigma$ and that $\sigma_3'$ is contained in the terminal (possibly partial) edge of $\sigma$.

It is convenient to work with lifts $\tilde{f}$ and $\tilde{\sigma} = \tilde{\sigma}_1'\tilde{\sigma}_2'\tilde{\sigma}_3'$. The set $\hat{S}_k = \{\tilde{x} \in \tilde{\sigma} : \tilde{f}^k(\tilde{x}) \in \tilde{f}_\#^k(\tilde{\sigma})\}$ is closed. Since $\tilde{\sigma}$ can be split at any point of $\cap_{k=1}^\infty \tilde{S}_k$ (Lemma 4.1.1(3)), this infinite intersection contains only the endpoints of $\tilde{\sigma}$. Thus there exists $M > 0$ so that $\cap_{k=1}^M \tilde{S}_k \subset \tilde{\sigma}_1' \cup \tilde{\sigma}_3'$. An easy induction argument shows that $f^N$ maps $\cap_{k=1}^N \tilde{S}_k$ onto $\tilde{f}_\#^N(\tilde{\sigma})$ for all $N \geq 1$. Since the $\tilde{E}_i$ that is the initial edge in $\tilde{f}(\tilde{\sigma})$ is not canceled when $\tilde{f}^k(\tilde{\sigma})$ is tightened to $\tilde{f}_\#^k(\tilde{\sigma})$, each $\tilde{S}_k$, and hence $\cap_{k=1}^M \tilde{S}_k$, contains an initial segment of $\tilde{E}_i$. Choose a point $\tilde{x} \in (\cap_{k=1}^M \tilde{S}_k) \cap \sigma_1'$ so that $\tilde{f}^M(\tilde{x})$ is as close to the terminal end of $\tilde{f}_\#^M(\sigma)$ as possible and let $\tilde{\sigma}_1$ be the initial segment of $\tilde{E}_i$ terminating at $\tilde{x}$. The choice of $\tilde{x}$ guarantees that $\tilde{f}^M(\tilde{x})$ is a vertex. Moreover, $\tilde{f}_\#^M(\tilde{\sigma}_1)$ is a proper subinterval of $f_\#^M(\tilde{\sigma})$: if not, then $f_\#^M(\tilde{\sigma}) = f_\#^l(E_i\mu^*)$ for some $l$ and some initial segment $\mu^*$ of $u_i$. This contradicts part (3) of Proposition 5.4.3, Lemma 4.1.1(5) and the assumption that $\sigma$ has no splittings. There are points of $\cap_{k=1}^M \tilde{S}_k$ in $\sigma_3'$ that map arbitrarily close to $\tilde{f}^M(\tilde{x})$. Since $\cap_{k=1}^M \tilde{S}_k$ is closed there exists $\tilde{y} \in \cap_{k=1}^M \tilde{S}_k$ in $\sigma_3'$ such that $\tilde{f}^M(\tilde{y}) = \tilde{f}^M(\tilde{x})$. The subdivision at $\tilde{x}$ and $\tilde{y}$ defines the desired $M$-splitting. This completes the first step.

If $\sigma = E_i\gamma$, then, after replacing $\sigma$ by some $f_\#^l(\sigma)$ if necessary, we may assume that the last (possibly partial) edge of $f_\#^k(\sigma)$ is contained in the same stratum for all $k \geq 0$. An immediate consequence of step 1 is that the last edge of $\sigma$ is not pointwise fixed by $f$. Thus one of the following conditions is satisfied:

(i) The terminal endpoint of $\sigma$ is a vertex and the terminal edge is some non-exponentially-growing $\bar{E}_j$ with nontrivial $u_j$.

(ii) The last edge of $\sigma$ is contained in an exponentially-growing stratum $H_r$.

If $\sigma = E_i\gamma\bar{E}_i$, then (i) holds for $E_i = E_j$ without replacing $\sigma$ by $f_\#^l(\sigma)$. Suppose at first that (i) holds.

*Step* 2 (at least 3 blocks cancel).    Write $\sigma = E_i \gamma' \bar{E}_j$, where $\gamma = \gamma'$ if $j = i$. Define the ray $R_i$ to be the infinite path $u_i \cdot f_\#(u_i) \cdot f_\#^2(u_i) \cdot \ldots$ and define $R_i^m$ to be the initial segment $u_i \cdot f_\#(u_i) \cdot f_\#^2(u_i) \cdot \ldots \cdot f_\#^{m-1}(u_i)$. We refer to the $f_\#^k(u_i)$'s as the blocks of $R_i$. Define $R_j$ and $R_j^m$ similarly with $u_j$ replacing $u_i$. Then $f_\#^m(\sigma) = [E_i R_i^m f^m(\gamma') \bar{R}_j^m \bar{E}_j]$. We claim that if $m$ is sufficiently large, then a subpath of $R_i^m$ containing at least three blocks of $R_i$ cancels with a subpath of $\bar{R}_j^m$ containing at least three blocks of $\bar{R}_j$ when $E_i R_i^m f^m(\gamma') \bar{R}_j^m \bar{E}_j$ is tightened to $f_\#^m(\sigma)$.

By step 1, there are a positive integer $M$ and initial subpaths $\mu_M$ of $R_i^M$ and $\bar{\nu}_M$ of $\bar{R}_j^M$ so that $f_\#^M(\sigma) = E_i \mu_M \bar{\nu}_M \bar{E}_j$. For $m > M$, $f_\#^m(\sigma) = f_\#^{m-M}(E_i \mu_M \bar{\nu}_M \bar{E}_j)$. Since $E_i \mu_M = f_\#^l(E_i \mu^*)$ for some initial segment $\mu^*$ of $u_i$, part (3) of Proposition 5.4.3 implies that $E_i \mu_M = E_i \cdot \mu_M$; similarly $\bar{\nu}_M \bar{E}_j = \bar{\nu}_M \cdot \bar{E}_j$. Thus $f_\#^m(\sigma)$ is obtained from the concatenation of $E_i R_i^{m-M} f_\#^{m-M}(\mu_M)$ and $f_\#^{m-M}(\bar{\nu}_M) \bar{R}_j^{m-M} \bar{E}_j$ by cancelling at the juncture point. Step 1 implies that for sufficiently large $m$, long cancellation must occur in both $R_i^{m-M}$ and $\bar{R}_j^{m-M}$. The only way that this could happen is if long segments of $R_i^{m-M}$ and $\bar{R}_j^{m-M}$ cancel with each other. This verifies our claim.

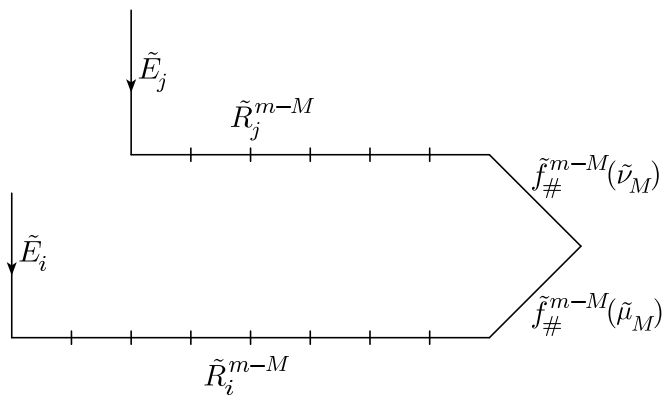

*Step* 3 ($H_s$ is non-exponentially-growing).    For future reference we formulate the third step as a sublemma. In step 5 we will apply the sublemma to $\delta = f_\#^m(\sigma)$ where $m$ is as in step 2.

SUBLEMMA 1.    *If $E_i$ and $E_j$ ($i \geq j$) have distinct lifts $\tilde{E}_i$ and $\tilde{E}_j$ whose corresponding rays $\tilde{R}_i$ and $\tilde{R}_j$ have a common subpath that contains at least three blocks in each ray and if $H_s$ is a non-exponentially-growing stratum, then the path $\tilde{\delta} \subset \Gamma$ that connects the initial endpoint of $\tilde{E}_i$ to the terminal endpoint of $\tilde{E}_j$ projects to an exceptional path $\delta \subset G$ of height $i$.*

*Proof of sublemma* 1. Let $\Gamma_{i-1} \subset \Gamma$ be the component of the full pre-image of $G_{i-1}$ that contains $\tilde{u}_i$ and $\tilde{u}_j$ and let $h_i : \Gamma_{i-1} \to \Gamma_{i-1}$ and $h_j : \Gamma_{i-1} \to \Gamma_{i-1}$ be the restricted lifts of $f$ that fix the initial endpoints of $\tilde{E}_i$ and $\tilde{E}_j$ respectively.

We consider first the case that $h_i = h_j$. Part (1) of Proposition 5.4.3 implies that $h_i$ is fixed-point free and hence that the initial endpoint of $\tilde{E}_j$ does not lie in $\Gamma_{i-1}$; thus $E_j = E_i$. The covering translation $S$ of $\Gamma_{i-1}$ that carries $\tilde{u}_i = \langle \tilde{v}, h_i(\tilde{v}) \rangle$ to $\tilde{u}_j = \langle S(\tilde{v}), h_j(S(\tilde{v})) \rangle$ commutes with $h_i$. Part (2) of Proposition 5.4.3 implies that $f_{\#}^k(u_i) = u_i$ for some $k > 0$ and hence by the inductive hypothesis that $f_{\#}(u_i) = u_i$; moreover, both $\tilde{R}_i$ and $\tilde{R}_j$ are contained in the axis of $S$. It follows that $u_i$ is a multiple of the indivisible circuit $\tau$ determined by $S$ and that the segment of the axis of $S$ that separates the terminal endpoints of $\tilde{E}_i$ and $\tilde{E}_j$ projects to $\tau^q$ or $\bar{\tau}^q$ for some $q \geq 0$. Thus $\delta = E_i \tau^q \bar{E}_i$ or $E_i \bar{\tau}^q \bar{E}_i$ is exceptional. This completes the proof in the special case and we assume from now on that $h_i \neq h_j$.

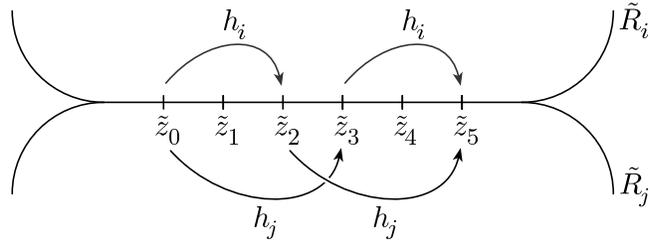

Let $\tilde{X} \subset \tilde{R}_i$ be the set of vertices that are either the initial endpoint of a lift of $E_s$ in $\tilde{R}_i$ or the terminal endpoint of a lift of $\bar{E}_s$ in $\tilde{R}_i$. Order the elements of $\tilde{X}$ so that $\tilde{x}_l < \tilde{x}_{l+1}$ in the orientation on $\tilde{R}_i$. Lemma 4.1.4 implies that $h_i(\tilde{x}_l) = \tilde{x}_{l+n_0}$ for all $l$ and some fixed $n_0$. Define $\tilde{Y} \subset \tilde{R}_j$ and $m_0$ similarly using $R_j$ and $h_j$ instead of $R_i$ and $h_i$. Then $\tilde{Z} = \tilde{X} \cap \tilde{Y} \subset \tilde{R}_i \cap \tilde{R}_j$ contains at least $n_0 + m_0 + 1$ consecutive elements $\tilde{z}_0, \ldots, \tilde{z}_{n_0+m_0}$ of $\tilde{X}$ and of $\tilde{Y}$ and $h_i h_j(\tilde{z}_0) = h_i(\tilde{z}_{m_0}) = \tilde{z}_{n_0+m_0} = h_j(\tilde{z}_{n_0}) = h_j h_i(\tilde{z}_0)$. Since $h_i h_j$ and $h_j h_i$ are lifts of $f^2$ that agree at a point, they are equal.

There is a nontrivial covering translation $S$ of $\Gamma_{i-1}$ such that $Sh_i = h_j$ and there is a covering translation $T$ of $\Gamma_{i-1}$ such that $Th_j = h_i S$. Then $h_i h_j = h_i S h_i = T h_j h_i$ so $T$ is the identity and $h_i$ commutes with $S$. A symmetric argument shows that $h_j$ also commutes with $S$. Part (2) of Proposition 5.4.3 implies that $u_i, u_j \subset G_{i-1}$ are periodic Nielsen paths and that $\tilde{R}_i$ and $\tilde{R}_j$ are contained in the axis of $S$. By the inductive hypothesis, $f_{\#}(u_i) = u_i$ and $f_{\#}(u_j) = u_j$. The covering translation $S_i$ of $\Gamma_{i-1}$ that carries the initial endpoint of $\tilde{u}_i$ to the terminal endpoint of $\tilde{u}_i$ and the covering translation $S_j$ of $\Gamma_{i-1}$ that carries the initial endpoint of $\tilde{u}_j$ to the terminal endpoint of $\tilde{u}_j$ both preserve the axis of $S$ and hence commute with $S$ and with each other. The segment of the axis of $S$ that separates the terminal endpoints of $\tilde{E}_i$ and $\tilde{E}_j$ projects to $\tau^q$ or $\bar{\tau}^q$ where $\tau$ is the circuit corresponding to the axis of $S$ and $q \geq 0$. After replacing $S$ with $S^{-1}$ if necessary, $u_i = \tau^k$ and $u_j = \tau^l$ or $\bar{\tau}^l$. Since $\tilde{R}_i$ and $\tilde{R}_j$ have a common subpath that contains blocks in both rays, $u_j = \tau^l$ and $\delta$ is exceptional.    □

*Step* 4 ($H_s$ is exponentially growing).    Suppose that $m$ is as in step 2 and that $H_s$ is exponentially growing. We will prove that $f^m_\#(\sigma)$ is an exceptional path of height $i$. Let $h_i : \Gamma_{i-1} \to \Gamma_{i-1}$ and $h_j : \Gamma_{i-1} \to \Gamma_{i-1}$ be as in step 3.

If the decomposition of $u_i$ given by Lemma 4.2.6 contains at least one $\rho_s$ or $\bar{\rho}_s$, then the argument of step 3 requires only one modification. Namely, $\tilde{X}$ is defined to be the set of lifts of $\rho_s$ or $\bar{\rho}_s$. These get 'translated' by $h_i$ and $h_j$ so the argument goes through exactly as before.

It remains to rule out the possibility that the decomposition of $u_i$ given by Lemma 4.2.6 contains no $\rho_s$ or $\bar{\rho}_s$. Suppose to the contrary. Then $u_i$ is *s*-legal. Since blocks of $\bar{R}_j$ cancel with segments of $R_i$, $u_j$ is also *s*-legal. The action of $f$ on $R_i$ and $R_j$ is more like an affine map than like a translation so we do not use the previous argument. By step 2, $f^m_\#(\sigma) = E_i \mu_m \bar{\nu}_m \bar{E}_j$ where $\mu_m$ and $\bar{\nu}_m$ are initial subpaths of $R_i$ and $\bar{R}_j$ and in particular are *s*-legal. Lemma 4.2.2 implies that $\mu_m$ and $\nu_m$ take on only finitely many values and hence that some $f^m_\#(\sigma)$ is a periodic Nielsen path, say $f^{m+p}_\#(\sigma) = f^m_\#(\sigma)$. There is a lift $\tilde{f} : \Gamma \to \Gamma$ whose restriction to $\Gamma_{i-1}$ equals $h_i$. Thus $\tilde{f}$ fixes the initial endpoint of $\tilde{E}_i$ and $\tilde{f}^p$ fixes the initial endpoint $\tilde{w}$ of $\tilde{E}_j$.

If $E_i \neq E_j$, then $\tilde{w} \in \Gamma_{i-1}$. Part (1) of Proposition 5.4.3 implies that $p > 1$. Let $\tilde{\gamma} \subset \Gamma$ be the path that connects $\tilde{w}$ to $h_i(\tilde{w})$ and let $\gamma$ be its image in $G_{i-1}$. Then $\gamma$ is a periodic Nielsen path and so by the inductive hypothesis is a Nielsen path. But then $[\gamma^p] = [\gamma f(\gamma) \dots f^{p-1}(\gamma)]$ lifts to the trivial path $[\tilde{\gamma} h_i(\tilde{\gamma}) \cdot \dots \cdot h_i^{p-1}(\tilde{\gamma})]$ which is impossible. We conclude that $E_i = E_j$. The covering translation $S$ of $\Gamma$ that carries the initial endpoint of $\tilde{E}_i$ to the initial endpoint of $\tilde{E}_j$ commutes with $\tilde{f}^p$. The restriction $S|\Gamma_{i-1}$ therefore commutes with $h_i^p$. Since $h_j^p = (S|\Gamma_{i-1})h_i^p(S|\Gamma_{i-1})^{-1}$, $h_i^p = h_j^p$. This implies that $\tilde{R}_i$ and $S(\tilde{R}_i) = \tilde{R}_j$ have an infinite end in common and hence that $\tilde{R}_i$ and the axis of $S$ have an infinite end in common. It follows that $h_i$ preserves the axis of $S$ and so commutes with $S$. This contradicts part (2) of Proposition 5.4.3 and the fact that $u_i$ is *s*-legal.

*Step* 5 ((2) and (3) are satisfied when (i) holds).    Choose a lift $\tilde{E}_i$ in the universal cover $\Gamma$ and choose $m$ as in step 2. There is a lift of $f^m_\#(\sigma)$ that begins with $\tilde{E}_i$ and ends at the inverse of some $\tilde{E}_j$. Let $\tilde{R}_i$ and $\tilde{R}_j$ be the lifts of $R_i$ and $R_j$ that begin at the terminal endpoints of $\tilde{E}_i$ and $\tilde{E}_j$ respectively. Our choice of $m$ guarantees that $\tilde{R}_i \cap \tilde{R}_j$ contains at least three blocks in each ray. By steps 3 and 4, $f^m_\#(\sigma)$ is an exceptional path of height $i$. If $i = j$, then $f_{\#}^m(\sigma)$ is fixed by $f_\#$. Since $\sigma$ and $f_{\#}^m(\sigma)$ have the same endpoints and the same image under $f_{\#}^m$, they must be equal. In particular, $\sigma$ is an exceptional path of height $i$.

*Step* 6 (Case (ii) does not occur).   Suppose that (ii) holds. Step 1 implies that for all sufficiently large $m$, $f_{\#}^m(E_i\gamma) = E_i\mu_m\nu_m$ where $\mu_m \subset R_i$ and where $\nu_m \subset G_r$ is $r$-legal. Step 1 also implies that if $(l - m)$ is sufficiently large, then $H_r$-edges of $f_{\#}^{(l-m)}(\nu_m)$ cancel with edges of $f_{\#}^{(l-m)}(\mu_m)$ when $f_{\#}^{(l-m)}(E_i\mu_m)f_{\#}^{(l-m)}(\nu_m)$ is tightened to $f_{\#}^l(\sigma)$; thus $r \leq s$. The symmetric argument shows that $s \leq r$ (so that $s = r$) and that $\mu_m$ must be $s$-legal. Lemma 4.2.2 implies that $\mu_m$ takes on only finitely many values. The same is true for $\nu_m$: the argument is essentially the same as the one in the proof of Lemma 4.2.5. There are no subtleties in applying this argument and we leave the details to the reader. We conclude that $f_{\#}^m(\sigma)$ is a periodic Nielsen path. But the argument of the last paragraph in the proof of step 4 proves that this is impossible. We conclude that (ii) does not occur.

*Step* 7 (Conclusion).   It remains to prove that every periodic Nielsen path $\rho \subset G_i$ has period one. Lemma 4.1.4 implies that $\rho$ splits into periodic Nielsen paths that are either contained in $G_{i-1}$ or are basic paths. Conditions (1)–(3) of this lemma (which we have already verified), imply that a basic path that is also a periodic Nielsen path splits into exceptional paths of height $i$ and periodic Nielsen paths in $G_{i-1}$. By induction and by examination of the exceptional paths, we conclude that $\rho$ has period one.   □

5.6. *Proof of Theorem* 5.1.5.  Choose (Lemma 5.2.3) a reduced $\mathcal{F}$-Nielsen minimizing relative train track map $f : G \to G$ and filtration $\emptyset = G_0 \subset G_1 \subset \cdots \subset G_K = G$ that represents $\mathcal{O}^k$ for some $k \geq 1$. We will modify $f : G \to G$ in various ways but will continue to call the resulting relative train track map $f : G \to G$. (Modifications of $f : G \to G$ may of course involve changes in $G$.)

Condition (5) in the definition of $\mathcal{F}$-Nielsen minimizing (Definition 5.2.1) implies that contractible components of $G_i$'s are unions of zero strata.

Lemma 5.2.5, Lemma 5.1.7 and Proposition 5.3.1 imply (eg-i), (eg-ii) and (eg-iii). If $H_r$ is an exponentially-growing stratum and $\rho \subset G_r$ is an indivisible periodic Nielsen path that intersects $H_r$ nontrivially, then since $f : G \to G$ is $\mathcal{F}$-Nielsen minimizing, $\rho$ has period one. These properties are unchanged by taking iterates.

After passing to an iterate if necessary, we may assume: that $f(v)$ is a fixed point for each vertex $v$; that the noncontractible components of the $G_i$'s are mapped to themselves by $f$; and that each non-exponentially-growing stratum $H_i$ consists of a single edge $E_i$ satisfying $f(E_i) = w_iE_iv_i$ for some paths $w_i, v_i$ in $G_{i-1}$. After subdividing at a fixed point in the interior of $E_i$ if necessary, we may assume that $w_i$ is the trivial path.

Suppose that $H_i$ is exponentially growing, that $C$ is a noncontractible component of $G_{i-1}$ and that $v \in H_i \cap C$. Choose a path $\alpha$ between $v$ and $f(v)$. Since $f|C : C \to C$ is a homotopy equivalence of $C$, there is a closed path $\beta$ based at $f(v)$ so that $f_\#(\beta) = f_\#(\alpha)$. Then $\delta = [\alpha\bar{\beta}]$ is a path between $v$ and $f(v)$ such that $f_\#(\delta)$ is trivial. The train track property for $H_i$ implies that $\delta$ is trivial and hence that $v = f(v)$ is a fixed point.

Apply Proposition 5.4.3 to the non-exponentially-growing strata in the filtration for $f : G \to G$ working upwards. At the end of this process, we have lost none of our previously acquired properties and arranged that if $H_i$ is a non-exponentially-growing stratum then $f(E_i) = E_i \cdot u_i$. Moreover, the endpoints of the edges in non-exponentially-growing strata are now periodic points; after passing to a further iterate, they are fixed points. We have therefore established (ne-ii). Suppose that $H_i$ is non-exponentially growing and that $\sigma = E_i\gamma, \gamma\bar{E}_i$ or $E_i\gamma\bar{E}_i$ is a basic path of height $i$ (Definition 4.1.3). If $\sigma$ splits at some point of $\gamma$, then $\sigma$ splits as a concatenation of either two basic paths of height $i$ or a basic path of height $i$ and a path in $G_{i-1}$. Lemma 5.5.1 implies that if $\sigma$ splits at a point in $E_i$ or $\bar{E}_i$ then some $f_\#^k(\sigma)$ splits as a concatenation of two subpaths one of which is $E_i$ or $\bar{E}_i$. Lemma 5.5.1 also implies that if $\sigma$ has no splittings, then some $f_\#^k(\sigma)$ is an exceptional path of height $i$. Thus condition (ne-iii) is satisfied. Applying Lemma 5.5.1 once again we see that every periodic Nielsen path has period one. The properties that we have verified so far are all stable under passing to an iterate of $f$.

We now turn our attention to the zero strata. We have already established that a contractible component of any $G_i$ is made up of edges in zero strata. We next show that if $H_i$ is a zero stratum, then $H_i$ is contained in the union of the contractible components of $G_i$. Suppose to the contrary that a zero stratum $H_i$ intersects a noncontractible component $C$ of $G_i$ and that $i$ is the largest positive integer for which this occurs. Then $f(C) \subset \mathrm{cl}(C \setminus (H_i \cap C))$ and so $\mathrm{cl}(C \setminus (H_i \cap C))$ is a proper subgraph of $C$ that has the same rank as $C$. We may therefore choose a vertex $v$ of $H_i$ that has valence one in $C$. Since $G$ does not have valence one vertices, $v$ is incident to an edge that is part of a higher stratum and so by assumption is not a zero stratum. But we have already shown that such vertices are fixed points. This contradicts the fact that the only edge incident to $v$ in $C$ maps off of itself. We conclude that a zero stratum $H_i$ is contained in the union of the contractible components of $G_i$.

We next reorganize the zero strata and push them up the filtration as high as possible. Assume that $H_i$ is a zero stratum. Since vertices in zero strata are not fixed points and so are not the $f$-image of a vertex, no edge in $G$ has $f$-image entirely contained in zero strata. In particular, the contractible components of $G_i$ are disjoint from $f(G_i)$. We may therefore amalgamate the edges in all of the contractible components of $G_i$ into a single zero stratum (still called $H_i$). If $H_{i+1}$ is a zero stratum, then we can amalgamate $H_i$ and $H_{i+1}$ into a single stratum (still called $H_i$). We may therefore assume that $H_{i+1}$ is not a zero stratum. If some components of $H_i$ are components of $G_{i+1}$, then they are not in the image of $G_{i+1}$, and we can remove these components from $H_i$ and consider them as a new zero stratum $H_{i+2}$. We may therefore assume that $G_{i+1}$ has no contractible components. Since the endpoints of edges in non-exponentially-growing strata are fixed points, $H_{i+1}$ must be an exponentially-growing stratum.

After performing these operations on each zero stratum, working upward through the filtration, we see that condition z-(i) is satisfied and $H_i$ is a zero stratum if and only if it is the union of the contractible components of $G_i$.

If $H_i$ is a zero stratum and $f|H_i$ is not an immersion, then we can fold a pair of edges in $H_i$. Tighten the images of the remaining edges. Since the $f$-image of an edge does not lie entirely in zero strata, no edge has trivial image after tightening. Folding and tightening in this manner do not undo any of our established properties (cf. 4.3.6 of [BH2]) and they reduce the total number of edges in the image of $H_i$. After finitely many steps, $f|H_i$ is an immersion. Perform this folding operation on each zero stratum, working up through the filtration so that the modifications made in one zero stratum do not undo the modifications made in the previous strata. At the end of this process, condition z-(ii) is satisfied.

If $v$ is a vertex in a zero stratum $H_i$, then it is not fixed by $f$ and so is not the endpoint of an edge in $H_j$ with $j > i + 1$. Since $G$ has no valence one vertices, $v$ must have valence at least two in $G_{i+1}$. If $v$ has valence two and if both incident edges are contained in $H_i$, then erase $v$ as a vertex. Since $v$ is not the image of any vertex the map remains simplicial. After erasing all such vertices, z-(iii) is satisfied.                                                                                                   □

5.7. **UPG**$(F_n)$. We assume in this section that $\mathbb{Z}^n$ is identified with the abelianization of $F_n$ and hence also with $H_1(G; \mathbb{Z})$ for any marked graph $G$. There is an induced homomorphism from $\mathrm{Out}(F_n)$ to $\mathrm{GL}(n, \mathbb{Z})$.

We also assume in this section that $\mathcal{F}$ is an $\mathcal{O}$-invariant free factor system.

*Definitions* 5.7.1. Denote $\{\mathcal{O} \in \mathrm{Out}(F_n) : \mathcal{L}(\mathcal{O}) = \emptyset\}$ by $\mathrm{PG}(F_n)$. Thus $\mathcal{O} \in \mathrm{PG}(F_n)$ if and only if some, and hence every, relative train track map representing $\mathcal{O}$ has no exponentially-growing strata. Recall that an element of $\mathrm{GL}(n, \mathbb{Z})$ is *unipotent* if it is conjugate to an upper triangular matrix with 1's on the diagonals. The subset of $\mathrm{PG}(F_n)$ consisting of elements $\mathcal{O}$ whose image in $\mathrm{GL}(n, \mathbb{Z})$ is unipotent is denoted $\mathrm{UPG}(F_n)$.

The following lemma shows that zero strata are not needed for elements of $\mathrm{PG}(F_n)$.

LEMMA 5.7.2.  *Every $\mathcal{O} \in \mathrm{PG}(F_n)$ is represented by a relative train track map $f : G \to G$ and filtration $\emptyset = G_0 \subset G_1 \subset \cdots \subset G_K = G$ such that:*

1. *$\mathcal{F} = \mathcal{F}(G_r)$ for some filtration element $G_r$.*

2. *Each vertex of each $G_i$ has valence at least two. In particular, there are no zero strata and all components of $G_i$ are noncontractible.*

*Proof.* Lemma 2.6.7 provides a relative train track map $f : G \to G$ and filtration $\emptyset = G_0 \subset G_1 \subset \cdots \subset G_K = G$ satisfying condition 1. Since $\mathcal{O} \in \mathrm{PG}(F_n)$, $f : G \to G$ has no exponentially-growing strata. If $v$ is a valence one vertex of some $G_i$ and if $E$ is the unique edge of $G_i$ that is incident to $v$, then perform a homotopy of $f$ by precomposing $f$ with the homotopy that slides $v$ across $E$ to the other endpoint of $E$ and then tightening. This homotopy only affects edges that are incident to $v$. If $E' \neq E$ is an edge of $H_l$ that is incident to $v$, then the homotopy changes the way $f(E')$ crosses edges in $G_i$, but the relative train track property is maintained since $l > i$. (Keep in mind that a topological representative of an element of $\mathrm{PG}(F_n)$ is a relative train track map if and only if it has no exponentially-growing strata.) The new image of $E$ is trivial so we can collapse $E$. This does not affect the relative train track property and does not change any $[[\pi_1(G_j)]]$, so that condition 1 is still satisfied. After finitely many such moves, condition 2 is satisfied. □

*Definition* 5.7.3.  Suppose that $f : G \to G$ and $\emptyset = G_0 \subset G_1 \subset \cdots \subset G_K = G$ are a relative train track map and filtration representing $\mathcal{O} \in \mathrm{PG}(F_n)$ and that there are no zero strata in the filtration. Choose a maximal tree $T$ for $G$ whose intersection $T_i$ with each $G_i$ is a maximal forest for $G_i$. For each edge $e$ in $G_i \setminus T_i$ choose an embedded circuit $\gamma \subset G_i$ that contains $e$ but is otherwise contained in $T_i$. We say that the basis $B$ for $H_1(G; \mathbb{Z})$ determined by the homology classes of the $\gamma$'s is *a PG basis determined by $f : G \to G$.* The element of $\mathrm{GL}(n; \mathbb{Z})$ determined by $B$ and the action of $f$ on $H_1(G, \mathbb{Z})$ is denoted $M_B$.

LEMMA 5.7.4.  *If $B$ is a PG basis determined by $f : G \to G$, then each diagonal entry in $M_B$ is either $-1, 0$ or $1$.*

*Proof.* An element $b \in B$ corresponds to an edge $e$ in some $G_i \setminus T_i$ and an embedded circuit $\gamma \subset G_i$. The lemma follows from the fact that $f_\#(\gamma)$ either crosses $e$ once or not at all. □

The following proposition relates the UPG property to relative train track maps.

PROPOSITION 5.7.5.   *The following are equivalent*:

(1)  $\mathcal{O} \in \mathrm{UPG}(F_n)$.

(2)  *If $B$ is a PG basis determined by $f : G \to G$, then each diagonal entry in $M_B$ equals* 1.

(3)  *There is relative train track map $f : G \to G$ and filtration $\emptyset = G_0 \subset G_1 \subset \cdots \subset G_K = G$ representing $\mathcal{O}$ such that*:

(3-a)  $\mathcal{F} = \mathcal{F}(G_r)$ *for some filtration element $G_r$.*

(3-b)  *Each vertex of each $G_i$ has valence at least two.*

(3-c)  *Each $H_i$ is a single edge $E_i$ satisfying $f(E_i) = v_i E_i u_i$ for paths $v_i, u_i \subset G_{i-1}$.*

*Proof.* If (3) is satisfied and $B$ is a PG basis determined by $f : G \to G$ and $\emptyset = G_0 \subset G_1 \subset \cdots \subset G_K = G$ then $M_B$ is upper triangular with 1's on the diagonal. Thus (3) $\Rightarrow$ (1).

If (1) is satisfied, then $M_B$ is conjugate to an upper triangular matrix with 1's on the diagonal and so has trace $n$. Condition (2) follows from Lemma 5.7.4. Thus (1) $\Rightarrow$ (2).

The remainder of the proof is dedicated to proving that (2) $\Rightarrow$ (3). Lemma 5.7.2 provides $f : G \to G$ and $\emptyset = G_0 \subset G_1 \subset \cdots \subset G_K = G$ satisfying (3-a) and (3-b). We will arrange (3-c) by induction on $i$. Note that if $\mathcal{O} \in \mathrm{UPG}(F_n)$ and $B$ is a PG basis for $f : G \to G$, then $\mathcal{O}^k \in \mathrm{UPG}(F_n)$ and $B$ is a PG basis for $\mathcal{O}^k$ with matrix $M_B^k$. We may therefore apply (2) to iterates of $f$.

Let $B$ be a PG basis for $f : G \to G$ . Since $f$ transitively permutes the edges in $G_1$, conditions (3-b) and (2) imply that $G_1$ is connected. Choose an element $b \in B$, let $e$ be the corresponding edge in $G_1$ and let $\gamma$ be the corresponding embedded circuit in $G_1$. If $G_1$ has rank at least two, then there is an edge $e' \subset G_1$ that is not in $\gamma$. After replacing $f$ by an iterate if necessary, we may assume that $f(e') = e$. But then the circuit $f(\gamma)$ does not contain $e$ and so the diagonal element corresponding to $b$ is 0. This contradicts (2) and we conclude that $G_1$ has rank one. By (3-b), $G_1$ is homeomorphic to a circle.

If $G_1$ contains more than one edge, then $f|G_1$ is a nontrivial rotation. Choose a homotopy $h_t : G \to G$ with support in a small neighborhood of $G_1$ so that $h_0 = $ identity and $h_1|G_1 = (f|G_1)^{-1}$. Define a new topological representative of $G$ by tightening $f h_1$. At this point, the map (which we still call $f$) fixes each edge in $G_1$; collapse all but one of these edges to arrange for $G_1$ to be a single edge $E_1$. As in the proof of Lemma 5.7.2, $f : G \to G$ is a relative train track map and (3-a) is still satisfied as is condition (3-b). Now (2) rules out the possibility that $f(E_1) = \bar{E}_1$.

We assume now that (3-c) holds for $G_{i-1}$ and prove it for $G_i$. Suppose that $H_i = \{E^1, \ldots, E^p\}$. (We reserve the notation $E_j$ for the edges of $G$ after (3-c) has been satisfied.) Then $f(E^j) = v_j E^{j+1} u_j$ or $f(E^j) = v_j \bar{E}^{j+1} u_j$ for subpaths $v_j, u_j \subset G_{i-1}$, where the indices of the $E^j$'s are taken mod $p$ and where some reordering may have been necessary. If the $E^j$'s are disjoint from $G_{i-1}$, then they determine a component of $G_i$ and we may proceed exactly as in the $G_1$ case. We may therefore assume that some, and hence every, $E^j$ has at least one in endpoint $G_{i-1}$. Let $C_1$ be a component of $G_{i-1}$ that contains an endpoint of each $E^j$. There are three cases to consider, depending on the location of the other endpoint of $E^j$.

Suppose at first that both endpoints of each $E^j$ lie in $C_1$. In this case, each $E^j$ determines an element of $B$. The embedded circuit $\gamma$ corresponding to $E^1$ is the concatenation of $E^1$ and a subpath in $G_{i-1}$. If $p > 1$ then the circuit $f_\#(\gamma)$ is therefore a concatenation of $E_2$ or $\bar{E}_2$ with a subpath in $G_{i-1}$. But then the diagonal element associated to $E^1$ is 0 in contradiction to (2). If $p = 1$ and $f(E^1) = u_1 \bar{E}_1 v_1$ then the diagonal element associated to $E^1$ is $-1$. This also contradicts (2) so (3-c) is satisfied.

Suppose next that each $E^j$ has an endpoint in a component $C_2 \neq C_1$ of $G_{i-1}$. In this case $f(E^j) = v_j E_{j+1} u_j$ and we need only show that $p = 1$. We may assume that the intersection of the maximal tree $T$ with $H_i$ is $E^1$. If $p > 1$, then the embedded circuit $\gamma$ corresponding to $E^p$ intersects $H_i$ is $\bar{E}^1$ and $E^p$. The image circuit $f_\#(\gamma)$ intersects $H_i$ in $\bar{E}_2$ and $E^1$. Thus the diagonal entry of $M_B$ corresponding to $E_2$ is 0 if $p > 2$ and $-1$ if $p = 2$. This contradiction to (2) verifies that $p = 1$.

Finally, we rule out the possibility that $E^1$ has an endpoint $x$ that is not in $G_{i-1}$. Since $G_i$ does not have valence-one vertices, there must be at least one other edge of $H_i$, say $E^2$, with an endpoint at $x$. We may assume that $E^1 \subset T$ but that no other $E^j$ with an endpoint at $x$ is contained in $T$. As in the previous case, $E^2$ determines an element of $B$ whose diagonal entry in $M_B$ is either 0 or $-1$. This contradiction to (2) completes the proof. $\square$

COROLLARY 5.7.6.   *If $\mathcal{O} \in \mathrm{PG}(F_n)$ is contained in the kernel of the natural homomorphism*

$$\mathrm{Out}(F_n) \to \mathrm{GL}(n; \mathbb{Z}) \to \mathrm{GL}(n, \mathbb{Z}/3\mathbb{Z}),$$

*then $\mathcal{O} \in \mathrm{UPG}(F_n)$. In particular, every subgroup of $\mathrm{PG}(F_n)$ contains a finite index subgroup in $\mathrm{UPG}(F_n)$.*

*Proof.* This is an immediate consequence of Lemma 5.7.4, condition (2) of Proposition 5.7.5 and the obvious fact that $-1$ and 0 are not congruent to 1 mod 3. $\square$

COROLLARY 5.7.7.   *Suppose that $\mathcal{O} \in \mathrm{UPG}(F_n)$ and that $\mathcal{F}$ is an $\mathcal{O}$-invariant free factor system.   Then there are a relative train track map $f : G \to G$   and filtration $\emptyset = G_0 \subset G_1 \subset \cdots \subset G_K = G$ that represent $\mathcal{O}$ and that satisfy conditions (1), (2), (3) and (4)   of Theorem 5.1.8.*

*Proof.* Choose $f : G \to G$ and $\emptyset = G_0 \subset G_1 \subset \cdots \subset G_K = G$ as in part (3) of Proposition 5.7.5. Then Condition (1) is satisfied and all periodic points of $f$ are fixed points. We arrange (2) by induction on $i$ as follows. The $i = 1$ case follows from the fact that $f|G_1 = $ identity. Suppose then that (2) holds for $1 \leq j \leq i - 1$. If both $u_i$ and $v_i$ are nontrivial, then subdivide $E_i$ at the unique fixed point in the interior of $E_i$ to create two edges with either $u_i$ or $v_i$ trivial. If $v_i$ is not trivial, replace $E_i$ by $\bar{E}_i$. We may therefore assume that $v_i$ is trivial.   Apply Proposition 5.4.3 to arrange that $f(E_i) = E_i \cdot u_i$ where the initial endpoint of $u_i$ is periodic and hence fixed. Since $f|G_{i-1}$ is unchanged, (2) is now satisfied for $1 \leq j \leq i$. This completes the induction step.

Condition (2) implies condition (3). Lemma 5.5.1 implies condition (4).   □

A useful corollary of condition (3) of Theorem 5.1.8 is that for any path $\tau$ with endpoints at vertices, there is a unique path $\sigma$ with endpoints at vertices such that $f_\#(\sigma) = \tau$. Since every exceptional path is the image of an exceptional path, a nonexceptional path cannot have an image that is an exceptional path.

We say that a path $\sigma$ with endpoints at vertices has *height $i$* if it crosses $E_i$ but does not cross $E_j$ for any $j > i$. By Lemma 4.1.4, every path of height $i$ with endpoints at vertices has a splitting whose pieces are either basic paths of height $i$ or paths with height less than $i$. This provides a recursive splitting of $\sigma$ into basic pieces. We will use the following notion of complexity as the basis for induction arguments.

*Definition* 5.7.8. We subdivide basic paths of height $i$ into two types. Those of the form $E_i\gamma$ and $\gamma\bar{E}_i$ are called *type 1 basic paths of height $i$* and those of the form $E_i\gamma\bar{E}_i$ are called *type 2 basic paths of height $i$*. The *complexity of a basic path* is the ordered pair specifying its height and its type; the pairs are ordered lexicographically. Thus a type 2 basic path of height $i$ has greater complexity than a type one basic path of height $i$ and lower complexity than a type one basic path of height $i + 1$. If $\sigma$ has height $i$ and $\sigma = \sigma_1 \cdot \ldots \cdot \sigma_l$ is the splitting of Lemma 4.1.4, then the *complexity of $\sigma$* is the highest complexity that occurs among those $\sigma_j$ that are basic paths of height $i$.

*Proof of Theorem* 5.1.8. Choose $f : G \to G$ and $\emptyset = G_0 \subset G_1 \subset \cdots \subset G_K = G$ as in Corollary 5.7.7.

*Step* 1 (Property (5)).   We prove (5) by induction on the complexity of $\sigma$ and assume without loss of generality that $\sigma$ is a basic path of height $i$ and that $\sigma = E_i\gamma$ or $\sigma = E_i\gamma_i\bar{E}_i$. In either case, if $\sigma = \sigma_1 \cdot \sigma_2$ is a splitting into subpaths with endpoints at vertices, then $\sigma_2$ has strictly smaller complexity than $\sigma$.

Since $f|G_1$ is the identity, $M(\sigma) = 0$ for $\sigma$ with height 1. We may therefore assume that (5) holds for paths with complexity lower than $\sigma$. In particular, $M(E_j) = 0$ for all $j < i$. After replacing $u_i$ by $f_\#^{M(u_i)}(u_i)$ if necessary, (which requires sliding $E_j$ as in the proof of Proposition 5.5.1) we may assume that $M(u_i) = 0$ and hence that $M(E_i) = 0$.

If $\sigma$ has no splittings, then parts (2) and (3) of Lemma 5.5.1 imply that some $f_\#^k(\sigma)$ is exceptional and so $M(\sigma) \leq k$. Suppose then that $\sigma$ can be split. If $\sigma$ can be split at some point in its initial edge, then part (1) of Lemma 5.5.1 implies that $f_\#^k(\sigma) = E_i \cdot \beta$ for some path $\beta$ and some $k \geq 0$. Thus $M(\sigma) \leq k + M(\beta)$ and, since $\beta$ has lower complexity than $\sigma$, induction completes the proof. Finally, if $\sigma$ cannot be split at any point in its initial edge, then split as close to the initial vertex as possible. This yields $\sigma = E_i\mu_1 \cdot \mu_2$ where $E_i\mu_1$ has no splittings. Parts (2) and (3) of Lemma 5.5.1 imply that some $f_\#^k(\sigma) = \alpha \cdot \beta$ where $\alpha$ is an exceptional path. Thus $M(\sigma) \leq k + M(\beta)$ and induction completes the proof of (5). $\qquad\square$

*Step* 2 (extending rays to lines). Let $s = \text{height}(u_i)$. For $l \geq 0$, denote $f_\#^l(u_i)$ by $B_l$ and define $B_{-l} \subset G_s$ to be the unique path with endpoints at vertices such that $f_\#^l(B_{-l}) = u_i$. Clearly $f_\#(B_l) = B_{l+1}$ for all $l$. For the remainder of this proof we will refer to $E_s$ and $\bar{E}_s$ as $s$-edges. Since $f_\#^k$ preserves highest edges, the number of $s$-edges in $B_l$ and the number of $s$-edges in $[B_lB_{l+1}]$ are independent of $l$. Since $[B_0B_1] = B_0B_1$, each $[B_lB_{l+1}]$ has twice as many $s$-edges as each $B_l$. In particular, no $s$-edges are canceled when $B_lB_{l+1}$ is tightened to $[B_lB_{l+1}]$. By construction (see the proof of Proposition 5.4.3) $B_0$, and hence each $B_l$, either begins with $E_s$ or ends with $\bar{E}_s$. Thus $f_\#(B_lB_{l+1}) = f_\#(B_l)f_\#(B_{l_1})$ and we have shown that $B_lB_{l+1} = B_l \cdot B_{l+1}$ for all $l$. The union of the $B_l$'s is an $f_\#$-invariant line in $G_s$.

We say that a bound independent of $\sigma$ is a uniform bound. As in the proof of (5), we argue by induction on complexity and the height-one case is obvious. We may therefore assume that (6) holds for paths with complexity lower than $\sigma$. There is no loss in assuming that $u_i$ is nontrivial and that $\sigma = E_i\gamma$ or $\sigma = E_i\gamma\bar{E}_i$ for some $\gamma \subset G_{i-1}$. To handle both of these cases simultaneously we write $\sigma = E_i\mu$ where either $\mu = \gamma$ or $\mu = \gamma\bar{E}_i$. In either case, the complexity of $\mu$ is smaller than the complexity of $\sigma$. We argue by induction on $j = \text{height}(\mu)$.

*Step* 3 (The case $j < s$). Suppose that $j < s$. If $k > 0$ then $B_{k-1}$ is not entirely canceled when $f^k(\sigma) = E_iB_0 \ldots B_{k-1}f^k(\mu)$ is tightened to $f_\#^k(\sigma)$. It follows (Lemma 4.1.1(6)) that $f_\#(\sigma) = E_i \cdot \mu'$ where $\mu' = [B_0f(\mu)]$ and hence that $M(\sigma) \leq M(\mu') + 1$. Since the length of $\mu'$ is a uniformly bounded multiple of the length of $\mu$, induction on complexity completes the proof.

*Step* 4 (The case $j = s$). Suppose that $j = s$. We will use the observation that $f_{\#}^k(\sigma) = E_i[B_0 \ldots B_{k-1} f^k(\mu)] = E_i f_{\#}^k([B_{-k} \ldots B_{-1}\mu])$. Since $f_{\#}^k$ preserves highest edges, any cancellation of $s$-edges when $f^k(\sigma)$ is tightened to $E_i f_{\#}^k([B_{-k} \ldots B_{-1}\mu])$ occurs when $B_{-k} \ldots B_{-1}\mu$ is tightened to $[B_{-k} \ldots B_{-1}\mu]$. We now restrict to $k$ so large that $B_{-k} \ldots B_{-1}$ contains more $s$-edges than $\mu$ does. Decompose $\mu$ into subpaths $\mu = \mu_1\mu_2$ where $\mu_1$ is the shortest initial segment that contains each $s$-edge of $\mu$ that is canceled when $B_{-k} \ldots B_{-1}\mu$ is tightened. Then $\mu_1 = \bar{B}_{-1}\bar{B}_{-2} \ldots \bar{B}_{-(r-1)}\bar{B}_{-r}^*$ where $\bar{B}_{-r}^*$ is an initial segment of $\bar{B}_{-r}$ and $r \geq 0$. Now $\sigma$ can be obtained from $(E_i\bar{B}_{-1} \ldots \bar{B}_{-r-1})[B_{-r-1} \ldots B_{-1}\mu]$ by tightening at the indicated juncture. Applying $f_{\#}^k$, we see that $f_{\#}^k(\sigma)$ is obtained from $(E_i \cdot B_0 \cdot B_1 \cdot \ldots \cdot B_{k-r-2})f_{\#}^k([B_{-r-1} \ldots B_{-1}\mu])$ by cancelling at the indicated juncture. By construction, the latter term begins with $B_{k-r-1}$ so there is no cancellation at the indicated juncture. Lemma 4.1.1(6) implies that $f_{\#}^{r+2}(\sigma) = E_i \cdot B_0 \cdot f_{\#}^{r+2}([B_{-r-1} \ldots B_{-1}\mu])$ and hence that $M(\sigma) \leq \max\{r + 2, M([B_{-r-1} \ldots B_{-1}\mu])\}$. The path $B_{-r-1} \ldots B_{-1}\mu$ is obtained from $\mu$ by removing $\mu_1$, adding $B_{-r-1}$ and perhaps adding part of $B_{-r}$. If $r \leq 1$, then the additional edges have a uniformly bounded length. If $r > 1$, then $\mu_1$, and hence $\mu$, contains $\bar{B}_{-r+1}$ and so the length of the additional edges is a uniformly bounded multiple of the length of $\mu$ . In either case, induction on complexity completes the proof.

*Step* 5 (The case $j > s$). We now assume that $j > s$. In this case the splittings of $\mu$ produced by Lemma 4.1.4 extend to splittings of $\sigma$. Since all but the first subpath in any such splitting of $\sigma$ have lower complexity than $\sigma$, we may assume that $\sigma = E_i\nu\bar{E}_j$ where height$(\nu) < j$. For the same reason, we may also assume that $\sigma$ does not split at any vertex. Let $t = $ height$(u_j)$ and let $q = $height$(\nu)$. For the remainder of the argument we make no use of the fact that $i > j$. We may therefore argue symmetrically on $i$ and $j$.

Denote $f_{\#}^m(u_j)$ by $C_m$ for $m \geq 0$ and define $C_m$ for $m < 0$ as we did for $B_l$. We refer to the $B_l$'s as '$B$-blocks' and the $C_m$'s as '$C$-blocks'. Now $f_{\#}^k(\sigma)$ is obtained from $(E_i B_0 \cdots B_{k-1})f_{\#}^k(\nu)(\bar{C}_{k-1} \cdots \bar{C}_0 \cdot \bar{E}_j)$ by cancellation at the two indicated junctures. Suppose that as part of this cancellation process a subpath of $B_0 \cdots B_{k-1}$ that contains at least three $B$-blocks cancels with a subpath of $\bar{C}_{k-1} \cdots \bar{C}_0$ that contains at least three $C$-blocks. Sublemma 1 implies that $f_{\#}^k(\sigma)$, and hence $\sigma$, is exceptional. In that case $M(\sigma) = 0$ and there is nothing to prove. We may therefore assume that no such cancellation occurs.

If $q > \max\{s, t\}$, then splittings of $\nu$ produced by Lemma 4.1.4 extend to splittings of $\sigma$ in contradiction to our assumption that $\sigma$ does not split at any vertex. After interchanging the roles of $i$ and $j$ if necessary, we may assume that $s \geq \max\{q, t\}$ and that any subpath of $B_0 \cdots B_{k-1}$ that cancels with a subpath of $\bar{C}_{k-1} \cdots \bar{C}_0$ contains fewer than three $B$-blocks.

If $q < s$ then all cancellation of $s$-edges in $B_0 \cdots B_{k-1}$ must be with $s$-edges in $\bar{C}_{k-1} \cdots \bar{C}_0$. It follows that $f_\#^4(\sigma) = [E_i B_0 B_1 B_2 B_3 f^4(\nu \bar{E}_j)]$ splits as $E_i \cdot [B_0 B_1 B_2 B_3 f^4(\nu \bar{E}_j)]$ and the proof concludes as in the $j < s$ case above.

If $q = s$, then we argue as in the $j = s$ case above. Restrict to $k$ so large that $B_{-k} \ldots B_{-1}$ contains more $s$-edges than $\nu$ does. Decompose $\nu$ into subpaths $\nu = \nu_1 \nu_2$ where $\nu_1$ is the shortest initial segment that contains each $s$-edge of $\nu$ that is canceled when $B_{-k} \ldots B_{-1}\nu$ is tightened. Then $\nu_1 = \bar{B}_{-1}\bar{B}_{-2} \ldots \bar{B}_{-(r-1)}\bar{B}_{-r}^*$ where $\bar{B}_{-r}^*$ is an initial segment of $\bar{B}_{-r}$ and $r \geq 0$. Now $\sigma$ can be obtained from $(E_i \bar{B}_{-1} \ldots \bar{B}_{-r-5})[B_{-r-5} \ldots B_{-1}\nu]\bar{E}_j$ by tightening at the indicated junctures. Applying $f_\#^k$ with $k \geq r + 6$, we see that $f_\#^k(\sigma)$ is obtained from $(E_i \cdot B_0 \cdot B_1 \cdot \ldots \cdot B_{k-r-6}) f_\#^k([B_{-r-5} \ldots B_{-1}\nu])$ $(\bar{C}_{k-1} \cdot \ldots \cdot \bar{C}_0 \bar{E}_j)$ by cancelling at the indicated junctures. By construction, $f_\#^k([B_{-r-5} \ldots B_{-1}\nu])(\bar{C}_{k-1} \cdot \ldots \cdot \bar{C}_0 \bar{E}_j)$ tightens to a path that begins with $B_{k-r-5}$. It follows that $f_\#^{r+6}(\sigma) = E_i \cdot B_0 \cdot f_\#^{r+6}([B_{-r-5} \ldots B_{-1}\nu \bar{E}_j])$ and hence that $M(\sigma) \leq \max\{r + 6, M([B_{-r-5} \ldots B_{-1}\nu \bar{E}_j])\}$. The proof concludes as in the $j = s$ case above.                                                         $\square$

We will need the following technical results in [BFH2].

LEMMA 5.7.9. *Suppose that* $f : G \to G$ *is as in Theorem* 5.1.8. *There is a constant* $C_1$ *so that if* $\omega$ *is a closed path that is not a Nielsen path,* $\sigma = \alpha \omega^k \beta$ *is a path and* $n > 0$, *then at most* $C_1$ *copies of* $f_\#^n(\omega)$ *are canceled when* $f_\#^n(\alpha) f_\#^n(\omega^k) f_\#^n(\beta)$ *is tightened to* $f_\#^n(\sigma)$.

*Proof.* It is sufficient to consider the case where $\beta$ is empty. The proof is by induction on the height of $\alpha$. If height$(\alpha) \leq$ height$(\omega)$, then none of the highest edges in $f_\#^n(\omega)$ are canceled during the tightening. We may therefore assume that height$(\alpha) >$ height$(\omega)$ and that we have verified the lemma for all $\alpha$ of lower height. By Lemma 4.1.4 it is sufficient to consider the case that $\alpha = E_j \alpha'$ where height$(\alpha') < j$. Let $B_l = f_\#^l(u_j)$ and $R_j = B_0 \cdot B_1 \cdot \ldots \cdot B_l \cdot \ldots$; the $B_l$'s are called the blocks of $R_j$. Then $f_\#^n(\sigma)$ is obtained from $(E_j B_0 \ldots B_{n-1})[f_\#^n(\alpha') f_\#^n(\omega^k)]$ by tightening at the indicated juncture. By the inductive hypothesis, the number of copies of $f_\#^n(\omega)$ that are canceled when $f_\#^n(\alpha') f_\#^n(\omega^k)$ is tightened to $[f_\#^n(\alpha') f_\#^n(\omega^k)]$ is bounded independently of $n, k, \alpha'$ and $\omega$. Sublemma 5.7.11 below therefore completes the proof.                            $\square$

SUBLEMMA 5.7.10. *If* $\gamma^N \subset R_j$ *and* $\gamma^N$ *contains at least three blocks of* $R_j$ *for some* $N \geq 3$, *then* $\gamma$ *is a Nielsen path.*

*Proof.* Lifting to the universal cover, we see that there exists $\tilde{R}_j = \tilde{B}_0 \tilde{B}_1 \ldots$ where $\tilde{f}(\tilde{B}_l) = \tilde{B}_{l+1}$ and there exists a covering translation $T$ corresponding to $\gamma$ so that the axis of $T$ has an interval $I$ in common with $\tilde{R}_j$ that contains three blocks of $\tilde{R}_j$ and three fundamental domains of the axis $A(T)$ of $T$. Since $T$ and $\tilde{f}$ both 'translate' the initial segments of the highest edges in $I$ (cf. the proof of Sublemma 1) the two lifts $T\tilde{f}$ and $\tilde{f}T$ of $f$ agree at a point and so must be equal. We conclude that $T$ commutes with $\tilde{f}$ and hence that $\gamma$ is Nielsen. ☐

SUBLEMMA 5.7.11. *There is a constant $K$ with the following property*: *If* $\gamma^K \subset R_j$, *then $\gamma$ is a Nielsen path.*

*Proof.* The proof is by induction on $j$, the $j = 1$ case being trivial.

If $u_j$ has $m$ edges, then $B_l$ can be written as a concatenation of at most $m$ edges and at most $m$ subpaths of rays $R_i$ with $i < j$. By the inductive hypothesis, there is a constant $K_0$ so that if $\gamma^{K_0}$ is contained in some $B_l$ then $\gamma$ is a Nielsen path.

Let $K = 5K_0$. If $\gamma^K \subset R_j$ then either $\gamma^{K_0}$ is contained in some $B_l$ or $\gamma^K$ contains three blocks of $R_j$. In either case $\gamma$ is a Nielsen path. ☐

# 6. The weak attraction theorem

*We assume throughout this section that $\mathcal{O}_{\#}$ fixes each element of $\mathcal{L}(\mathcal{O})$.*

Each exponentially-growing stratum $H_r$ of an improved relative train track map $f : G \to G$ has a canonically associated (see Lemma 4.2.5) finite set $P_r$ of paths in $G_r$. If $P_r$ is nonempty, then it contains a preferred element $\rho_r$. We assign a path $\hat{\rho}_r$ to each exponentially-growing stratum $H_r$ as follows. If $P_r \neq \emptyset$, then $\hat{\rho}_r = \rho_r$. If $P_r = \emptyset$, then choose a vertex in $H_r$ and define $\hat{\rho}_r$ to be the trivial path at that vertex. If $Z$ is a subgraph of $G$ and $\hat{\rho}_r$ is trivial, then a nontrivial path in $G$ is contained in $\langle Z, \hat{\rho}_r \rangle$ (see Definition 5.1.6) if and only if it is contained in $Z$.

An element of $\mathcal{L}(\mathcal{O})$ is said to be *topmost* if it is not contained in any other element of $\mathcal{L}(\mathcal{O})$. Let $\Lambda^-$ be the expanding lamination for $\mathcal{O}^{-1}$ that is paired (Lemma 3.2.4) with $\Lambda^+$. The following weak attraction theorem is the main result of this section. It is an explicit description of the basin of attraction of $\Lambda^+$ in the birecurrent elements of $\mathcal{B}$. In the next section we will exploit the fact that 'most' birecurrent paths are attracted to $\Lambda^+$.

THEOREM 6.0.1. *Suppose that $\Lambda^+$ is a topmost element of $\mathcal{L}(\mathcal{O})$, that $f : G \to G$ is an improved relative train track map representing $\mathcal{O}$ and that $H_r$ is the exponentially-growing stratum that determines $\Lambda^+$. Then there exists a subgraph $Z$ such that $Z \cap G_r = G_{r-1}$ and such that every birecurrent path $\gamma \subset G$ satisfies exactly one of the following.*

1. $\gamma$ is a generic line for $\Lambda^-$.

2. $\gamma \in \langle Z, \hat{\rho}_r \rangle$.

3. $\gamma$ is weakly attracted to $\Lambda^+$.

*Remark* 6.0.2. Suppose that $G = G_r$ and hence that $Z = G_{r-1}$. If $H_r$ is not a geometric stratum, then Lemma 5.1.7 and (eg-ii) imply that $\langle Z, \hat{\rho}_r \rangle$ contains the same bi-infinite paths as $G_{r-1}$. If $H_r$ is a geometric stratum, then (eg-iii) implies that $\hat{\rho}_r = \rho_r$ and that the set of bi-infinite paths in $\langle Z, \rho_r \rangle$ is the union of the bi-infinite paths in $G_{r-1}$ with the circuit $\rho_r$.

*Remark* 6.0.3. Suppose that $\phi : S \to S$ is a homeomorphism of a compact surface in Thurston normal form and that $S_0 \subset S$ is a proper subsurface that is a pseudo-Anosov component of $\phi$. The geometric analog of Theorem 6.0.1 implies that the expanding measured foliation $F$ for $\phi | S_0$ weakly attracts every closed curve that is not entirely contained in $S \backslash S_0$. If we work in the projective foliation space PF rather than in $\mathcal{B}$, then $F$ is only certain to attract every simple closed curve that is entirely contained in $S_0$. Thus the basin of attraction in $\mathcal{B}$ can be larger than the basin of attraction in PF.

The following proposition is one of the two main steps in proving Theorem 6.0.1.

PROPOSITION 6.0.4. *Suppose that $f : G \to G$ is an improved relative train track map representing $\mathcal{O}$, that $\Lambda^+$ is a topmost element of $\mathcal{L}(\mathcal{O})$ and that $H_r$ is the exponentially-growing stratum that determines $\Lambda^+$. Then there is a subgraph $Z$ of $G$ such that*:

1. $Z \cap G_r = G_{r-1}$.

2. *$Z$ contains every zero stratum and every exponentially-growing stratum other than $H_r$.*

3. *$f(E) \in \langle Z, \hat{\rho}_r \rangle$ for each edge $E$ in $Z$.*

4. *Suppose that $\sigma \subset G$ is a finite path whose endpoints are fixed by $f$. Then $\sigma \in \langle Z, \hat{\rho}_r \rangle$ if and only if $\sigma$ is not weakly attracted to $\Lambda^+$.*

We will state and prove some preliminary results before beginning the proof of Proposition 6.0.4.

LEMMA 6.0.5. *Suppose that $f : G \to G$ is an improved relative train track map, that $H_r$ is an exponentially-growing stratum, and that $X$ is a subgraph of $G$ that does not contain any edges in $H_r$. Then the set of bi-infinite paths in $\langle X, \hat{\rho}_r \rangle$ determines a closed subset of $\mathcal{B}_G$.*

*Proof.* If $\hat{\rho}_r$ is trivial then $\langle X, \hat{\rho}_r \rangle$ contains the same bi-infinite paths as $X$ and the lemma is obvious. We may therefore assume that $\hat{\rho}_r = \rho_r$.

A (not necessarily finite) path $\sigma \in \langle X, \rho_r \rangle$ has a *locally defined* canonical decomposition into subpaths in $\langle X, \rho_r \rangle$. More precisely, suppose that $\sigma = \ldots a_i a_{i+1} a_{i+2} \ldots$ is the decomposition into single edges of $G$. The $a_i$'s can be grouped into subpaths $b_j$ that are either single edges in $X$ or $\rho_r$ or $\bar{\rho}_r$. The $b_j$'s are uniquely determined by the following rule. If $(\bar{a}_i, a_{i+1})$ is an illegal turn in $H_r$, then $(\bar{a}_i, a_{i+1})$ is the illegal turn in either $\rho_r$ or $\bar{\rho}_r$, and $a_i, a_{i+1}$ and some adjacent edges are grouped into a $b_j$ that is $\rho_r$ or $\bar{\rho}_r$ respectively. All edges not so grouped are in $X$ and determine $b_j$'s. Let $M$ be the number of edges in $\rho_r$. Then the endpoints of the $b_j$'s are separated by at most $M$ edges and the subpath $a_{i-M} \cdots a_{i+M}$ determines whether or not the endpoint shared by $a_i$ and $a_{i+1}$ is an endpoint of some $b_j$. We say that the endpoints of the $b_j$'s are *cutting vertices* for $\sigma$. Any subpath of $\sigma$ that is bounded by cutting vertices is contained in $\langle X, \rho_r \rangle$.

Suppose that $\mu_i \to \gamma$ in $\mathcal{B}_G$ and that $\mu_i \in \langle X, \rho_r \rangle$. Write $\gamma$ as an increasing union of finite subpaths $\gamma_i$. After passing to a subsequence, we may assume that $\gamma_i$ is a subpath of $\mu_i$. There is a common subpath of $\mu_i$ and $\gamma_i$ that is bounded by cutting vertices in $\mu_i$ and that covers all of $\gamma_i$ with the exception of at most $M$ edges at the beginning and end. Thus $\gamma$ can be written as an increasing union of subpaths in $\langle X, \rho_r \rangle$. The canonical decompositions of these subpaths agree on their overlap so that $\gamma \in \langle X, \rho_r \rangle$. $\square$

The following lemma is due to Peter Scott.

LEMMA 6.0.6. *If $\Phi : F \to F$ is an automorphism of a finitely generated free group and $H$ is a finitely generated subgroup of $F$ such that $\Phi(H) \subset H$, then $\Phi(H) = H$. In particular $\Phi|H$ is an automorphism.*

*Proof* (Peter Scott). Suppose at first that $H$ is a free factor of $F$. Then $\Phi(H)$ is a free factor of $F$ and the Kurosh subgroup theorem (see also Lemma 2.6.2) implies that $\Phi(H)$ is a free factor of $H$. Since $\Phi(H)$ and $H$ have the same rank, $\Phi(H) = H$.

For general $H$, the LERF property implies that $H$ is a free factor of some finite index subgroup $F'$ of $F$. The subgroup $F'' = \cap_i \Phi^i(F')$ is an intersection of subgroups of fixed finite index and so itself has finite index in $F$. By the Kurosh subgroup theorem, $H \cap F''$ is a free factor of $F''$. By construction $\Phi(F'') = F''$. Thus $\Phi(H \cap F'') \subset H \cap F''$ and, by the special case considered above, $\Phi(H \cap F'') = H \cap F''$. Since this subgroup has the same finite index in both $\Phi(H)$ and $H$, the index of $\Phi(H)$ in $H$ must be one. In other words, $\Phi(H) = H$. $\square$

COROLLARY 6.0.7.    *Suppose that $f : G \to G$ is an improved relative train track map and that $X$ is a subgraph of $G$ that does not contain any edges in $H_r$. If $f(E) \in \langle X, \hat{\rho}_r \rangle$ for each edge $E$ of $X$, then $f_\#$ restricts to a bijection on the set of bi-infinite paths in $\langle X, \hat{\rho}_r \rangle$ and to a bijection on the set of finite paths in $\langle X, \hat{\rho}_r \rangle$ whose endpoints are fixed by $f$.*

*Proof.* Suppose that $\tau$ is a path in $\langle X, \hat{\rho}_r \rangle$ with fixed endpoints and that $x$ is one of the endpoints of $\tau$. Let $H$ be the subgroup of $\pi_1(G, x)$ consisting of those elements represented by closed paths in $\langle X, \hat{\rho}_r \rangle$ with basepoint at $x$ and let $\Phi : \pi_1(G, x) \to \pi_1(G, x)$ be the automorphism determined by $f$. Our hypotheses imply that $\Phi_\#(H) \subset H$, and so Lemma 6.0.6 implies that $\Phi_\#|H$ is an automorphism. Let $\beta$ be the element of $H$ determined by $\tau f(\bar{\tau})$ and let $\sigma$ be the closed path in $\langle X, \hat{\rho}_r \rangle$ with basepoint at $x$ whose corresponding element $\alpha \in H$ satisfies $\Phi(\alpha) = \beta$. Then $f_\#(\sigma\tau) = \tau$ and we have shown that $f_\#$ restricts to a surjection on the set of finite paths in $\langle X, \hat{\rho}_r \rangle$ whose endpoints are fixed by $f$. Injectivity is an immediate consequence of the fact that $f$ is a homotopy equivalence.

Since $f_\#$ induces a bijection on the set of bi-infinite paths in $X$, we may assume for the rest of the proof that $\hat{\rho}_r = \rho_r$. Noncontractible components of $X$ that do not contain either endpoint of $\rho_r$ are permuted by $f$. The restriction of $f$ to the union of these components is a homotopy equivalence and so induces a bijection on the set of bi-infinite paths that they carry. We may therefore assume that each component of $X$ contains an endpoint of $\rho_r$. Let $x$ be the initial endpoint of $\rho_r$ and let $H$ and $\Phi$ be defined as above.

A circuit in $G$ can be written as a concatenation of subpaths in $\langle X, \rho_r \rangle$ if and only if the conjugacy class that it determines in $\pi_1(G, x)$ contains an element of $H$. Since $\Phi|H$ is a homotopy equivalence, $f_\#$ induces a bijection of these circuits and hence on the set $\mathcal{C}$ of periodic bi-infinite paths in $\langle X, \rho_r \rangle$ that they determine. Every bi-infinite path in $\langle X, \rho_r \rangle$ can be approximated by an element of $\mathcal{C}$. Lemma 6.0.5 therefore implies that the set of bi-infinite paths in $\langle X, \rho_r \rangle$ is the closure of $\mathcal{C}$ and hence that $f_\#$ restricts to a surjection on the set of bi-infinite paths in $\langle X, \rho_r \rangle$ . Injectivity of the restriction follows from the fact that $f_\#$ induces a bijection on the set of bi-infinite paths in $G$.    □

*Proof of Proposition* 6.0.4.    We build $Z_i = Z \cap G_i$ in stages, beginning with $Z_r = G_{r-1}$.

*Step* 1 (on verifying (4)).    Before beginning the construction of $Z$, we make two general observations about proving (4) for $Z_i$ once (3) for $Z_i$ is known.

If $\sigma \in \langle Z_i, \hat{\rho}_r \rangle$, then $f_\#^k(\sigma) \in \langle Z_i, \hat{\rho}_r \rangle$ for all $k \geq 0$. Since $\langle Z_i, \hat{\rho}_r \rangle$ does not contain any $r$-legal paths with more $H_r$-edges than are contained in $\hat{\rho}_r$, $\sigma$ is not weakly attracted to $\Lambda^+$. Thus we need only prove the 'if' half of (4). The second observation is that we can replace $\sigma$ by $f_\#^k(\sigma)$ without loss of generality. This replacement clearly has no effect on being weakly attracted to $\Lambda^+$, so we need only check that if $f_\#^k(\sigma) \in \langle Z_i, \hat{\rho}_r \rangle$ for some $k > 0$, then $\sigma \in \langle Z_i, \hat{\rho}_r \rangle$. This follows from Corollary 6.0.7 and the fact that $\sigma$ is the only path in $G$ whose endpoints are fixed by $f$ and that has $f_\#^k$-image equal to $f_\#^k(\sigma)$.  $\square$

*Step* 2 ($Z_r$). Condition (3) for $Z_r$ follows from the invariance of $G_{r-1}$. Suppose that $\sigma \subset G_r$ has fixed endpoints and is not weakly attracted to $\Lambda^+$. After $\sigma$ is replaced by some $f_\#^k(\sigma)$ if necessary, we may assume that each $f_\#^i(\sigma)$ contains the same number, say $m$, of illegal turns in $H_r$. If $m = 0$, then $\sigma$ is $r$-legal; since $\sigma$ is not weakly attracted to $\Lambda^+$, $\sigma \subset G_{r-1} = Z_{r-1}$. Suppose then that $m > 0$. Lemma 4.2.6 implies that $P_r \neq \emptyset$ and hence that $\hat{\rho}_r = \rho_r$. Moreover, with $k$ increased if necessary, Lemma 4.2.6 and the fact that $\rho_r$ is the unique periodic element of $P_r$ imply that $\sigma$ splits into pieces that are either $r$-legal or equal to either $\rho_r$ or $\bar{\rho}_r$. Since $\sigma$ is not weakly attracted to $\Lambda^+$, each of the $r$-legal pieces must lie in $G_{r-1} = Z_{r-1}$. Thus $\sigma \in \langle Z_{r-1}, \rho_r \rangle$ and (4) is satisfied.

Assume now that we have defined $Z_{i-1}$ satisfying (1)–(4) for some $i > r$. To complete the inductive step we will define $Z_i$ maintaining (2) and then verify (3) and (4).

*Step* 3 (the case that $H_i$ is non-exponentially growing). When $H_i$ is non-exponentially growing, $H_i$ is a single edge $E_i$ and $f(E_i) = E_i \cdot u_i$ for some closed path $u_i$ whose basepoint is fixed by $f$. If $u_i$ is weakly attracted to $\Lambda^+$, then define $Z_i = Z_{i-1}$. Condition (3) for $Z_{i-1}$ implies condition (3) for $Z_i$. To verify (4), we must show that if $\sigma \subset G_i$ has fixed endpoints and is not weakly attracted to $\Lambda^+$, then $\sigma \subset G_{i-1}$. Suppose to the contrary that $\sigma$ crosses $E_i$. Since the endpoints of $\sigma$ are fixed, they are not in the interior of $E_i$. Lemma 4.1.4 implies that $\sigma$ splits into subpaths, at least one of which is a basic path of height $i$. Since $u_i$ is weakly attracted to $\Lambda^+$, Lemma 3.1.16 implies that $u_i$ is not a Nielsen path. After $\sigma$ is replaced by an iterate, ne-(iii) implies that there is a further splitting of $\sigma$ into pieces, at least one of which is $E_i$ or $\bar{E}_i$. This contradicts Corollary 4.2.4 and so verifies (4).

If $u_i$ is not weakly attracted to $\Lambda^+$, then define $Z_i = Z_{i-1} \cup E_i$. The inductive hypothesis implies that $u_i \in \langle Z_{i-1}, \hat{\rho}_r \rangle$ so (3) is satisfied. Assume that $\sigma \subset G_i$ has fixed endpoints and is not weakly attracted to $\Lambda^+$. Since $\sigma$ splits into subpaths that are either entirely contained in $G_{i-1}$ or are basic pieces of height $i$, we may assume that $\sigma$ is a basic piece of height $i$. We may assume further that no iterate of $\sigma$ splits as a concatenation of two basic paths or as a concatenation of one basic path and a path in $G_{i-1}$. Condition ne-(iii) of Theorem 5.1.5 implies, after $\sigma$ is replaced by an iterate if necessary, that $\sigma$ is an exceptional path of height $i$; i.e. $\sigma = E_i \tau^k \bar{E}_j$ or $E_i \bar{\tau}^k \bar{E}_j$ where $j \leq i$, $\tau$ is a Nielsen path, $u_i = \tau^l$ and $u_j = \tau^m$. The inductive hypothesis implies that $\tau \in \langle Z_{i-1}, \hat{\rho}_r \rangle$. If $i \neq j$, then the inductive hypothesis also implies that $E_j \in \langle Z_{i-1}, \hat{\rho}_r \rangle$. Thus $\sigma \in \langle Z_i, \hat{\rho}_r \rangle$ and (4) is satisfied.

*Step* 4 (the case that $H_i$ is exponentially growing). Suppose that $H_i$ is exponentially growing. Define $Z_i = Z_{i-1} \cup H_i$. For each edge $E_j$ of $H_i$, $f(E_j) = a_1 \cdot b_1 \cdot a_2 \cdot \ldots \cdot a_m$ where each $a_l$ is a subpath in $H_i \subset Z_i$ and each $b_l \subset G_{i-1}$. Since $\Lambda^+$ is topmost, each $b_l$ is not weakly attracted to $\Lambda^+$. If the endpoints of $b_l$ are fixed, then the inductive hypothesis implies that $b_l \in \langle Z_{i-1}, \hat{\rho}_r \rangle$. If the endpoints of $b_l$ are not fixed, then $b_l$ is contained in a zero stratum and so is contained in $Z_{i-1}$ by the inductive hypothesis. This verifies (3).

Suppose that $\sigma \subset G_i$ has fixed endpoints and is not weakly attracted to $\Lambda^+$. As in Step 2, Lemma 4.2.6 implies, after $\sigma$ is replaced by an iterate if necessary, that either $P_i = \emptyset$ and $\sigma$ is $i$-legal or $\hat{\rho}_i = \rho_i$ and $\sigma$ splits into pieces that are either $i$-legal or equal to $\rho_i$ or $\bar{\rho}_i$. By Lemma 4.2.1, an $i$-legal path splits into subpaths $c_j \subset H_i \subset Z_i$ and $d_j \subset G_{i-1}$; as in the preceding argument the inductive hypothesis implies that $d_j \in \langle Z_{i-1}, \hat{\rho}_r \rangle$. It therefore suffices to assume that $\hat{\rho}_i = \rho_i$ and to show that $\rho_i \in \langle Z_i, \hat{\rho}_r \rangle$. Decompose $\rho_i = \alpha_i \beta_i$ where $\alpha_i$ and $\beta_i$ are $i$-legal. If the initial endpoint of $\rho_i$ is a vertex, let $E$ be the initial edge of $\rho_i$; otherwise, let $E$ be the edge that contains the initial endpoint of $\rho_i$. For sufficiently large $k$, $\alpha_i$ is a subpath of $f^k_\#(E)$. Since $f^k_\#(E) \in \langle Z_i, \hat{\rho}_r \rangle$ and the initial and terminal segments of $\alpha_i$ are in $H_i$, and so in particular are not in $G_r$, $\alpha_i \in \langle Z_i, \hat{\rho}_r \rangle$. A similar argument holds for $\beta_i$. Thus $\rho_i$ and $\bar{\rho}_i$ are contained in $\langle Z_i, \hat{\rho}_r \rangle$. This completes the proof in the case that $H_i$ is exponentially growing.

*Step* 5 (the case that $H_i$ is a zero stratum). If $H_i$ is a zero stratum, then define $Z_i = Z_{i-1} \cup H_i$. Property (4) for $G_i$ follows from property (4) for $G_{i-1}$ and the observation that a path in $G_i$ with fixed endpoints cannot be contained in, and hence cannot intersect, a component of $G_i$ that is mapped off of itself.

Theorem 5.1.5 implies that $H_{i+1}$ is an exponentially-growing stratum and that $H_i$ is a forest. For each edge $E_j$ of $H_{i+1}$, $f(E_j) = a_1 \cdot b_1 \cdot a_2 \cdot \ldots \cdot a_m$ where each $a_l$ is a subpath in $H_{i+1}$ and each $b_l \subset G_i$. Since $\Lambda^+$ is topmost, each $b_l$ is not weakly attracted to $\Lambda^+$. Let $\{\beta_l\}$ be the set of paths in $H_i$ that occur as a $b_k$ in the above decomposition for some $f(E_j)$.

Suppose that $P, Q \in H_{i+1} \cap H_i$. We write $P \sim Q$ if $P$ and $Q$ belong to the same component of $H_i$ or equivalently if there is a path in $H_i$ that connects $P$ and $Q$. There is also an equivalence relation generated by $P \sim' Q$ if and only if $P$ and $Q$ are the endpoints of some $\beta_l$. If we collapse each component of the forest $H_i$ to a point then the image of $G_{i+1}$ is a graph $\hat{G}_{i+1}$ that is homotopy equivalent to $G_{i+1}$. We can also view $\hat{G}_{i+1}$ as being obtained from $G_{i+1}\setminus \operatorname{int}(H_i)$ by identifying $\sim$ equivalent points in $H_{i+1}\cap H_i$. Let $\hat{G}'_{i+1}$ be the graph obtained from $G_{i+1}\setminus \operatorname{int}(H_i)$ by identifying $\sim'$ equivalent points in $H_{i+1}\cap H_i$. Then $\hat{G}_{i+1}$ is obtained from $\hat{G}'_{i+1}$ by identifying the elements in certain finite subsets. Each $\sim$ equivalence class that contains $m$ $\sim'$ equivalence classes determines a set of $m$ points in $\hat{G}'_{i+1}$ that must be identified to form $\hat{G}_{i+1}$. In particular, $\hat{G}'_{i+1}$ corresponds to a free factor system $\mathcal{F}$ of $\mathcal{F}(G_{i+1})$ that lies between $\mathcal{F}(G_i) = \mathcal{F}(G_{i-1})$ and $\mathcal{F}(G_{i+1})$. A bi-infinite path is carried by $\mathcal{F}$ if and only if it is contained in $\langle (G_{i+1}\setminus H_i), \cup\beta_l\rangle$. This collection of bi-infinite paths is mapped into itself by $f_\#$, so $\mathcal{F}$ is invariant under the action of an iterate of $\mathcal{O}$. Since $f : G \to G$ is reduced and $\mathcal{F}$ carries the expanding lamination determined by $H_{i+1}$, $\mathcal{F} = \mathcal{F}(G_{i+1})$. This implies that $\sim$ and $\sim'$ must be the same relation.

Condition z-(ii) of Theorem 5.1.5 and the fact that each $f(v)$ is fixed imply that $f(\beta_l)$ is a path in $G_{i-1}$ with fixed endpoints. The inductive hypothesis implies that $f(\beta_l) \in \langle Z_{i-1}, \hat{\rho}_r\rangle$.

If $P_r = \emptyset$, then $f(\beta_l) \subset Z_{i-1}$; since each edge in $H_i$ is crossed by some $\beta_l$, (3) is satisfied. We may therefore assume that $\hat{\rho}_r = \rho_r$ and that $\beta_l$ has a decomposition as an alternating concatenation of subpaths $\mu_j$ that map into $Z_{i-1}$ and $\nu_j$ that map to either $\rho_r$ or $\bar{\rho}_r$.

Suppose that $\delta_1$ and $\delta_2$ are paths in $H_i$ with endpoints in $H_{i+1}\cap H_i$ and that $\delta_1$ and $\delta_2$ have decompositions into $\mu_j$'s and $\nu_j$'s as above. We claim that if $\delta_1$ and $\delta_2$ have a common initial endpoint, then $[\delta_1^{-1}\delta_2]$ has a decomposition into $\mu_j$'s and $\nu_j$'s as above. It suffices to prove that the maximum common initial segment $\alpha$ of $\delta_1$ and $\delta_2$ contains every $\nu_j$ that it intersects. If this fails, then the image of the initial segment of both $\delta_1\setminus\alpha$ and $\delta_2\setminus\alpha$ would complete the partial crossing of $\rho_r$ or $\bar{\rho}_r$ begun in the image of $\alpha$. By condition eg-(i) of Theorem 5.1.5, $\rho_r$ and $\bar{\rho}_r$ have different initial edges so the partial image can only be completed in one way. Thus the initial segments of $\delta_1\setminus\alpha$ and $\delta_2\setminus\alpha$ have the same image, in contradiction to the fact (z-(ii) of Theorem 5.1.5) that $f|H_i$ is an immersion. This verifies our claim.

Since $\sim$ equals $\sim'$, every path $\delta \subset H_i$ with endpoints in $H_{i+1}\cap H_i$ can be expressed as $[b_1 b_2 \ldots b_m]$. The previous paragraph and induction imply that each $\delta_k = [b_1 \ldots b_k]$, and in particular $\delta$, has a decomposition into $\mu_j$'s and $\nu_j$'s and that each $\nu_j$ occurs in the decomposition of some $b_l$.

We next check that each $\nu_j$ is contained in a single edge of $H_i$. Suppose to the contrary that some $\nu_j$ crosses a vertex $w$. Condition z-(iii) of Theorem 5.1.5 implies that there is a (possibly trivial) path $\tau$ that starts at $w$, ends in $H_{i+1}\cap H_i$ and intersects $\nu_j$ only in $w$. Choose $\beta_l$ that contains $\nu_j$. The unique path $\delta \subset H_i$ that connects the initial endpoint $P$ of $\beta_l$ to the terminal endpoint $Q$ of $\tau$ agrees with $\beta_l$ up to $w$ and then follows $\tau$. But then $\delta$ deviates from $\beta_l$ in the middle of $\nu_j$, in contradiction to our observation in the previous paragraph that the maximum common initial subinterval of $\delta$ and $\beta_l$ contains each $\nu_j$ that it intersects. We conclude that each $\nu_j$ is contained in a single edge.

Given an edge $e$ of $H_i$, choose $\beta_l = \mu_1 \nu_1 \mu_2 \nu_2 \ldots$ that contains it. The endpoints of $e$ are not contained in the interior of any $\nu_j$ so $f(e) \in \langle Z_{i-1}, \rho_r \rangle$. This completes the proof of (3). $\qquad\square$

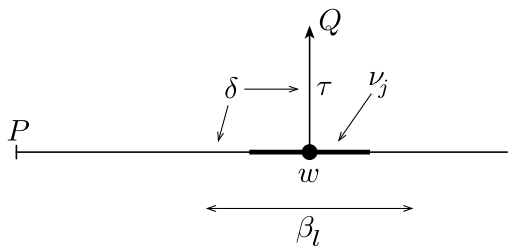

The following proposition is the second main piece of the proof of Theorem 6.0.1.

PROPOSITION 6.0.8. *Suppose that $H_s$ is an exponentially-growing stratum of an improved relative train track map $f : G \to G$, that $\gamma \subset G_s$ is a birecurrent path that is not contained in $G_{s-1}$ and that $\gamma$ is not weakly attracted to the expanding lamination $\Lambda_s^+$ associated to $H_s$. If $H_s$ is nongeometric, then $\gamma$ is a generic line for $\Lambda_s^-$. If $H_s$ is geometric, then either $\gamma$ is a generic line for $\Lambda_s^-$ or $\gamma = \rho_s$.*

*Proof of Proposition* 6.0.8 *in the geometric case.* We use the notation of Definition 5.1.4; in particular, $\phi : S \to S$ is a pseudo-Anosov homeomorphism, $Q$ is a graph, $\mathcal{A}$ is a collection of annuli, $Y = Q \cup \mathcal{A} \cup S$, $\Phi : (Y, Q) \to (G_s,$ noncontractible components of $G_{s-1})$ is a homotopy equivalence and $h : Y \to Y$ is a homotopy equivalence that satisfies $\Phi h \simeq f \Phi$. Let $\lambda$ be a generic line of $\Lambda_s^+$.

For every bi-infinite path $\sigma \subset G_s$ there is a bi-infinite path $\sigma^* \subset Y$ that intersects $\partial S$ transversely, that intersects each $A_i$ (if at all) in arcs that run from one component of $\partial A_i$ to the other and that satisfies $\Phi_\#(\gamma^*) = \gamma$. If $\sigma$ is birecurrent then either $\sigma^* \subset G_{s-1}$, $\sigma^* \subset S$ or its intersection with $\partial S$ decomposes $\sigma^*$ into an alternating concatenation of finite paths $a_i^* \subset Y \setminus \mathrm{int}(S)$ and finite geodesics $b_i^* \subset S$; each $a_i^*$ represents a nontrivial element in $\pi_1(Y \setminus \mathrm{int}(S))$ and each $b_i^*$ represents a nontrivial element in $\pi_1(S, \partial S)$.

We assume that $h_\#(\sigma^*) = ((f_\#(\sigma))^*$. In other words, we assume that $h_\#(\sigma^*)$ intersects $\partial S$ and the $A_i$'s as above.

For each $b_i^*$ and $k > 0$, $h_\#^k(\sigma^*)$ contains a path in $S$ that is homotopic rel $\partial S$ to $\phi^k(b_i^*)$. In particular, if $\sigma$ corresponds to a circuit (and so is periodic), then the length of the components of $h_\#^k(\sigma^*) \cap S$ tend to infinity as $k \to \infty$. Thus $\lambda^* \cap S$ cannot contain any finite components and we conclude that $\lambda^* \subset S$.

We say that two paths $\mu_1$ and $\mu_2$ in $S$ are $\delta$-parallel if there is a map of a rectangle into $S$ that agrees with $\mu_1$ on its upper horizontal edge, agrees with $\mu_2$ on its lower horizontal edge and maps vertical fibers into arcs with length less than $\delta$. We say that $\sigma^*$ is weakly attracted to $\lambda^*$ if for each finite subpath $\lambda_0^*$ of $\lambda^*$ and each $\delta > 0$ there exists $k > 0$ so that $h_\#^k(\sigma^*)$ contains a subpath that is $\delta$-parallel to $\lambda_0^*$. It is easy to check that $\sigma^*$ is weakly attracted to $\lambda^*$ if and only if $\sigma$ is weakly attracted to $\lambda$ and we leave this to the reader. Thus generic lines for $\Lambda_s^+$ correspond to generic leaves of the expanding lamination for $\phi$. The lemma now follows from well-known properties of the pseudo-Anosov map $\phi$.                                                                    $\square$

The following lemma is a refinement of Lemma 4.2.5 and is needed for the proof of Proposition 6.0.8 in the nongeometric case.

LEMMA 6.0.9. *Suppose that $H_s$ is an exponentially-growing stratum of an improved relative train track map $f : G \to G$.*

- *If $P_s = \emptyset$, then there is a positive integer $K_1$ with the following property. If $k \geq K_1$ and $\sigma \subset G_s$ is a path with exactly one illegal turn in $H_s$, then $f_\#^k(\sigma)$ is $s$-legal.*

- *If $P_s \neq \emptyset$ and $\rho_s'$ is a noninitial, nonterminal subpath of $\rho_s$, then there exist positive integers $L < K_1$ with the following property. If $k \geq K_1$ and if both $\sigma \subset G_s$ and $f_\#^k(\sigma)$ have exactly one illegal turn in $H_s$, then $f_\#^L(\sigma)$ contains $\rho_s'$ as an unoriented subpath.*

*Proof.* After subdividing if necessary, we may assume that each $\rho \in P_s$ has endpoints at vertices.

Suppose that for each $k \geq 1$ there are paths $\sigma_k \subset G_s$ such that both $\sigma_k$ and $f_\#^k(\sigma_k)$ have exactly one illegal turn in $H_s$. Decompose $\sigma_k$ as a concatenation $\sigma_k = \alpha_k \beta_k$ of $s$-legal subpaths. If the edge lengths of the $\alpha_k$'s are unbounded, then after passing to a subsequence, we may assume that the $\bar\alpha_k$'s converge to an infinite ray $\bar\alpha^*$, in the sense that the length of the maximal common initial segment of $\bar\alpha_k$ and $\bar\alpha^*$ goes to infinity as $k$ goes to infinity. If the length of the $\alpha_k$'s is bounded, then after passing to a subsequence, we may assume that each $\bar\alpha_k \subset \bar\alpha_{k+1}$; in this case, let $\bar\alpha^* = \cup_{k=1}^\infty \bar\alpha_k$. Define $\beta^*$ similarly and let $\rho = \bar\alpha^* \beta^*$. Then $f_\#^k(\rho)$ has one illegal turn in $H_s$ for all $k$ so Lemma 4.2.6 implies that $\rho$ contains an element of $P_s$.

The first item follows immediately. The second item follows from this construction and the fact that if $L$ is the cardinality of $P_s$, then $f_\#^L(\rho) = \rho_s$ for each $\rho \in P_s$.                                                                    $\square$

*Proof of Proposition* 6.0.8 *in the nongeometric case.* There is no loss in assuming that $G = G_s$. After passing to an iterate if necessary, we may assume that there is an improved relative train track map and filtration representing $\mathcal{O}^{-1}$ such that $\mathcal{F}(G_{s-1})$ is realized by a filtration element. If $P_s \neq \emptyset$, then we may assume that the endpoints of $\rho_s$ are vertices.

The first step in the proof is to show that $\gamma$ loses illegal turns in $H_s$ at an exponential rate under the action of the $f_\#$. This will be made explicit during the course of the proof.

If $P_s = \emptyset$, let $K_1$ be the constant of the first item of Lemma 6.0.9 and let $L = 0$. If $P_s \neq \emptyset$, let $\rho'_s$ be a subpath of $\rho_s$ that contains all but a proper initial segment of the first edge of $\rho_s$ and a proper terminal segment of the last edge of $\rho_s$ and let $L$ and $K_1$ be the constants of the second item of Lemma 6.0.9. By Lemma 4.2.2 and Corollary 4.2.4, there exists $K \geq K_1$ so that for any edge $E$ in $H_s$, every path in $G$ that contains $f_\#^{K-L}(E)$ as a subpath is weakly attracted to $\Lambda^+$.

Suppose that $\gamma_0$ is a finite subpath of $\gamma$ such that $f_\#^K(\gamma_0)$ is a subpath of $f_\#^K(\gamma)$. We claim that if $\gamma_0$ contains $m$ illegal turns in $H_s$, then $f_\#^K(\gamma_0)$ contains at most $m - \lceil \frac{m}{3} \rceil$ illegal turns in $H_s$. Let $l$ be the number of illegal turns that $f_\#^K(\gamma_0)$ has in $H_s$. Write $\gamma_0$ as a concatenation of subpaths $\gamma_0 = \alpha_1 \alpha_2 \cdots \alpha_l$ where $f_\#^K(\gamma_0) = f_\#^K(\alpha_1) f_\#^K(\alpha_2) \dots f_\#^K(\alpha_l)$ and where each $f_\#^K(\alpha_i)$ contains one illegal turn in $H_s$. Each $\alpha_i$ must contain at least one illegal turn in $H_s$ since the image of an $s$-legal path is $s$-legal. It therefore suffices to show that $\alpha_{i-1} \alpha_i \alpha_{i+1}$ contains at least four illegal turns in $H_s$ for each $2 \leq i \leq l-1$. By Lemma 6.0.9, we may assume that $P_s \neq \emptyset$ and hence that $\hat{\rho}_s = \rho_s$.

Let $\beta_i = f_\#^L(\alpha_i)$. If $\alpha_{i-1} \alpha_i \alpha_{i+1}$ contains exactly three illegal turns in $H_s$, then $\beta_{i-1} \beta_i \beta_{i+1}$ contains exactly three illegal turns in $H_s$. Lemma 6.0.9 implies that, with the possible exception of short initial and terminal segments, $\beta_{i-1} \beta_i \beta_{i+1}$ contains a subpath of the form $\rho^1 \sigma_1 \rho^2 \sigma_2 \rho^3$, where each $\rho^i$ is $\rho_s$ or $\bar{\rho}_s$. Moreover, all of the $H_s$ edges that are canceled when $f^{K-L}(\beta_{i-1} \beta_i \beta_{i+1})$ is tightened to $f_\#^K(\alpha_{i-1} \alpha_i \alpha_{i+1})$ are contained in $\rho^1, \rho^2$ and $\rho^3$. Lemma 5.1.7 and eg-(ii) imply that the endpoints of $\rho_s$ are distinct and not both contained in noncontractible components of $G_{s-1}$. Since these endpoints are fixed points, they cannot be contained in contractible components of $G_{s-1}$. It follows that at least one of $\sigma_1$ or $\sigma_2$ must contain an edge $E$ of $H_s$. But then $f_\#^K(\alpha_{i-1} \alpha_i \alpha_{i+1})$ and hence $f_\#^K(\gamma)$ contains $f_\#^{K-L}(E)$ and so is weakly attracted to $\Lambda_s^+$. This contradiction verifies our claim and completes the first step in the proof.

The second step in the proof is to show that for any finite subpath $\gamma_1$ of $\gamma$ there exists $\tau \subset G$ with a uniformly bounded (i.e. bounded independently of $\gamma$ and $\gamma_1$) number of $H_s$ edges such that $\gamma_1$ is a subpath of $\mathcal{O}_\#^{-k}(\tau)$ in $G$ for some $k \geq 0$. If $G_{s-1} = \emptyset$, then $\tau$ will be a circuit; if $G_{s-1} \neq \emptyset$, then $\tau$ will be a bi-infinite path with both ends in $G_{s-1}$.

After extending $\gamma_1$ to a larger subpath of $\gamma$ if necessary, we may assume that $f_\#^K(\gamma_1)$ is a subpath of $\gamma' = f_\#^K(\gamma)$. For future reference, note that if $C'$ is a positive integer so that no path with edge length greater than $C'$ has trivial $f_{\#}^K$-image, then at most $C'$ initial and $C'$ terminal edges need to be added to $\gamma_1$ to arrange this property. Let $C$ be the bounded cancellation constant for $f^K$ and let $\gamma_2'$ be the subpath of $\gamma'$ that is obtained from $f_{\#}^K(\gamma_1)$ by adding $2C$ initial edges and $2C$ terminal edges. We claim that $f_{\#}^K(N(\gamma_1)) \supset N(\gamma_2')$.

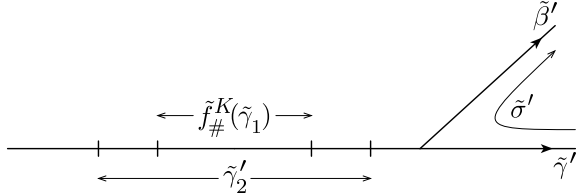

To see this, it is convenient to work in the universal cover $\Gamma$. Choose lifts $\tilde{f}^K : \Gamma \to \Gamma$, $\tilde{\gamma}_1 \subset \tilde{\gamma}$ and $\tilde{f}_{\#}^K(\tilde{\gamma}_1) \subset \tilde{\gamma}_2' \subset \tilde{\gamma}' = \tilde{f}_{\#}^K(\tilde{\gamma})$. Given $\beta' \in N(\gamma_2')$, choose a lift $\tilde{\beta}'$ that contains $\tilde{\gamma}_2'$. There is a unique bi-infinite path $\tilde{\beta} \subset \Gamma$ such that $\tilde{f}_{\#}^K(\tilde{\beta}) = \tilde{\beta}'$. Let $\tilde{\sigma} \subset \Gamma$ be the bi-infinite path connecting the forward end of $\tilde{\gamma}$ to the forward end of $\tilde{\beta}$. Then $\tilde{\sigma}' = \tilde{f}_{\#}^K(\tilde{\sigma})$ is the bi-infinite path connecting the forward end of $\tilde{\gamma}'$ to the forward end of $\tilde{\beta}'$. In particular, $\tilde{\sigma}'$ is disjoint from $\tilde{\gamma}_2'$. The bounded cancellation lemma therefore implies that $\tilde{f}^K(\tilde{\sigma})$ is disjoint from $\tilde{f}^K(\tilde{\gamma}_1)$ and hence that $\tilde{\sigma}$ is disjoint from $\tilde{\gamma}_1$. A symmetric argument on the backward ends implies that $\tilde{\beta} \subset N(\gamma_1)$ as claimed.

After increasing the number of edges in $\gamma_2' \subset \gamma'$ by at most $2C'$, we may assume that $f_{\#}^K(\gamma_2')$ is a subpath of $f_{\#}^K(\gamma')$.

The difference between the number of illegal turns of $f_{\#}^K(\gamma_1)$ in $H_s$ and the number of illegal turns of $\gamma_2'$ in $H_s$ is at most $4C + 2C'$. If $\gamma_1$ contains $m$ illegal turns in $H_s$ and $m$ is sufficiently large, say $m > M$, then by the first step, $\gamma_2'$ contains fewer than $m$ illegal turns in $H_s$. Iterating this, we conclude that for any $\gamma_1$, there exist $k$ and a finite subpath $\hat{\gamma}_2 \subset \hat{\gamma} = f_{\#}^k(\gamma)$ such that $f_{\#}^k(N(\gamma_1)) \supset N(\hat{\gamma}_2)$ and such that $\hat{\gamma}_2$ contains at most $M$ illegal turns in $H_s$. Lemma 4.2.2 and Corollary 4.2.4 imply that the number of $H_s$-edges in $\hat{\gamma}_2$ is bounded independently of $\gamma$ and $\gamma_1$. If $G_{s-1} = \emptyset$, then extend $\hat{\gamma}_2$ to a circuit $\tau$ that crosses each edge at most one more time than $\hat{\gamma}_2$ does. If $G_{s-1} \neq \emptyset$, then extend $\hat{\gamma}_2$ to a bi-infinite path $\tau$ that crosses an $H_s$ edge at most two more times than $\hat{\gamma}_2$ does. Since $f_{\#}^k(N(\gamma_1)) = \mathcal{O}_{\#}^k(N(\gamma_1))$, we have shown that $\mathcal{O}_{\#}^{-k}(\tau) \in \mathcal{O}_{\#}^{-k}(N(\hat{\gamma}_2)) \subset N(\gamma_1)$. This completes the second step.

Let $f' : G' \to G'$ be an improved relative train track map representing $\mathcal{O}^{-1}$ such that $\mathcal{F}(G_{s-1})$ is realized by a filtration element. For the final step, we consider first the case that $G_{s-1} \neq \emptyset$. Since $\mathcal{F}(\Lambda^-) = \mathcal{F}(\Lambda^+)$ is not carried by $\mathcal{F}(G_{s-1})$ and since there are no $\mathcal{O}^{-1}$-invariant free factor systems between $\mathcal{F}(G_{s-1})$ and $\{[[F_n]]\}$, $\Lambda_s^-$ is associated to the highest stratum $H_{s'}'$ and $\mathcal{F}(G_{s'-1}') = \mathcal{F}(G_{s-1})$. Choose a homotopy equivalence $h : G \to G'$ that respects the markings and so induces the natural identification of $\mathcal{B}(G)$ with $\mathcal{B}(G')$. Since $\mathcal{F}(G'_{s'-1}) = \mathcal{F}(G_{s-1})$, the bounded cancellation lemma implies that the number of $H'_{s'}$ edges in $h_\#(\alpha)$ is bounded independently of the path $\alpha \subset G_{s-1}$. Let $\gamma' = h_\#(\gamma) \subset G'$. Given a finite subpath $\gamma'_1 \subset \gamma'$, choose a finite subpath $\gamma_1 \subset \gamma$ such that $h_\#(N(\gamma_1)) \subset N(\gamma'_1)$. By the second step, there exist $\tau$ and $k$ for $\gamma_1$ as above so that $\mathcal{O}_\#^{-k}(\tau) \subset N(\gamma_1)$. Let $\tau' = h_\#(\tau)$. Then $(f')^k_\#(\tau') = \mathcal{O}_\#^{-k} h_\#(\tau) = h_\# \mathcal{O}_\#^{-k}(\tau) \in N(\gamma'_1)$, or equivalently $\gamma'_1$ is a subpath of $(f')^k_\#(\tau')$. There is a positive integer $M'$, independent of $\gamma$ and $\gamma_1$, such that $\tau'$ contains fewer than $M'$ edges in $H'_{s'}$.

Let $\delta'_0 \subset \gamma'$ be any finite subpath that crosses an edge in $H'_{s'}$. Choose a finite subpath $\gamma'_1 \subset \gamma'$ that contains at least $M'+1$ copies of $\delta'_0$ and let $k$ and $\tau'$ be as in the preceding paragraph. At least one of the copies of $\delta'_0$ must be contained in $(f')^k(E')$ where $E'$ is a single edge of $H'_{s'}$ in $\tau'$. This implies that $\delta'_0$ is contained in every generic line of $\Lambda^-$. Since $\delta'_0$ was arbitrary, $\gamma'$ is a line of $\Lambda^-$. Since $\gamma'$ is birecurrent and is not contained in $G'_{s'-1}$, Lemma 3.1.15 implies that $\gamma'$ is a generic line of $\Lambda_s^-$. Since $\gamma \subset G$ and $\gamma' \subset G'$ determine the same line in $\mathcal{B}$, $\gamma$ is a generic line of $\Lambda_s^-$.

It remains to consider the case that $G_{s-1} = \emptyset$. As there are only finitely many possibilities for $\tau$, we may assume that $\tau$ is independent of $\gamma_1$. Choose $h : G \to G'$ that respects the markings, let $\gamma' = h_\#(\gamma)$, let $\tau' = h_\#(\tau)$ and let $M'$ be the number of edges in $\tau'$. For any finite subpath $\delta'_0 \subset \gamma'$, there is a finite subpath $\gamma'_1 \subset \gamma'$ that contains at least $2M'+1$ copies of $\delta'_0$.

If $\mathcal{O}_\#^{-k}(\tau')$ takes on only finitely many values, then $\gamma$ is the periodic bi-infinite path determined by an $\mathcal{O}$-invariant circuit. This contradicts the first step in the proof and we conclude that the number of edges in the circuit $\mathcal{O}_\#^{-k}(\tau')$ tends to $\infty$ as $k$ tends to $\infty$. It follows that for sufficiently large $k$, $\gamma'_1$ is a subpath of the bi-infinite path determined by $\mathcal{O}_\#^{-k}(\tau')$ that intersects at most two fundamental domains. In particular $\gamma'_1$ is contained in a subpath that is a concatenation of at most $2M'$ segments, each of which is a subset of $(f')^k(E')$ for a single edge $E'$ of $G'$. The proof now concludes as in the previous case.    $\square$

The following corollary will be strengthened at the end of the section after we complete the proof of Theorem 6.0.1. This partial result is used to prove Theorem 6.0.1.

COROLLARY 6.0.10. *Suppose that $f : G \to G$ is an improved relative train track map representing $\mathcal{O}$, that $\Lambda^+ \in \mathcal{L}(\mathcal{O})$ is associated to the exponentially-growing stratum $H_s$ and that $\Lambda^- \in \mathcal{L}(\mathcal{O}^{-1})$ is paired with $\Lambda^+$. If $\gamma \subset G_s$ is a bi-infinite path that is not a generic line of either $\Lambda^+$ or $\Lambda^-$, then $\gamma$ is weakly attracted to $\Lambda^+$ under the action of $\mathcal{O}$ if and only if $\gamma$ is weakly attracted to $\Lambda^-$ under the action of $\mathcal{O}^{-1}$.*

*Proof.* We may assume that $G = G_s$, and, after passing to an iterate if necessary, that there is an improved relative train track map $f' : G' \to G'$ for $\mathcal{O}^{-1}$ such that $\mathcal{F}(G_{s-1})$ is realized by a filtration element. Since $\mathcal{F}(\Lambda^-) = \mathcal{F}(\Lambda^+)$ is not carried by $\mathcal{F}(G_{s-1})$ and since there are no $\mathcal{O}^{-1}$-invariant free factor systems between $\mathcal{F}(G_{s-1})$ and $\{[[F_n]]\}$, $\Lambda^-$ is associated to the highest stratum $H'_{s'}$ and $\mathcal{F}(G'_{s'-1}) = \mathcal{F}(G_{s-1})$. If $\gamma$ is carried by $\mathcal{F}(G_{s-1}) = \mathcal{F}(G'_{s'-1})$ or if $\gamma$ is an $\mathcal{O}$-invariant circuit, then $\gamma$ is not weakly attracted to either $\Lambda^+$ or $\Lambda^-$. In all other cases, Proposition 6.0.8 implies that $\gamma$ is weakly attracted to both $\Lambda^+$ and $\Lambda^-$. $\qquad\square$

COROLLARY 6.0.11. *If $\Lambda^- \in \mathcal{L}(\mathcal{O}^{-1})$ is paired with a topmost lamination $\Lambda^+ \in \mathcal{L}(\mathcal{O})$, then $\Lambda^-$ is topmost.*

*Proof.* Choose an improved relative train track map $f : G \to G$ representing an iterate of $\mathcal{O}$, let $H_r$ be the exponentially-growing stratum that determines $\Lambda^+$ and let $Z$ be the subgraph of Proposition 6.0.4. Suppose that $\Lambda_s^+$ and $\Lambda_s^-$ are paired laminations associated to an exponentially-growing stratum $H_s$ with $s > r$. Conditions (2) and (3) of Proposition 6.0.4 imply that each generic line $\gamma$ of $\Lambda_s^+$ is contained in $\langle Z \cap G_s, \hat{\rho}_r \rangle$. Choose an $H_s$-edge $E$ that occurs infinitely often (counting orientation) in $\gamma$. A subpath of $\gamma$ that starts at an occurrence of $E$ and ends just before an occurrence of $E$ defines an $s$-legal circuit $\gamma_0$ that is contained in $\langle Z \cap G_s, \hat{\rho}_r \rangle$ (because the endpoints of $\gamma_0$ are cutting vertices as defined in the proof of Lemma 6.0.5). Lemma 4.2.1 and Corollary 4.2.4 imply that $\gamma_0$ is weakly attracted to $\Lambda_s^+$ under the action of $\mathcal{O}$. Lemma 3.1.16 and Corollary 6.0.10 therefore imply that $\gamma_0$ is weakly attracted to $\Lambda_s^-$ under the action of $\mathcal{O}^{-1}$. Corollary 6.0.7 implies that the set of bi-infinite paths in $\langle Z \cap G_s, \hat{\rho}_r \rangle$ is $\mathcal{O}^{-1}$-invariant. Lemma 6.0.5 therefore implies that each generic line of $\Lambda_s^-$ is contained in $\langle Z, \hat{\rho}_r \rangle$. On the other hand, the bi-infinite paths of $\langle Z \cap G_r, \hat{\rho}_r \rangle$ differ from the bi-infinite paths of $G_{r-1}$ in at most a circuit, so that $\langle Z, \hat{\rho}_r \rangle$ does not contain any generic lines of $\Lambda^-$. We conclude that a generic line of $\Lambda^-$ is not contained in the closure of a generic line of $\Lambda_s^-$ and hence that $\Lambda^-$ is topmost. $\qquad\square$

*Proof of Theorem* 6.0.1. Let $f : G \to G$, $\mathcal{O}$, $H_r$ and $Z$ be as in Proposition 6.0.4.

If a generic line of $\Lambda^-$ is weakly attracted to $\Lambda^+$, then $\Lambda^-$ (being closed and $\mathcal{O}$-invariant) would contain $\Lambda^+$. Lemma 3.1.15 implies that a bicurrent line of $\Lambda^-$ is either carried by a strictly smaller free factor system than $\mathcal{F}(\Lambda^-)$ or is generic. Since $\mathcal{F}(\Lambda^+) = \mathcal{F}(\Lambda^-)$, each generic line of $\Lambda^+$ would be a generic line of $\Lambda^-$ and so $\Lambda^+ = \Lambda^-$. This contradicts Proposition 3.3.3 and thereby shows that (1) and (3) are mutually exclusive. Since the bi-infinite paths of $\langle G_{r-1}, \hat{\rho}_r \rangle$ and the bi-infinite paths of $G_{r-1}$ differ by at most a peripheral curve, (1) and (2) are mutually exclusive. If $\gamma \in \langle Z, \hat{\rho}_r \rangle$, then $\mathcal{O}_\#^k(\gamma) \in \langle Z, \hat{\rho}_r \rangle$ for all $k > 0$. Since there is a uniform bound to the number of $H_r$-edges in an $r$-legal path in $\langle Z, \hat{\rho}_r \rangle$, $\gamma$ is not weakly attracted to $\Lambda^+$. We have now shown that (1), (2) and (3) are mutually exclusive.

Let $s$ be the smallest positive value for which $\gamma \subset G_s$. Since $G_{r-1} \subset Z$ we may assume that $s \geq r$. The $s = r$ case is proved by Proposition 6.0.8.

Suppose that $s > r$. If $\gamma$ splits into finite paths whose endpoints are fixed by $f$, then Proposition 6.0.4 completes the proof. We may therefore assume (Lemma 4.1.4) that $H_s$ is exponentially growing and that $\gamma \notin \langle G_{s-1}, \hat{\rho}_s \rangle$. Let $\Lambda_s^+$ and $\Lambda_s^-$ be the lamination pair associated to $H_s$. During the proof of Corollary 6.0.11 we showed that each generic line of $\Lambda_s^-$ is contained in $\langle Z, \hat{\rho}_r \rangle$. We may therefore assume that $\gamma$ is not a generic line of $\Lambda_s^-$ and hence, by Proposition 6.0.8, is weakly attracted to $\Lambda_s^+$. As in the proof of Proposition 6.0.4, Corollary 6.0.7 allows us to replace $\gamma$ by any $f_\#^k(\gamma)$. We may therefore assume that $\gamma$ contains subpaths of a generic line $\lambda$ of $\Lambda_s^+$ with arbitrarily many $H_s$-edges. Since $f$ maps the set of endpoints of edges in $H_s$ into itself, at least one of the vertices is the image of a vertex and is therefore a fixed point. Lemma 4.2.2 and birecurrence imply that $\gamma$ has a splitting into finite paths whose endpoints are fixed by $f$ and again Proposition 6.0.4 completes the proof. □

COROLLARY 6.0.12. *If $\Lambda^+ \in \mathcal{L}(\mathcal{O})$ and $\Lambda^- \in \mathcal{L}(\mathcal{O}^{-1})$ are paired topmost expanding laminations and $\gamma$ is a bi-recurrent path that is not a generic line of either $\Lambda^+$ or $\Lambda^-$, then $\gamma$ is weakly attracted to $\Lambda^+$ under the action of $\mathcal{O}$ if and only if $\gamma$ is weakly attracted to $\Lambda^-$ under the action of $\mathcal{O}^{-1}$.*

*Proof.* Choose an improved relative train track map $f : G \to G$ representing an iterate of $\mathcal{O}$, let $H_r$ be the exponentially-growing stratum that determines $\Lambda^+$ and let $Z$ be as in Proposition 6.0.1. If $\gamma$ is not weakly attracted to $\Lambda^+$ under the action of $\mathcal{O}$, then the same is true for all $\mathcal{O}_\#^{-k}(\gamma)$ and so $\mathcal{O}_\#^{-k}(\gamma) \in \langle Z, \hat{\rho}_r \rangle$ for all $k \geq 0$. Lemma 6.0.5 and the fact that generic leaves of $\Lambda^-$ are not contained in $\langle Z, \hat{\rho}_r \rangle$ imply that $\gamma$ is not weakly attracted to $\Lambda^-$ under the action of $\mathcal{O}^{-1}$. The symmetric argument with the roles of $\Lambda^+$ and $\Lambda^-$ reversed completes the proof. □

## 7. Reduction to UPG($F_n$)

In this section we reduce the Tits alternative for Out($F_n$) to Theorem 1.0.2. More precisely, we prove the following theorem.

THEOREM 7.0.1. *Suppose that $\mathcal{H}$ is a subgroup of Out($F_n$) that does not contain a free subgroup of rank 2. Then there are a finite index subgroup $\mathcal{H}_0$ of $\mathcal{H}$, a finitely generated free abelian group $A$, and a map $\Phi : \mathcal{H}_0 \to A$ such that Ker($\Phi$) is UPG.*

We begin by using the weak attraction theorem and Corollary 3.4.3 to analyze the stabilizers of topmost laminations.

PROPOSITION 7.0.2.   *Suppose that $\mathcal{H}$ is a subgroup of $\mathrm{Out}(F_n)$ and that $\Lambda^+ \in \mathcal{L}(\mathcal{O})$ and $\Lambda^- \in \mathcal{L}(\mathcal{O}^{-1})$ are paired topmost laminations for some $\mathcal{O} \in \mathcal{H}$. Then either $\mathcal{H}$ contains a free subgroup of rank two, or at least one of the subgroups, $\mathrm{Stab}_{\mathcal{H}}(\Lambda^+)$ or $\mathrm{Stab}_{\mathcal{H}}(\Lambda^-)$, has finite index.*

The proof of Proposition 7.0.2 reduces to the following technical lemma. Before proving the lemma, we use it to prove the proposition.

LEMMA 7.0.3.   *Suppose that $\mathcal{H}$ is a subgroup of $\mathrm{Out}(F_n)$, that $\Lambda^+ \in \mathcal{L}(\mathcal{O})$ and $\Lambda^- \in \mathcal{L}(\mathcal{O}^{-1})$ are paired topmost laminations for $\mathcal{O} \in \mathcal{H}$ and that $f : G \to G$ is an improved relative train track map representing $\mathcal{O}$ such that (see Definition 3.2.3) $\{F(\Lambda^+)\} = \{F(\Lambda^-)\} = \mathcal{F}(G_r)$ for some filtration element $G_r$. Let $Z$ and $\hat{\rho}_r$ be as in Theorem 6.0.1 and let $\lambda^{\pm}$ be generic lines for $\Lambda^{\pm}$. Then $\mathcal{H}$ has a finite index subgroup $\mathcal{H}_0$ such that $\psi(\lambda^+), \psi(\lambda^-) \notin \langle Z, \hat{\rho}_r \rangle$ for each $\psi \in \mathcal{H}_0$.*

*Proof of Proposition* 7.0.2. By Theorem 5.1.5, there is an improved relative train track map $f : G \to G$ and filtration representing an iterate of $\mathcal{O}$ such that $\{F(\Lambda^+)\} = \{F(\Lambda^-)\} = \mathcal{F}(G_r)$ for some filtration element $G_r$. Let $\mathcal{H}_0$ be the finite index subgroup of Lemma 7.0.3 and let $\psi \in \mathcal{H}_0$. Theorem 6.0.1 and Corollary 6.0.12 imply that one of the following conditions is satisfied:

(i)  $\psi(\lambda^+)$ is generic for $\Lambda^+$ or for $\Lambda^-$.

(ii)  $\psi(\lambda^+)$ is weakly attracted to $\Lambda^+$ under the action of $\mathcal{O}$ and to $\Lambda^-$ under the action of $\mathcal{O}^{-1}$.

The same statement holds for $\psi(\lambda^-)$. Thus, each $\psi \in \mathcal{H}_0$ satisfies one of the following two conditions.

1.  The four laminations $\psi(\Lambda^+), \psi(\Lambda^-), \Lambda^+$ and $\Lambda^-$ are not all distinct.

2.  $\psi(\lambda^+)$ and $\psi(\lambda^-)$ are weakly attracted to $\Lambda^+$ under the action of $\mathcal{O}$ and to $\Lambda^-$ under the action of $\mathcal{O}^{-1}$

Condition 1 either holds for both $\psi$ and $\psi^{-1}$ or fails for both $\psi$ and $\psi^{-1}$. Moreover, if condition 2 holds for both $\psi$ and $\psi^{-1}$, then $\psi$ satisfies the hypotheses of Corollary 3.4.3. Thus either condition 1 holds for all $\psi \in \mathcal{H}_0$, or the hypotheses of Corollary 3.4.3 are satisfied for some $\psi \in \mathcal{H}_0$. In the former case, Lemma 7.0.4 below implies that either $\mathrm{Stab}_{\mathcal{H}_0}(\Lambda^+)$ or $\mathrm{Stab}_{\mathcal{H}_0}(\Lambda^-)$ has finite index. In the latter case, Corollary 3.4.3 implies that $\mathcal{H}$ contains a free subgroup of rank two.   □

LEMMA 7.0.4.  *Suppose that a group $\mathcal{H}$ acts on a set $Y$ and that there are points $x, y \in Y$ such that*

$$(*) \qquad\qquad \{\psi(x), \psi(y)\} \cap \{x, y\} \neq \emptyset$$

*holds for all $\psi \in \mathcal{H}$. Then either $\mathrm{Stab}_{\mathcal{H}}(x)$ or $\mathrm{Stab}_{\mathcal{H}}(y)$ has finite index.*

*Proof.* If $\mathrm{Stab}_{\mathcal{H}}(x)$ does not have finite index, then there exist $h_i \in \mathcal{H}$, $i \geq 1$, such that the $h_i(x)$'s are all distinct; we may assume without loss of generality that $h_i(x) \notin \{x, y\}$. Condition $(*)$ implies that each $h_i(y)$ is either $x$ or $y$. Passing to a subsequence, we may assume that either each $h_i(y) = x$ or each $h_i(y) = y$. In the former case, each $h_1 h_i(y) \notin \{x, y\}$ and there are at most two values of $i$ for which $h_1 h_i(x) \in \{x, y\}$. This contradicts $(*)$ and we conclude that each $h_i(y) = y$.

By a completely symmetric argument, we conclude that if $\mathrm{Stab}_{\mathcal{H}}(y)$ does not have finite index, then there exist $g_j \in \mathcal{H}$, $j \geq 1$, such that the $g_j(y)$'s are distinct elements of $X \setminus \{x, y\}$ and such that each $g_j(x) = x$.

But then each $g_1 h_i(y) = g_1(y) \notin \{x, y\}$ and there are at most two values of $i$ for which $g_1 h_i(x) \in \{x, y\}$. This contradicts $(*)$ and completes the proof.  $\square$

The proof of Lemma 7.0.3 divides into the geometric and nongeometric cases. We consider the nongeometric case first, using the fact (Lemma 5.1.7) that $\langle Z, \hat{\rho}_r \rangle$ carries the same bi-infinite paths as a free factor system. The proof is particularly simple when $\mathrm{Rank}(H_1(G_r)) > \mathrm{Rank}(H_1(G_{r-1})) + 1$ and the reader may wish to focus on this case first. For the general case, we pass to finite covers via Lemma 7.0.5 below (which holds trivially with $k = 1$ in the case that $\mathrm{Rank}(H_1(G_r)) > \mathrm{Rank}(H_1(G_{r-1})) + 1$).

If $\tilde{G}$ is a cover of $G$ and $X$ is a subgraph of $G$, then we denote the full pre-image of $X$ by $\tilde{X}$. We denote the Euler characteristic by $\chi$.

LEMMA 7.0.5.  *Suppose that $f : G \to G$ is an improved relative train track map, that $H_r$ is an exponentially-growing stratum and that $G_r$ is connected. Then there exist $k > 0$ and a regular connected $k$-fold cover $\tilde{G}$ such that $\mathrm{Rank}(H_1(\tilde{G}_r)) > \mathrm{Rank}(H_1(\tilde{G}_{r-1})) + k$. Moreover, we can arrange that every outer automorphism of $F_n \cong \pi_1(G)$ lifts to an outer automorphism of $\pi_1(\tilde{G})$.*

*Proof.* If $G_{r-1}$ has contractible components then $H_{r-1}$ is a zero stratum and is the union of the contractible components of $G_{r-1}$. In that case, redefine the filtration by declaring each edge of $H_{r-1}$ to be an edge of $H_r$. This may destroy the relative train track property but $f : G \to G$ is still a topological representative with respect to this new shortened filtration. We may therefore assume that $f : G \to G$ is a topological representative and that each component of $G_{r-1}$ is noncontractible.

We may also assume that $\text{Rank}(H_1(G_r)) \leq \text{Rank}(H_1(G_{r-1})) + 1$ for otherwise the lemma is trivially satisfied with $k = 1$ and $\tilde{G} = G$. Let $m = \text{Rank}(H_0(G_{r-1}))$ be the number of components in $G_{r-1}$. Corollary 3.2.2 implies that either $m \geq 3$ or $\text{Rank}(H_1(G_r)) = \text{Rank}(H_1(G_{r-1})) + 1$ and $m = 2$.

Choose $k_1 > m$ and connected $k_1$-fold covering spaces for each component of $G_{r-1}$. Extend this to a $k_1$-fold covering space $\hat{G}$ of $G$. The key point here is that, independently of $k_1$, the full pre-image $\hat{G}_{r-1}$ of $G_{r-1}$ has $m$ components.

Since $F_n$ has only finitely many subgroups of index $k_1$, the intersection $N$ of all such subgroups is a normal subgroup of finite index. Let $k$ be the index of $N$, let $\tilde{G}$ be the regular connected $k$-fold cover of $G$ corresponding to $N$, and let $k_2$ be the integer $\frac{k}{k_1}$. Then $\tilde{G}$ is a $k_2$-fold cover of $\hat{G}$ and $\tilde{G}_{r-1}$ has at most $k_2 m$ components. It is easy to see that for all $\mathcal{O} \in Out(F_n)$, $N$ is invariant under the action induced by $\mathcal{O}$ on normal subgroups. In particular, every outer automorphism of $\pi_1(G)$ lifts to an outer automorphism of $\pi_1(\tilde{G})$. We say that $N$ is *characteristic*.

To verify the conclusions of the lemma, note that

$$\text{Rank}(H_1(\tilde{G}_i)) = k \times \text{Rank}(H_1(G_i)) - k \times \text{Rank}(H_0(G_i)) + \text{Rank}(H_0(\tilde{G}_i)).$$

(This follows from $\chi = \text{Rank}(H_0) - \text{Rank}(H_1)$ and $\chi(\tilde{G}_i) = k \times \chi(G_i)$). Thus,

$$\begin{aligned}
\text{Rank}&(H_1(\tilde{G}_r)) - \text{Rank}(H_1(\tilde{G}_{r-1})) \\
= \;\; & k \times [\text{Rank}(H_1(G_r)) - \text{Rank}(H_1(G_{r-1}))] \\
& + k \times [\text{Rank}(H_0(G_{r-1})) - \text{Rank}(H_0(G_r))] \\
& + [\text{Rank}(H_0(\tilde{G}_r)) - \text{Rank}(H_0(\tilde{G}_{r-1}))] \\
\geq \;\; & k \times [\text{Rank}(H_1(G_r)) - \text{Rank}(H_1(G_{r-1}))] + k(m-1) + (1 - mk_2) \\
= \;\; & k \times [\text{Rank}(H_1(G_r)) - \text{Rank}(H_1(G_{r-1}))] + k[m - 1 + \frac{1}{k} - \frac{m}{k_1}].
\end{aligned}$$

Since $m \geq 2$ and $k_1 > m$, both terms are nonnegative. If $m \geq 3$, then the second term is strictly larger than $k$. If $m = 2$, then the first term equals $k$ and the second term is positive. $\qquad\square$

*Proof of Lemma 7.0.3 in the nongeometric case.* Since $\{F(\Lambda^+)\} = \{F(\Lambda^-)\} = \mathcal{F}(G_r)$, $H_r$ is exponentially growing and $G_r$ is connected. Let $\tilde{G}$ be a $k$-fold cover of $G$ as in Lemma 7.0.5; if $\hat{\rho}_r = \rho_r$, let $\tilde{\rho}_r^1, \ldots, \tilde{\rho}_r^k \subset \tilde{G}$ be the lifts of $\rho_r$. Define $\mathcal{H}_0 \subset \mathcal{H}$ to be the finite index subgroup of elements whose lifts to $\pi_1(\tilde{G})$ act by the identity on $H_1(\tilde{G}; \mathbb{Z}_2)$. Suppose $\beta$ is a circuit in $\langle G_{r-1}, \hat{\rho}_r \rangle$ that lifts to a circuit $\tilde{\beta}$ in $\tilde{G}$. If $\hat{\rho}_r$ is trivial, then $\tilde{\beta} \subset \tilde{G}_{r-1}$; otherwise, $\tilde{\beta}$ can be decomposed into subpaths that are either single edges in $\tilde{G}_{r-1}$, some $\tilde{\rho}_r^i$ or the inverse of some $\tilde{\rho}_r^i$. In either case, the $\mathbb{Z}_2$-homology classes generated by all such $\tilde{\beta}$ are contained in a subspace of $H_1(\tilde{G}_r; \mathbb{Z}_2)$ of dimension at most $\text{Rank}(H_1(\tilde{G}_{r-1})) + k$. Lemma 7.0.5 implies that there is a circuit $\tilde{\alpha} \subset \hat{G}_r$ whose $\mathbb{Z}_2$-homology class is not represented by a lift of a circuit in $\langle G_{r-1}, \hat{\rho}_r \rangle$. Since $\tilde{Z} \cap \tilde{G}_r = \tilde{G}_{r-1}$, the $\mathbb{Z}_2$-homology class of $\tilde{\alpha}$ is not represented by a lift of a circuit in $\langle Z, \hat{\rho}_r \rangle$. Let $\alpha \subset G_r$ be the projected image of $\tilde{\alpha}$. For each $\psi \in \mathcal{H}_0$, the $\mathbb{Z}_2$-homology class determined by $\tilde{\psi}_\#(\tilde{\alpha})$ cannot be represented by the lift of a circuit in $\langle Z, \hat{\rho}_r \rangle$, and so $\psi_\#(\alpha) \notin \langle Z, \hat{\rho}_r \rangle$. Since $\psi_\#(\mathcal{F}(G_r)) = \{F(\psi_\#(\Lambda^+))\} = \{F(\psi_\#(\Lambda^-))\}$, every free factor system that contains $\psi_\#(\lambda^+)$ or $\psi_\#(\lambda^-)$ must contain $\psi_\#(\alpha)$ for every circuit $\alpha \subset G_r$. Lemma 5.1.7 therefore implies that $\psi_\#(\lambda^+) \notin \langle Z, \hat{\rho}_r \rangle$ and $\psi_\#(\lambda^-) \notin \langle Z, \hat{\rho}_r \rangle$. $\square$

We now turn to the proof of Lemma 7.0.3 in the case $H_r$ is a geometric stratum. The main difference between the cases is that we can no longer use Lemma 5.1.7 to conclude that if $\psi_\#(\lambda^+) \in \langle Z, \hat{\rho}_r \rangle$ or if $\psi_\#(\lambda^-) \in \langle Z, \hat{\rho}_r \rangle$ then $\psi_\#(\alpha) \in \langle Z, \hat{\rho}_r \rangle$ for every circuit $\alpha \subset G_r$. We replace this with Corollary 7.0.8 below.

Suppose that $\mu \subset G$ is a bi-infinite path and that $\alpha \subset G$ is a circuit. Choose a lift $\tilde{\alpha}$ of $\alpha$ in the universal cover $\Gamma$ of $G$ and let $T : \Gamma \to \Gamma$ be the indivisible covering translation with axis equal to $\tilde{\alpha}$. We say that $\alpha$ *is in the span of* $\mu$ if for all positive integers $L$, there are lifts $\tilde{\mu}_i$ of $\mu$, $0 \le i \le m-1$, such that

(Sp) $\tilde{\mu}_0 \cap \tilde{\mu}_1, \tilde{\mu}_1 \cap \tilde{\mu}_2, \ldots, \tilde{\mu}_{m-2} \cap \tilde{\mu}_{m-1}$ and $\tilde{\mu}_{m-1} \cap T(\tilde{\mu}_0)$ each contain at least $L$ edges.

LEMMA 7.0.6. *Suppose that* $f : G \to G$ *,* $Z$ *and* $\hat{\rho}_r$ *are as in Theorem 6.0.1. If* $\mu \in \langle Z, \hat{\rho}_r \rangle$ *and* $\alpha \subset G$ *is in the span of* $\mu$, *then* $\alpha \in \langle Z, \hat{\rho}_r \rangle$.

*Proof.* The lemma is obvious if $\mu \in Z$ so we may assume that $\hat{\rho}_r = \rho_r$. Write $\tilde{\mu} = \ldots \tilde{b}_{-1} \tilde{b}_0 \tilde{b}_1 \ldots$ where each $b_j$ is either a single edge of $Z$ or is equal to $\rho_r$ or $\bar{\rho}_r$. We do not assume *a priori* that this is a decomposition into subpaths, but there is no loss in assuming that $b_j \ne \bar{b}_j$. Since the initial edges of $\rho_r$ and of $\bar{\rho}_r$ lie in $H_r$ and are distinct, no cancellation can occur at the juncture of $\tilde{b}_j$ and $\tilde{b}_{j+1}$; we conclude that the $\tilde{b}_j$'s are subpaths of $\tilde{\mu}$.

Let $M$ be the number of edges in $\rho_r$. Choose $\tilde{\mu}_i$ satisfying (Sp) with $L = M$; denote $T(\tilde{\mu}_0)$ by $\tilde{\mu}_m$. Since $\tilde{\mu}_{i-1} \cap \tilde{\mu}_i$, $1 \le i \le m$, contains at least $M$ edges, there is a vertex $p_i \in \tilde{\mu}_{i-1} \cap \tilde{\mu}_i$ that is a cutting vertex (see the proof of Lemma 6.0.5) for both $\tilde{\mu}_{i-1}$ and $\tilde{\mu}_i$. For $1 \le i \le m-1$, let $\tilde{\gamma}_i$ be the subpath of $\tilde{\mu}_i$ that is bounded by $\tilde{p}_i$ and $\tilde{p}_{i+1}$. Let $\tilde{\gamma}_m$ be the subpath of $\tilde{\mu}_m$ that is bounded by $\tilde{p}_m$ and $T(\tilde{p}_1)$. Then each $\gamma_i \subset \langle Z, \rho_r \rangle$ and $\alpha$ is the circuit obtained by tightening $\gamma_1 \cdot \ldots \cdot \gamma_m$. $\square$

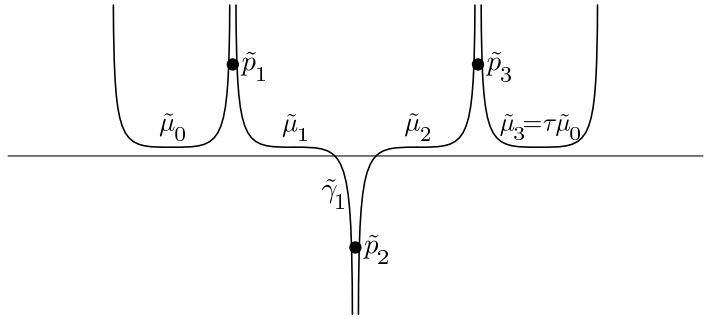

Let $\phi : S \to S$ and $\Phi : Y \to G_r$ be as in the definition of geometric stratum. For each closed geodesic curve $\alpha_S \subset S$, we say that $\Phi_\#(\alpha_S) \subset G_r$ is an $H_r$-*geometric circuit*.

LEMMA 7.0.7. *If $H_r$ is a geometric stratum with generic lines $\lambda^\pm$, then every $H_r$-geometric circuit in $G_r$ is in the span of $\lambda^+$ and in the span of $\lambda^-$.*

*Proof.* It suffices to consider $\lambda^+$. Suppose that an $H_r$-geometric circuit $\alpha$ and length $L$ are given. Let $\lambda_S^+$ and $\alpha_S$ be the geodesics is $S$ such that $\Phi_\#(\lambda_S^+) = \lambda^+$ and $\Phi_\#(\alpha_S) = \alpha^+$. We showed during the proof of Proposition 6.0.8 that $\lambda_S^+$ is a leaf of the expanding lamination $\Lambda_S^+$ for $\phi$. There exists $\varepsilon > 0$ and a length $L_S$ so that if a pair of lifts of $\lambda_S^+$ to the universal cover of $S$ contain $\varepsilon$-parallel subintervals of length $L_S$, then their $\Phi_\#$-images in the universal cover $\Gamma$ of $G$ have a common subinterval containing at least $L$ edges. There exists $\delta > 0$ so that any two lifts of $\lambda_S^+$ that have points within $\delta$ of each other, have $\varepsilon$-parallel subintervals of length $L_S$. The lemma now follows from the well-known fact that $\alpha_S$ is freely homotopic to a closed curve of the form $u_0 \cdot s_0 \cdot u_1 \cdot s_1 \cdot \ldots \cdot u_{k-1} \cdot s_{k-1}$ where each $u_i$ is an interval in $\lambda_S^+$ and where each $s_i$ has length at most $\delta$. (To prove this well-known fact, note that for sufficiently small $\delta$, the complement $S \setminus N_\delta(\Lambda_S^+)$ of the $\delta$ neighborhood of the expanding lamination $\Lambda_S^+$ in $S$ is a finite disjoint union of contractible or peripheral sets. The circuit $\alpha_S$ can therefore be homotoped into $N_\delta(\Lambda_S^+)$ and then further homotoped to have the desired decomposition into subpaths.) $\square$

COROLLARY 7.0.8. *Suppose that $f : G \to G$ , $Z$ and $\hat\rho_r$ are as in Theorem 6.0.1, that $H_r$ is a geometric stratum with generic lines $\lambda^\pm$ and that $\psi \in Out(F_n)$. If $\psi_\#(\lambda^+) \in \langle Z, \rho_r \rangle$ or $\psi_\#(\lambda^-) \in \langle Z, \rho_r \rangle$, then $\psi_\#(\alpha) \in \langle Z, \rho_r \rangle$ for each $H_r$-geometric circuit $\alpha \subset G_r$.*

*Proof.* Let $\tilde h : \Gamma \to \Gamma$ be a lift of a topological representative $h : G \to G$ representing $\psi$. For all $L$, there exists $L_1$ so that if $\tilde\beta \subset \Gamma$ is a path with edge length at least $L_1$, then $\tilde h_\#(\tilde\beta)$ is a path with edge length at least $L$. If $C$ is the constant of the bounded cancellation lemma applied to $\tilde h : \Gamma \to \Gamma$, and if $\tilde{\mu}_1$ and $\tilde{\mu}_2$ are bi-infinite paths such that $\tilde{\mu}_1 \cap \tilde{\mu}_2$ has edge length at least $L_1$, then $\tilde{h}_\#(\tilde{\mu}_1) \cap \tilde{h}_\#(\tilde{\mu}_2)$ has edge length at least $L - 2C$.

By Lemma 7.0.7, $\alpha$ is in the span of $\lambda^+$ and in the span of $\lambda^-$. The preceding argument shows that $\psi_\#(\alpha)$ is in the span of $\psi_\#(\lambda^+)$ and in the span of $\psi_\#(\lambda^-)$. Lemma 7.0.6 now completes the proof.    □

The following lemma is a modification of Lemma 7.0.5. We use the notation of Definition 5.1.4. Let $\partial_A S$ be the union of the components $\alpha_1^*, \ldots, \alpha_m^*$ of $\partial S$. Every regular cover $\tilde{G}$ of $G$ determines a regular cover $\tilde{S}$ of $S$; we denote the full pre-image of $\partial_A S$ by $\partial_A \tilde{S}$.

LEMMA 7.0.9.   *If $H_r$ is a geometric stratum, then there is a regular connected $k$-fold cover $\tilde{G}$ such that the induced cover $\tilde{S}$ of $S$ satisfies*

$$\mathrm{Rank}(H_1(\tilde{S})) - \mathrm{Rank}(H_1(\partial_A \tilde{S})) > k.$$

*Moreover, every outer automorphism of $F_n \cong \pi_1(G)$ lifts to an outer automorphism of $\pi_1(\tilde{G})$.*

*Proof.* If $S$ is $S^2$ with $m + 1$ disks removed, then

$$\mathrm{Rank}(H_1(S)) - \mathrm{Rank}(H_1(\partial_A S)) = 0.$$

If $S$ is a Mobius band with $m$ disks removed, then

$$\mathrm{Rank}(H_1(S)) - \mathrm{Rank}(H_1(\partial_A S)) = 1.$$

In all other cases, $\mathrm{Rank}(H_1(S)) - \mathrm{Rank}(H_1(\partial_A S)) \geq 2$ and we may choose $k = 1$ and $\tilde{G} = G$. Since $S$ supports a pseudo-Anosov homeomorphism, $m \geq 2$ and $m = 2$ only if $S$ is a Mobius band with 2 disks removed. (This last fact is well-known; it follows from Lemma 3.2.2 and the fact that if $S$ is $S^2$ with three disks removed or if $S$ is the Mobius band with one disk removed, then $S$ deformation retracts to a one complex made up of $\partial_A S$ and one edge.)

Choose elements $c_i \in \pi_1(G)$, $1 \leq i \leq m$, whose associated circuit is $\alpha_i$. For $j \geq 1$, denote the concatenation of $j$ copies of $c_i$ by $c_i^j$. Since $F_n$ is residually finite, there is a finite index normal subgroup $N$ that does not contain $c_i^j$ for $1 \leq i \leq m$ and $1 \leq j \leq 4$, and therefore does not contain any element conjugate to such $c_i^j$. If $\tilde{G}$ is the regular connected finite cover of $G$ corresponding to $N$ and if the closed path that goes $k$ times around $\alpha_i$ lifts to a closed circuit in $\tilde{G}$, then $k \geq 5$. Thus, if $q : \tilde{S} \to S$ is the covering space of $S$ induced by $\tilde{G}$, then the restriction of $q$ to any component of $\partial_A \tilde{S}$ is at least a five fold cover. It follows that $\mathrm{Rank}(H_0(\partial_A \tilde{S})) \leq \frac{k}{5} \mathrm{Rank}(H_0(\partial_A S))$. As noted in the proof of Lemma 7.0.5, after passing to a further cover if necessary, we may assume that $N$ is characteristic.

Let $m = \text{Rank}(H_0(\partial_A S))$ be the number of components of $\partial_A S$; since $S$ is connected, $\text{Rank} H_0(S) = 1$. As in the proof of Lemma 7.0.5,

$$\text{Rank}(H_1(\tilde{S})) - \text{Rank}(H_1(\partial_A \tilde{S}))$$

$$\begin{aligned}
= \quad & k \times [\text{Rank}(H_1(S)) - \text{Rank}(H_1(\partial_A S))] \\
& + \; k \times [\text{Rank}(H_0(\partial_A S)) - \text{Rank}(H_0(S))] \\
& + \; [\text{Rank}(H_0(\tilde{S})) - \text{Rank}(H_0(\partial_A \tilde{S}))] \\
\geq \quad & k \times [\text{Rank}(H_1(S)) - \text{Rank}(H_1(\partial_A S))] + k(m-1) + (1 - \frac{mk}{5}) \\
= \quad & k \times [\text{Rank}(H_1(S)) - \text{Rank}(H_1(\partial_A S))] + k[m - 1 + \frac{1}{k} - \frac{m}{5}].
\end{aligned}$$

Both terms are nonnegative. If $m \geq 3$, then the second term is strictly larger than $k$. If $m = 2$, then the first term equals $k$ and the second term is positive. $\qquad \square$

*Proof of Lemma* 7.0.3 *in the geometric case.* The proof is now essentially the same as in the nongeometric case. Let $\tilde{G}$ be a covering space of $G$ as in Lemma 7.0.9 and let $\mathcal{H}_0$ be the finite index subgroup of $\mathcal{H}$ whose lifts to $\pi_1(\tilde{G})$ act by the identity on $H_1(\tilde{G}; \mathbb{Z}_2)$. Since $\partial S = \partial_A S \cup \rho^*$, there are at most $k$ components in $\partial \tilde{S} \setminus \partial_A \tilde{S}$. It is therefore possible to choose a circuit $\tilde{\alpha}_S$ in $\tilde{S}$ whose $\mathbb{Z}_2$-homology class is not represented by peripheral curves. There are induced covering spaces $\tilde{Y}, \tilde{Q}$ and $\tilde{S}$ and there is an induced homotopy equivalence $\tilde{\Phi} : Y \to \tilde{G}$. The peripheral homology $H_1(\partial_A \tilde{S}; \mathbb{Z}_2)$ is a direct summand of the homology $H_1(\tilde{S}; \mathbb{Z}_2)$; $H_1(\tilde{Y}; \mathbb{Z}_2)$ is formed from $H_1(\tilde{S}; \mathbb{Z}_2) \bigoplus H_1(\tilde{Q}; \mathbb{Z}_2)$ by identification of $H_1(\partial_A \tilde{S}; \mathbb{Z}_2)$ with its image in $H_1(\tilde{Q}; \mathbb{Z}_2)$. It follows that $\tilde{\alpha} = \tilde{\Phi}_\#(\tilde{\alpha}_S)$ determines a nonzero element of $H_1(\tilde{G}_r; \mathbb{Z}_2)$ that is not represented by a lift of a circuit in $\langle G_{r-1}, \rho_r \rangle$ and hence is not represented by a lift of a circuit in $\langle Z, \rho_r \rangle$. If $\alpha$ is the projection of $\tilde{\alpha}$, then $\psi_\#(\alpha)$ is not in $\langle Z, \rho_r \rangle$ and Corollary 7.0.8 completes the proof. $\qquad \square$

LEMMA 7.0.10. *Assume that $\mathcal{H} \subset \text{Out}(F_n)$ does not contain a free subgroup of rank two. Then there are a finite collection $\mathcal{L}$ of attracting laminations for elements of $\mathcal{H}$ and a finite index subgroup $\mathcal{H}_0$ of $\mathcal{H}$ that stabilizes each element of $\mathcal{L}$ with the following feature. If $\psi \in \mathcal{H}_0$ and if $\Lambda^+ \in \mathcal{L}(\psi)$ and $\Lambda^- \in \mathcal{L}(\psi^{-1})$ are paired topmost laminations, then at least one of $\Lambda^+$ and $\Lambda^-$ is in $\mathcal{L}$.*

*Proof.* Among all free factor systems other than $\{[[F_n]]\}$ that are invariant under the action of a finite index subgroup of $\mathcal{H}$, choose one, $\mathcal{F}_1$, of maximal complexity (as defined in subsection 2.6). If $\mathcal{H}_0$ is a finite index subgroup of $\mathcal{H}$ that stabilizes $\mathcal{F}_1$, then, by induction applied to $\mathcal{H}_0|\mathcal{F}_1$, we may assume that there is a finite collection $\mathcal{L}_1$ of attracting laminations carried by $\mathcal{F}_1$ and a finite index subgroup (also called $\mathcal{H}_0$) that stabilizes each element of $\mathcal{L}_1$ with the following feature. If $\psi \in \mathcal{H}_0$ and if $\Lambda^+ \in \mathcal{L}(\psi)$ and $\Lambda^- \in \mathcal{L}(\psi^{-1})$ are paired topmost laminations carried by $\mathcal{F}_1$, then either $\Lambda^+ \in \mathcal{L}_1$ or $\Lambda^- \in \mathcal{L}_1$.

If every topmost lamination pair for elements of $\mathcal{H}_0$ is carried by $\mathcal{F}_1$, then we are done. Otherwise, choose $\mathcal{O} \in \mathcal{H}_0$ and a topmost lamination pair $\Gamma^+ \in \mathcal{L}(\mathcal{O})$ and $\Gamma^- \in \mathcal{L}(\mathcal{O})$ that is not carried by $\mathcal{F}_1$. Proposition 7.0.2 implies, after passing to a smaller finite index subgroup if necessary, that at least one of $\Gamma^+$ or $\Gamma^-$ is stabilized by $\mathcal{H}_0$. We assume without loss of generality that $\Gamma^+$ is stabilized by $\mathcal{H}_0$; if possible, choose $\mathcal{H}_0$ to stabilize both $\Gamma^+$ and $\Gamma^-$. Define $\mathcal{L}$ to be the union of $\mathcal{L}_1$ with $\Gamma^+$ and with $\Gamma^-$ if it is stabilized by $\mathcal{H}_0$.

Choose an improved relative train track map $f : G \to G$ for some iterate of $\mathcal{O}$ such that $\mathcal{F}_1 = \mathcal{F}(G_l)$ for some filtration element $G_l$. Since $\mathcal{F}_1$ and $\Gamma^+$ are both $\mathcal{H}_0$-invariant, so is the unique smallest free factor system containing $\mathcal{F}_1$ and carrying $\Gamma^+$. Our choice of $\mathcal{F}_1$ therefore guarantees that this smallest free factor system is $\{[[F_n]]\}$. It follows that $\Gamma^+$ is associated to the highest stratum $G_r = G$.

Suppose that $\psi \in \mathcal{H}_0$ and that $\Lambda^+ \in \mathcal{L}(\psi)$ and $\Lambda^- \in \mathcal{L}(\psi^{-1})$ are paired topmost laminations. If $\Lambda^{\pm}$ are carried by $\mathcal{F}_1$, then either $\Lambda^+ \in \mathcal{L}_1 \subset \mathcal{L}$ or $\Lambda^- \in \mathcal{L}_1 \subset \mathcal{L}$. Suppose then that $\Lambda^{\pm}$ are not carried by $\mathcal{F}_1$. Proposition 7.0.2 implies that either $\Lambda^+$ or $\Lambda^-$, say $\Lambda^-$, is stabilized by a finite index subgroup of $\mathcal{H}_0$. After replacing $f$ by an iterate if necessary, $f_{\#}$ stabilizes $\Lambda^-$.

Let $\lambda^-$ be a generic line for $\Lambda^-$. Since $\{[[F_n]]\}$ is the only free factor system that contains $\mathcal{F}_1$ and carries $\Lambda^-$, we have $\lambda^- \not\subset G_{r-1}$. Theorem 6.0.1 and Remark 6.0.2 therefore imply that either $\lambda^-$ is $\Gamma^+$-generic or $\lambda^-$ is weakly attracted to $\Gamma^+$ under the action of $f_{\#}$. In the former case, $\Lambda^- = \Gamma^- \in \mathcal{L}$. In the latter case, every subpath of a generic line $\gamma^+$ of $\Gamma^+$ is contained in some line $f_{\#}^m(\lambda^-)$ of $\Lambda^-$, so that $\gamma^+$ is a line in $\Lambda^-$. By our previous arguments, $\psi$ is represented by an improved relative train track map in which $\Lambda^-$ is associated to the highest stratum and in which the next highest stratum realizes $\mathcal{F}_1$. Since $\gamma^+$ is not carried by $\mathcal{F}_1$, Lemma 3.1.15 implies that $\gamma^+$ is a generic line for $\Lambda^-$ and hence that $\Gamma^+ = \Lambda^- \in \mathcal{L}$.        $\square$

*Proof of Theorem* 7.0.1. Let $\mathcal{L} = \{\Lambda_1, \dots \Lambda_k\}$ and $\mathcal{H}_0$ be as in Lemma 7.0.10. Define $\Phi = \bigoplus PF_{\Lambda_i^+} : \mathcal{H}_0 \to \mathbb{Z}^k$ where each $PF_{\Lambda_i^+}$ is as in Corollary 3.3.1. By Corollary 5.7.6, it suffices to show that $\mathrm{Ker}(\Phi)$ is contained in $PG(F_n)$. If $\psi \in \mathcal{H}_0$ is not in $PG(F_n)$, then there exist paired, topmost laminations $\Gamma^+ \in \mathcal{L}(\psi)$ and $\Gamma^- \in \mathcal{L}(\psi)$. Proposition 3.3.3 implies that neither $PF_{\Gamma^+}(\psi)$ nor $PF_{\Gamma^-}(\psi)$ is zero so that Lemma 7.0.10 implies that $\psi \notin \mathrm{Ker}(\Phi)$.        $\square$

University of Utah, Salt Lake City, UT
*E-mail address*: bestvina@math.utah.edu

Rutgers University, Newark, NJ
*E-mail address*: feighn@andromeda.rutgers.edu

Herbert H. Lehman College (CUNY)
*E-mail address*: michael@alpha.lehman.cuny.edu